\documentclass{amsart}
\usepackage{amsfonts,amsthm,amsmath,amssymb}
\usepackage{enumerate}
\usepackage{booktabs}
\usepackage[T1]{fontenc}
\usepackage[utf8]{inputenc}
\usepackage{xcolor}
\usepackage[hidelinks,pdfencoding=auto,psdextra]{hyperref}
\usepackage{microtype,placeins,etoolbox}
\AtBeginEnvironment{tikzpicture}{\microtypesetup{expansion=false}}
\usepackage{pgfplots}
\pgfplotsset{compat=1.18}

\newtheorem{theorem}{Theorem}[section]
\newtheorem{proposition}[theorem]{Proposition}
\newtheorem{lemma}[theorem]{Lemma}
\newtheorem{corollary}[theorem]{Corollary}
\newtheorem{conjecture}[theorem]{Conjecture}
\newtheorem*{maintheorem}{Main Theorem}
\theoremstyle{definition}
\newtheorem{definition}[theorem]{Definition}

\theoremstyle{remark}
\newtheorem{remark}[theorem]{Remark}
\newtheorem*{addendum}{Addendum}

\newcommand{\C}{\mathbb{C}}
\newcommand{\R}{\mathbb{R}}
\newcommand{\Z}{\mathbb{Z}}
\newcommand{\N}{\mathbb{N}}
\newcommand{\ii}{{\rm i}}
\renewcommand{\Re}{\operatorname{Re}}
\renewcommand{\Im}{\operatorname{Im}}
\newcommand{\supp}{\operatorname{supp}}
\newcommand{\Pot}[1]{\mathcal{L}_{#1}}
\newcommand{\M}{\mathcal{M}}
\newcommand{\Sweep}{\Sigma}
\newcommand{\W}{\mathcal{W}}

\numberwithin{equation}{section}
\definecolor{actionblue}{RGB}{62,128,181}
\definecolor{actionred}{RGB}{196,80,64}
\definecolor{actiongreen}{RGB}{64,140,108}
\definecolor{actiongold}{RGB}{211,151,49}
\pgfplotsset{paperplot/.style={scale only axis,axis equal image,
  xlabel={$\Re E$},ylabel={$\Im E$},label style={font=\footnotesize},
  tick label style={font=\scriptsize},title style={font=\small},
  tick align=outside,axis line style={gray!75},
  scaled ticks=false,/pgf/number format/1000 sep={}}}

\definecolor{spectralblue}{RGB}{25,113,166}
\definecolor{frozenpink}{RGB}{225,110,179}
\pgfplotsset{introplot/.style={scale only axis,axis equal image,
  tick label style={font=\scriptsize},title style={font=\small},
  axis line style={black!45},tick style={black!45},
  tick align=outside,major tick length=2pt,scaled ticks=false,
  /pgf/number format/1000 sep={},xtick pos=bottom,ytick pos=left}}

\hypersetup{pdftitle={The enigma of locally Toeplitz matrices. I. Discrete Schrödinger operators with complex potential},
  pdfauthor={Boris Shapiro}}

\begin{document}

\title[The enigma of locally Toeplitz matrices. I.]
{The enigma of locally Toeplitz matrices. I.\\ Discrete Schr\"odinger operators with complex potential}

\author[B.~Shapiro]{Boris Shapiro}
\address[B.~Shapiro]{Department of Mathematics, Stockholm University, SE-106 91 Stockholm, Sweden}
\email{shapiro@math.su.se}
\address[B.~Shapiro]{Department of Mathematics with Computer Science,\newline Guangdong Technion--Israel
Institute of Technology, 241 Daxue Road, Shantou 515063, P.R. China}
\email{boris.shapiro@gtiit.edu.cn}

\date{November 4, 2026; revised November 25, 2026}

\begin{abstract}
A banded matrix whose diagonals vary slowly is locally close to a
Toeplitz matrix. Freezing its coefficients at each position gives a
family of limiting Toeplitz eigenvalue distributions. Their average,
which we call the standard measure, typically fills a planar domain,
whereas the eigenvalues of a non-normal matrix of this kind appear to
accumulate on finitely many curves inside that domain. We conjecture
that, for generic symbols, the eigenvalue distributions converge to a
mother body of the standard measure: a measure on finitely many analytic
arcs with the same exterior logarithmic potential and with a potential
bounded above by that of the standard measure. For every symbol we prove
that each subsequential limit has the last two properties, and we show
by an example that without genericity the eigenvalue distributions need
not converge.

We then prove the conjecture for finite discrete Schr\"odinger operators
with a complex potential $V$, Dirichlet boundary conditions and a
constant hopping coefficient $\alpha>0$. When $\Im V$ is strictly
monotone and $V$ admits a suitable holomorphic extension, the eigenvalue
distributions converge to a measure on finitely many analytic arcs.
Its logarithmic potential falls below that of the standard measure by
twice the smallest of four explicit tunnelling actions. The proof rests
on a four-term asymptotic expansion of the Dirichlet determinant,
obtained by the complex WKB method for difference equations. When
$\Im V$ is not monotone, an exact transfer recursion gives a conditional
description of the limit. For a circular potential arc we verify the
geometric switching hypothesis of this description, leaving
non-resonance as the only remaining assumption.
\end{abstract}


\maketitle

\section{Locally Toeplitz matrices and the mother-body conjecture}\label{sec:intro}

Let $a(x,z)=\sum_{|j|\le d}a_j(x)z^j$ be a Laurent polynomial in $z$
whose coefficients are continuous functions of $x\in[0,1]$. We call $a$
a \emph{banded symbol}. The sequence of \emph{locally Toeplitz
matrices} associated with $a$ is
\begin{equation}\label{eq:lt-array}
 T_N(a)=\left[a_{j-i}\!\left(\frac{i+j}{2N}\right)\right]_{0\le i,j<N},
 \qquad N\ge1 .
\end{equation}
The diagonal with offset $j$ samples the coefficient $a_j$ along
$[0,1]$, so $T_N(a)$ is an ordinary Toeplitz matrix when $a$ does not
depend on $x$. Matrices of this type, also called twisted Toeplitz
matrices \cite{TC04,TE05}, arise as discretizations of differential
operators with variable coefficients. In the tridiagonal case they are
finite sections of Jacobi operators with slowly varying coefficients;
in particular, they include discrete Schr\"odinger operators with a
slowly varying potential.

A large banded matrix with slowly varying diagonals locally resembles a
Toeplitz matrix. Near the row with index $Nx$ its diagonals are
approximately constant, with values $a_j(x)$. Using these constants
throughout gives the Toeplitz matrices with symbol $a(x,\cdot)$. As their
size tends to infinity, their eigenvalue distributions converge to a
probability measure $\nu_x$, described by Schmidt--Spitzer and Hirschman
\cite{SS60,Hi67}; see also \cite{BG05,DK08}. We compare two probability
measures associated with $a$,
\[
 \mu_N(a)=\frac1N\sum_{j=1}^N\delta_{\lambda_j(T_N(a))},
 \qquad \M(a)=\int_0^1\nu_x\,dx ,
\]
where the eigenvalues $\lambda_j(T_N(a))$ are counted with algebraic
multiplicity. We call $\M(a)$ the \emph{standard measure} of $a$. It is
determined entirely by the family of frozen problems. If $a$ does not
depend on $x$, then $\mu_N(a)$ converges to $\M(a)$ by the definition
of $\nu_x$. For Hermitian locally Toeplitz matrices the same conclusion
is a spectral distribution theorem \cite{Kac53,Tilli}; for real Jacobi
matrices with varying recurrence coefficients see also \cite{KVA99}.

For non-normal matrices the picture is entirely different. The support
of each $\nu_x$ is a finite union of analytic arcs. As $x$ runs through
$[0,1]$, these supports typically sweep out a two-dimensional region of
the complex plane, and $\M(a)$ then has a density with respect to area
measure there. Yet the eigenvalues of $T_N(a)$ appear to follow just a
few curves inside this region. Figures~\ref{fig:clouds}
and~\ref{fig:standard} in Section~\ref{sec:examples} illustrate this
phenomenon for three examples. The basic question of this paper is
therefore the following.
\begin{quote}
\emph{How is the limiting eigenvalue distribution of $T_N(a)$ related to
the standard measure $\M(a)$ as $N\to\infty$?}
\end{quote}

\subsection{What the traces tell us}
Some information survives freezing even for non-normal matrices. For
every integer $m\ge1$,
\begin{equation}\label{eq:lt-moments}
 \lim_{N\to\infty}\frac1N\operatorname{tr}T_N(a)^m
 =\int_0^1 [z^0]\,a(x,z)^m\,dx
 =\int_{\C} w^m\,d\M(a)(w),
\end{equation}
where $[z^0]$ denotes the constant term of a Laurent polynomial. To see
the first equality, expand a diagonal entry of $T_N(a)^m$ as a sum over
paths of $m$ steps that return to their starting site. Such a path stays
within $md$ sites of its start. Away from the endpoints, replacing each
coefficient on the path by its value at the starting site introduces an
error $o(1)$, by uniform continuity. The sum of the frozen path weights
is exactly $[z^0]\,a(x,z)^m$. Since only $O(md)$ starting sites see an
endpoint, the normalized trace tends to the displayed integral. The same
calculation at fixed $x$ gives the second equality.

By the row-sum bound, all eigenvalues lie in the disk of radius
$\sum_j\|a_j\|_\infty$. Hence every subsequential limit of $\mu_N(a)$
has the same holomorphic moments $\int w^m$, $m\ge0$, as $\M(a)$, while
the traces say nothing about the mixed moments $\int w^m\bar w^k$.
Expanding $\log|E-w|$ in powers of $w/E$ shows that the logarithmic
potentials of the two measures agree for large $|E|$, and hence, by
harmonicity, on the unbounded component of the complement of both
supports. In the bounded components of that complement the traces give
no information. This still leaves room for very different distributions
of mass. For instance, a unit point mass at the origin and the uniform
probability measure on the unit disk both have exterior potential
$\log|E|$. The question is which distribution of mass the matrices
themselves select.

\subsection{The mother-body conjecture}
For a finite positive measure $\nu$ with compact support we write
\[
 \Pot{\nu}(E)=\int_{\C}\log|E-w|\,d\nu(w),
\]
so that $\Delta\Pot{\nu}=2\pi\nu$ in the sense of distributions.
Potentials take values in $[-\infty,\infty)$. The relevant notion comes
from potential theory and the inverse problem of geophysics, where a
mother body of a solid body is a measure on a lower-dimensional skeleton
that produces the same exterior field \cite{Gus98,BoSh16}. We use it
in the following form.

\begin{definition}\label{def:MB}
Let $U$ be a bounded open set with piecewise smooth boundary. Let $\Xi$
be a positive, absolutely continuous measure supported on $\overline U$
whose density is real-analytic away from finitely many real-analytic arcs.
A \emph{mother body} of $\Xi$ is a positive measure $\nu$ of the same
mass with the following properties.
\begin{enumerate}[(i)]
\item $\supp\nu\subseteq\overline U$ is a finite union of compact
real-analytic arcs and points.
\item $\Pot{\nu}=\Pot{\Xi}$ on $\C\setminus\overline U$.
\item $\Pot{\nu}\le\Pot{\Xi}$ throughout $\C$.
\end{enumerate}
\end{definition}

Thus the mass on the arcs produces the same exterior potential as the
planar mass, and condition (iii) fixes the direction of the change
inside. About a decade ago the author suggested that the eigenvalues of
locally Toeplitz matrices select a measure of this kind.

\begin{conjecture}[The locally Toeplitz mother-body conjecture]\label{conj:LT}
Let $a$ be a generic banded symbol with continuous, piecewise analytic
coefficients on $[0,1]$, and use the sampling in \eqref{eq:lt-array}.
\begin{enumerate}[(a)]
\item The measures $\mu_N(a)$ converge weakly to a probability measure
$\mu(a)$.
\item $\supp\mu(a)$ is a finite union of real-analytic arcs and isolated
points.
\item Suppose that $U=\operatorname{int}\supp\M(a)$ is bounded with
piecewise smooth boundary, that $\supp\M(a)\subseteq\overline U$, and that
the density of $\M(a)$ is real-analytic away from finitely many
real-analytic arcs, so that $\Xi=\M(a)$ satisfies the assumptions of
Definition~\ref{def:MB}. Then $\mu(a)$ is a mother body of $\M(a)$.
\end{enumerate}
\end{conjecture}

Genericity cannot be dropped in part~(a). By
Theorem~\ref{thm:LTfacts}(5) below, the measures $\mu_N(a)$ need not
converge even for a symbol with polynomial coefficients; in that example
$a_{\pm1}=0$, so that $T_N(a)$ splits into two blocks. On the other
hand, several assertions of the conjecture hold for every symbol. By
Theorem~\ref{thm:LTfacts}(3), part~(c) follows from part~(b); more
precisely, the conclusion of~(c) holds for every weak subsequential
limit of $\mu_N(a)$ whose support is as in~(b). If $\supp\M(a)$ has zero
planar area, then $\mu_N(a)\to\M(a)$ for every symbol, by
Theorem~\ref{thm:LTfacts}(4).

\begin{theorem}\label{thm:LTfacts}
Let $a$ be a banded symbol with continuous, piecewise analytic
coefficients on $[0,1]$, and use the sampling in \eqref{eq:lt-array}.
\begin{enumerate}[(1)]
\item All $\mu_N(a)$ are supported in a fixed compact set. Hence weak
subsequential limits exist, and they are probability measures.
\item Every weak subsequential limit $\nu$ of $\mu_N(a)$ satisfies
$\supp\nu\subseteq\supp\M(a)$, $\Pot{\nu}=\Pot{\M(a)}$ on
$\C\setminus\supp\M(a)$, and $\Pot{\nu}\le\Pot{\M(a)}$ on $\C$.
\item Let $U=\operatorname{int}\supp\M(a)$ and suppose that
$\supp\M(a)\subseteq\overline U$. If $\nu$ is a weak subsequential limit
of $\mu_N(a)$ and $\supp\nu$ is a finite union of compact real-analytic
arcs and points, then $\nu$ has properties (i)--(iii) of
Definition~\ref{def:MB} with $\Xi=\M(a)$.
\item If $\supp\M(a)$ has zero planar area, then $\mu_N(a)\to\M(a)$
weakly.
\item For $d=2$, $a_0(x)=2x-1$, $a_2(x)=1$, $a_{-2}(x)=-x(1-x)$ and
$a_{\pm1}=0$, the measures $\mu_N(a)$ do not converge. In this example
$\supp\M(a)$ has positive area.
\end{enumerate}
\end{theorem}

The proof of Theorem~\ref{thm:LTfacts} is given in
Section~\ref{sec:proofLT}. Parts (1)--(4) combine the
Schmidt--Spitzer--Hirschman description of the frozen measures $\nu_x$
with a transfer-matrix representation of $\det(E-T_N(a))$ and
elementary potential theory. Part~(5) is an explicit computation.
By parts (2) and (3), the conjecture reduces to the existence of the
limit and the structure of its support. The real difficulty lies inside
$\supp\M(a)$, where Theorem~\ref{thm:LTfacts} gives only an inequality,
and the size of the deficit $\Pot{\M(a)}-\Pot{\nu}$ has to be
determined.

\subsection{Three examples in the symmetric tridiagonal case}\label{sec:examples}
Apart from Theorem~\ref{thm:LTfacts}, which holds for all banded
symbols, the rest of the paper is devoted to the symmetric tridiagonal
case
\[
 a(x,z)=b(x)+c(x)\bigl(z+z^{-1}\bigr),
\]
in which we are able to obtain partial rigorous results. From
Section~\ref{sec:operator} on we specialize further to a constant
off-diagonal coefficient, that is, to discrete Schr\"odinger operators
with a complex potential.

For such a symbol the frozen problem is explicit. At position $x$ the
eigenvalues of the constant $m\times m$ tridiagonal matrix are
\[
 b(x)+2c(x)\cos\frac{j\pi}{m+1},\qquad 1\le j\le m.
\]
Thus $\nu_x$ is the arcsine law on the segment
$S_x=b(x)+c(x)[-2,2]$, interpreted as a point mass when $c(x)=0$, and
the union of these segments is the support of $\M(a)$.

In the figures we take $N=100$ and use a slightly different sampling:
the diagonal entry at site $k$ is $b(k/N)$, and both entries on the
edge from $k$ to $k+1$ are $c(k/N)$. This left-endpoint rule moves each
sampling point of \eqref{eq:lt-array} by at most $1/(2N)$ and has the
same frozen segments and the same standard measure. For constant $c$
the two rules give exactly the same matrices, so nothing changes from
Section~\ref{sec:operator} on. For non-constant $c$ the matrices differ,
and for non-normal matrices such an $O(1/N)$ change is not harmless a
priori; at $N=100$ the two finite spectra differ visibly. The figures
therefore illustrate the conjecture for the left-endpoint variant, and
we expect parts (a)--(c) of Conjecture~\ref{conj:LT} to be insensitive
to such $O(1/N)$ changes of the sampling points.

\begin{figure}[!htbp]
\centering
\begin{tabular}{@{}c@{\hspace{1.0em}}c@{}}
\begin{tikzpicture}[baseline=(current bounding box.center)]
\begin{axis}[introplot,width=4.2cm,height=3.5cm,xtick={-1,0,1},xmin=-1.08,xmax=1.08,ymin=-.50,ymax=.50,ytick={-.4,0,.4},title={(a)}]
\addplot[only marks,mark=*,mark size=.64pt,spectralblue] coordinates {
(-0.8023532,-0.2453297) (-0.8023532,0.2453297) (-0.6858159,-0.3031349) (-0.6858159,0.3031349) (-0.6012374,-0.3396733)
(-0.6012374,0.3396733) (-0.5322898,-0.3650686) (-0.5322898,0.3650686) (-0.4731118,-0.3833274) (-0.4731118,0.3833274)
(-0.4208171,-0.3965289) (-0.4208171,0.3965289) (-0.373735,-0.4059172) (-0.373735,0.4059172) (-0.3307965,-0.4123064)
(-0.3307965,0.4123064) (-0.2912689,-0.4162638) (-0.2912689,0.4162638) (-0.2546241,-0.4182029) (-0.2546241,0.4182029)
(-0.2204665,-0.4184365) (-0.2204665,0.4184365) (-0.188491,-0.4172077) (-0.188491,0.4172077) (-0.1584562,-0.4147096)
(-0.1584562,0.4147096) (-0.1301672,-0.4110991) (-0.1301672,0.4110991) (-0.1034643,-0.4065054) (-0.1034643,0.4065054)
(-0.0782144,-0.4010364) (-0.0782144,0.4010364) (-0.0543054,-0.3947837) (-0.0543054,0.3947837) (-0.0316417,-0.3878252)
(-0.0316417,0.3878252) (-0.0101413,-0.3802284) (-0.0101413,0.3802284) (0.0102668,-0.3720522) (0.0102668,0.3720522)
(0.0296448,-0.3633479) (0.0296448,0.3633479) (0.048047,-0.3541612) (0.048047,0.3541612) (0.0655216,-0.3445325)
(0.0655216,0.3445325) (0.0821112,-0.3344981) (0.0821112,0.3344981) (0.0978538,-0.3240905) (0.0978538,0.3240905)
(0.1127835,-0.3133394) (0.1127835,0.3133394) (0.1269306,-0.3022714) (0.1269306,0.3022714) (0.1403227,-0.2909109)
(0.1403227,0.2909109) (0.1529844,-0.2792806) (0.1529844,0.2792806) (0.1649382,-0.267401) (0.1649382,0.267401)
(0.1762042,-0.2552911) (0.1762042,0.2552911) (0.1868008,-0.2429689) (0.1868008,0.2429689) (0.1967446,-0.2304508)
(0.1967446,0.2304508) (0.2060505,-0.2177522) (0.2060505,0.2177522) (0.2147324,-0.2048878) (0.2147324,0.2048878)
(0.2228026,-0.1918711) (0.2228026,0.1918711) (0.2302722,-0.1787152) (0.2302722,0.1787152) (0.2371513,-0.1654324)
(0.2371513,0.1654324) (0.2434491,-0.1520343) (0.2434491,0.1520343) (0.2491736,-0.1385323) (0.2491736,0.1385323)
(0.2543322,-0.124937) (0.2543322,0.124937) (0.2589312,-0.1112588) (0.2589312,0.1112588) (0.2629764,-0.0975077)
(0.2629764,0.0975077) (0.2664725,-0.0836935) (0.2664725,0.0836935) (0.2694239,-0.0698256) (0.2694239,0.0698256)
(0.2718339,-0.0559133) (0.2718339,0.0559133) (0.2737055,-0.0419657) (0.2737055,0.0419657) (0.2750408,-0.0279918)
(0.2750408,0.0279918) (0.2758414,-0.0140003) (0.2758414,0.0140003) (0.2761081,0) (1,0)};
\end{axis}
\end{tikzpicture} &
\begin{tikzpicture}[baseline=(current bounding box.center)]
\begin{axis}[introplot,width=5.0cm,height=3.5cm,xmin=-11.7,xmax=-2.6,xtick={-10,-8,-6,-4},ymin=-1.45,ymax=1.45,ytick={-1,0,1},title={(b)}]
\addplot[only marks,mark=*,mark size=.64pt,spectralblue] coordinates {
(-11.1535211,-0.52988) (-11.1535211,0.52988) (-11.1461654,-0.3645794) (-11.1461654,0.3645794) (-11.1302954,-0.7146111)
(-11.1302954,0.7146111) (-11.1231487,-0.2201891) (-11.1231487,0.2201891) (-11.0904708,-0.0959711) (-11.0904708,0.0959711)
(-11.0619376,0) (-11.0197145,-0.9179413) (-11.0197145,0.9179413) (-10.9587306,0) (-10.8212504,0)
(-10.6967252,-1.2464771) (-10.6967252,1.2464771) (-10.6832277,0) (-10.5440232,0) (-10.4033085,0)
(-10.2609732,0) (-10.1170781,0) (-9.9718047,0) (-9.8254093,0) (-9.6781892,0)
(-9.530459,0) (-9.3825336,0) (-9.2347188,0) (-9.0873051,0) (-8.940564,0)
(-8.7947466,0) (-8.6500834,0) (-8.5067846,0) (-8.3650406,0) (-8.2250235,0)
(-8.0868883,0) (-7.9507738,0) (-7.8168045,0) (-7.6850914,0) (-7.5557335,0)
(-7.4288188,0) (-7.3044257,0) (-7.1826237,0) (-7.0634745,0) (-6.9470331,0)
(-6.8333485,0) (-6.7224642,0) (-6.6144194,0) (-6.5092495,0) (-6.4069864,0)
(-6.3076592,0) (-6.2112945,0) (-6.1179159,0) (-6.0275443,0) (-5.9401955,0)
(-5.8558783,0) (-5.7745878,0) (-5.6962944,0) (-5.6209212,0) (-5.5483003,0)
(-5.4780904,0) (-5.4096453,0) (-5.3418922,0) (-5.3192172,-0.173665) (-5.3192172,0.173665)
(-5.3012224,-0.0924914) (-5.3012224,0.0924914) (-5.2842862,-0.015443) (-5.2842862,0.015443) (-5.23833,0)
(-5.1847646,0) (-5.1303635,0) (-5.0742707,0) (-5.0164012,0) (-4.9567383,0)
(-4.8952982,0) (-4.8321143,0) (-4.7672271,0) (-4.70068,0) (-4.6325159,0)
(-4.5627762,0) (-4.4914999,0) (-4.4187238,0) (-4.344482,0) (-4.2688062,0)
(-4.191726,0) (-4.1132686,0) (-4.0334597,0) (-3.9523228,0) (-3.8698802,0)
(-3.7861524,0) (-3.7011589,0) (-3.6149176,0) (-3.5274456,0) (-3.4387588,0)
(-3.3488723,0) (-3.2578002,0) (-3.1655559,0) (-3.0721521,0) (-3,0)};
\end{axis}
\end{tikzpicture}\\[.4ex]
\multicolumn{2}{c}{\begin{tikzpicture}[baseline=(current bounding box.center)]
\begin{axis}[introplot,width=4.0cm,height=3.7cm,xtick={-10,-5,0},ytick={-10,-5,0,5},xmin=-13.2,xmax=4,ymin=-10.6,ymax=7.2,title={(c)}]
\addplot[only marks,mark=*,mark size=.64pt,spectralblue] coordinates {
(-12.6057936,-9.8535116) (-11.7213367,-8.4787525) (-11.4466865,-2.4941792) (-11.0527424,-7.4289292) (-10.9953295,-2.0318437)
(-10.6610656,-1.66813) (-10.501689,-6.5545871) (-10.3892186,-1.3556911) (-10.2114958,5.0952374) (-10.2105801,5.7347621)
(-10.1966172,4.5711469) (-10.1815155,6.5877296) (-10.1711193,4.1236426) (-10.1578489,-1.0763504) (-10.1378579,3.7320188)
(-10.0986403,3.3834424) (-10.0547095,3.0691778) (-10.0295235,-5.7969858) (-10.0069702,2.7828859) (-9.9561076,2.5197504)
(-9.9552078,-0.8208595) (-9.9026547,2.2759888) (-9.847033,2.0485579) (-9.7895775,1.8349686) (-9.7740658,-0.5838171)
(-9.7305537,1.6331616) (-9.6701689,1.441422) (-9.6158692,-5.125125) (-9.609603,-0.36179) (-9.608582,1.2583156)
(-9.5459109,1.0826412) (-9.4822399,0.9133919) (-9.4584373,-0.1524627) (-9.417626,0.7497237) (-9.3521062,0.5909301)
(-9.3181231,0.0458191) (-9.2856088,0.4364426) (-9.2482768,-4.520018) (-9.2217115,0.2806227) (-9.1657244,0.2043175)
(-8.9942189,0.1524907) (-8.9184225,-3.9689204) (-8.8190858,0.0850872) (-8.6418994,0.0157931) (-8.6203926,-3.4627256)
(-8.4626283,-0.0552994) (-8.349798,-2.9946229) (-8.2813301,-0.128113) (-8.1032739,-2.5593378) (-8.0980648,-0.2025711)
(-7.9128938,-0.2785968) (-7.8781756,-2.1526704) (-7.7258803,-0.3561138) (-7.6723835,-1.7711995) (-7.537085,-0.4350458)
(-7.484173,-1.4120831) (-7.3463982,-0.5152507) (-7.312123,-1.0729201) (-7.1545391,-0.7510461) (-7.1481236,-0.5903036)
(-6.957149,-0.6000669) (-6.7824084,-0.5733349) (-6.6036379,-0.5403212) (-6.4224999,-0.5038343) (-6.2390989,-0.4636579)
(-6.0533735,-0.4197165) (-5.8652678,-0.3719444) (-5.6747222,-0.3202782) (-5.4816731,-0.2646565) (-5.2860529,-0.2050196)
(-5.0877898,-0.1413087) (-4.8868071,-0.0734648) (-4.6830229,-0.0014282) (-4.4763494,0.0748628) (-4.2666913,0.1554733)
(-4.053945,0.2404728) (-3.8379966,0.329937) (-3.6187203,0.4239497) (-3.3959754,0.5226045) (-3.1696033,0.6260078)
(-2.939424,0.7342814) (-2.705231,0.8475662) (-2.4667855,0.9660267) (-2.223809,1.0898564) (-1.9759733,1.219285)
(-1.7228888,1.3545872) (-1.4640879,1.4960943) (-1.1990037,1.64421) (-0.9269407,1.7994315) (-0.6470351,1.9623781)
(-0.3581976,2.1338328) (-0.0590313,2.3148015) (0.252293,2.5066038) (0.5782316,2.7110139) (0.9221906,2.9304955)
(1.2890976,3.1686167) (1.6865137,3.4308583) (2.1271019,3.7264063) (2.6353656,4.0730398) (3.274716,4.5167238)};
\end{axis}
\end{tikzpicture}}
\end{tabular}
\caption{Eigenvalues of the three unperturbed tridiagonal samples,
$N=100$. Each blue dot represents an eigenvalue in the complex plane.
The horizontal and vertical axes give its real and imaginary parts.
The profiles $(b,c)$, $(b_2,c_2)$ and $(b_3,c_3)$ of (a)--(c) are given
in the text.}
\label{fig:clouds}
\end{figure}
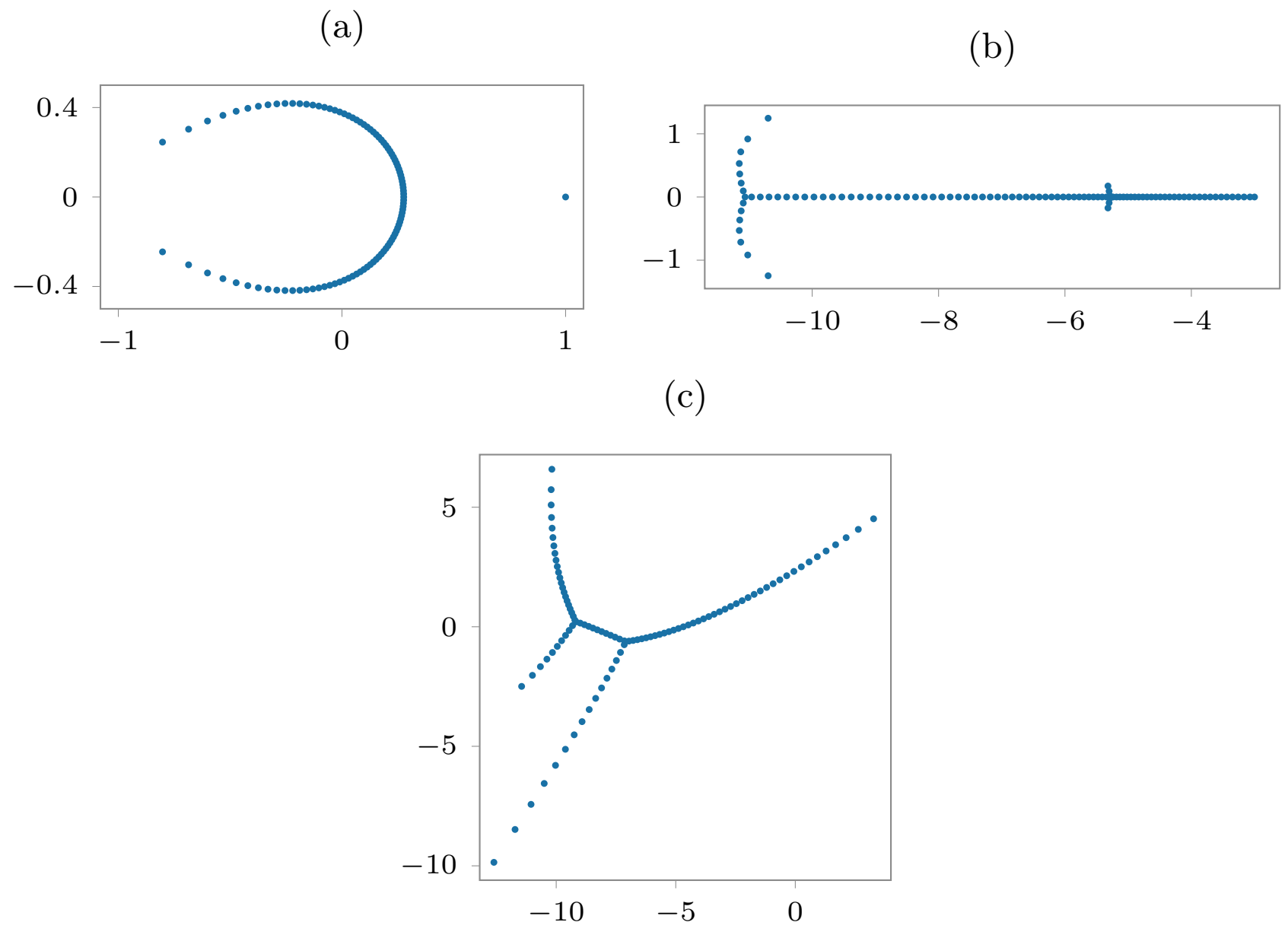

\begin{figure}[!htbp]
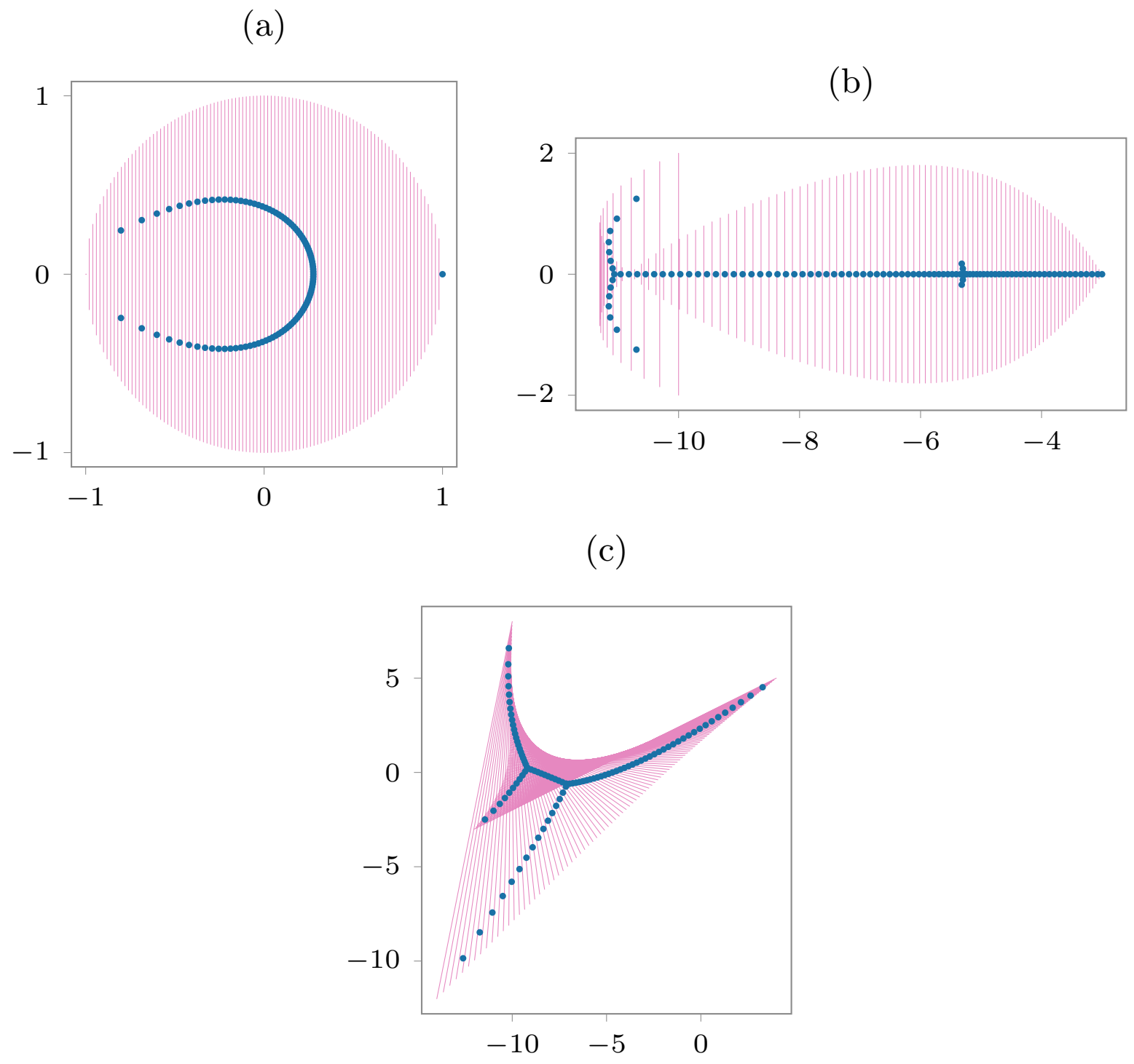

\centering

\caption{The eigenvalues of Figure~\ref{fig:clouds} (blue) overlaid on
a sample of the frozen segments $S_x$ (pink). Their union over all
$x\in[0,1]$ is the support of the standard measure. In (a) the segments
fill the unit disk. In (c) their direction changes along the lattice.}
\label{fig:standard}
\end{figure}

Figure~\ref{fig:clouds} shows the finite spectra alone, and
Figure~\ref{fig:standard} places the same points over the frozen
segments; the panels are labelled in the same order in both figures.
In (a) the choice $b(x)=1-2x$, $c(x)=\ii\sqrt{x(1-x)}$ makes the
geometry especially transparent. Writing $u=1-2x$, the endpoints of
$S_x$ are $u\pm\ii\sqrt{1-u^2}$, so the frozen segments are the vertical
chords of the unit disk. The eigenvalues follow two complex conjugate
arcs. The isolated point at $1$ has a simpler explanation: since
$c(0)=0$, the first site is decoupled and contributes its diagonal entry
$b(0)=1$ with weight $1/N$ in the counting measure. This is a feature
of the left-endpoint rule. With the midpoint samples
\eqref{eq:lt-array} the first edge carries $c(1/(2N))\ne0$, and the
corresponding eigenvalue is an outlier near $0.9318$.

In (b) the segments are again vertical, but their lengths now follow a
polynomial profile, and many eigenvalues cluster near a real interval.
In (c) both the direction and the centre of the segments vary, and the
eigenvalues follow several branches where the segments overlap.
For reference, the remaining profiles are
\[
\begin{aligned}
 b_2(x)&=P(x),\qquad c_2(x)=\ii(4x-4x^2-x^3),\\
 b_3(x)&=(-4+\ii)-2x-(3+\ii)x^2-(3+2\ii)x^3,\\
 c_3(x)&=(-4-2\ii)+(5+5\ii)x-(4+3\ii)x^2+(4+5\ii)x^3 .
\end{aligned}
\]
Here $P(x)=\sum_{j=0}^{14}p_jx^j$ with coefficients in ascending order
\[
 (p_0,\ldots,p_{14})=(-3,-5,-4,-1,-1,1,-4,-3,3,5,3,-2,0,-3,4).
\]

The region swept out by the frozen segments also has a spectral meaning
of its own. At its interior points the resolvent norm of a non-normal
locally Toeplitz matrix typically grows exponentially in $N$
\cite{TC04,TE05}, so the eigenvalues are extremely sensitive to
perturbations. In our earlier paper on tridiagonal locally Toeplitz
matrices \cite{GiSh26} we show that small structured random
perturbations typically lead to a two-dimensional limiting eigenvalue
distribution governed by the standard measure; for perturbations of
Toeplitz and twisted Toeplitz matrices by Gaussian noise see
\cite{BPZ19,SV21}. Here we study the unperturbed matrices. If their
limiting measure is carried by curves inside the same region, what
determines these curves and the density along them?

\subsection{Discrete Schr\"odinger operators with complex potential}\label{sec:intro-schr}
For real potentials the spectral theory of one-dimensional discrete
Schr\"odinger operators
\[
 (Hu)_k=u_{k-1}+u_{k+1}+V_ku_k
\]
is a highly developed subject, with landmarks such as Anderson
localization \cite{And58,AW15} and the theory of quasi-periodic
operators behind Hofstadter's butterfly \cite{Hof76,Jit99,AJ09,Avi15};
see \cite{DF22} for a modern account. Most of its tools rely on
self-adjointness: the spectral theorem and spectral measures,
Weyl--Titchmarsh theory, oscillation theory, and the variational
characterization and interlacing of eigenvalues. For complex potentials
none of these is available, and the spectral theory of discrete
Schr\"odinger operators with complex potential has been, and remains, a
serious challenge. The spectrum of a non-self-adjoint operator can be
highly unstable under small perturbations, and the resolvent can be
very large far away from the spectrum \cite{Dav02,Dav07,TE05}. Even the
location of the spectrum is a delicate question: spectral enclosures
and eigenvalue bounds of Lieb--Thirring type for complex Jacobi and
discrete Schr\"odinger operators are a subject of current research
\cite{GoKu07,BGK09,IbSt,BoSt21,KLS22}; see also the survey
\cite{Bec01}. Non-Hermitian random models, beginning with the work of
Hatano and Nelson, showed that the loss of self-adjointness changes the
spectral picture qualitatively. Eigenvalues of periodic truncations
leave the real axis and accumulate on curves \cite{HN96,GK98,GK03}, and
the spectra of random non-self-adjoint operators can be two-dimensional
sets with an intricate structure \cite{Dav01a,Dav01b,CWD12}. The
spectra of non-self-adjoint quasi-periodic operators, studied since the
work of Sarnak \cite{Sar82}, are still being determined; see
\cite{Lon19,BO26,KO26} for recent progress.

Finite truncations add a further difficulty. For non-normal operators
the eigenvalues of large truncations need not approximate the spectrum
of the infinite operator. Already for banded Toeplitz matrices they
accumulate on the Schmidt--Spitzer curves, a finite union of arcs inside
the spectrum of the Toeplitz operator, which consists of the symbol
curve together with the points around which it has non-zero winding
number \cite{SS60,RT92,BG05}. Boundary conditions therefore matter, and
the eigenvalue distribution of large finite matrices is an object of
its own.

In this paper the potential is neither random nor quasi-periodic but
slowly varying. We sample a fixed function, $V_k=V(k/N)$, on a lattice
of mesh $1/N$, keep the hopping coefficient constant, and impose
Dirichlet conditions at both ends. This is the symmetric tridiagonal
case of Section~\ref{sec:examples} with $b=V$ and constant $c$, and a
semiclassical regime with Planck constant $h=1/N$. Our main tool is the
complex WKB method for difference equations developed by Buslaev,
Fedotov, Klopp and others \cite{BF95,FK19,Fed23}; see also \cite{HS88}.
For real $V$ the matrices are Hermitian, and their eigenvalue
distributions converge to the standard measure \cite{Kac53,Tilli,KVA99}.
For complex $V$, little is known beyond exactly solvable cases such as
the Lommel and Charlier polynomials studied by Blaschke and \v{S}tampach
\cite{BlSt,BlStC}. For the class of potentials considered below we show
that the limit is a mother body of the standard measure, carried by
finitely many analytic arcs, and we give an explicit formula for its
logarithmic potential.

\subsection{Main results}\label{sec:mainresults}
Let $\alpha>0$, let $V\colon[0,1]\to\C$ be continuous, and let $H_N$ be
the $N\times N$ matrix with diagonal entries $V(k/N)$, $0\le k<N$, and
all entries on the two adjacent diagonals equal to $\alpha$. It is the
locally Toeplitz matrix of the symbol $a(x,z)=V(x)+\alpha(z+z^{-1})$,
and also the restriction to $N$ sites, with Dirichlet conditions, of the
discrete Schr\"odinger operator with potential $V(k/N)$. The frozen
segments $V(x)+[-2\alpha,2\alpha]$ are horizontal and all have length
$4\alpha$, and their centres trace the curve $V([0,1])$. The standard
measure $\M$ is therefore carried by the tube $\Sweep$ swept out by
these segments.

\emph{Monotone imaginary part.} Suppose that $(\Im V)'>0$ on $[0,1]$
and that $V$ extends holomorphically to a suitable complex neighbourhood
of $[0,1]$; these are hypotheses (H1) and (H2) of
Section~\ref{sec:operator}. The Main Theorem in
Section~\ref{sec:operator} states that the eigenvalue distributions
$\mu_N$ of $H_N$ converge weakly to a probability measure $\mu$,
supported on finitely many real-analytic arcs, with
\begin{equation}\label{eq:intro-main}
 \Pot{\mu}=\Pot{\M}-2\min\{A_L,A_R,U_+,U_-\}\quad\text{on }\Sweep,
 \qquad
 \Pot{\mu}=\Pot{\M}\quad\text{off }\Sweep .
\end{equation}
Here $A_L$, $A_R$, $U_+$ and $U_-$ are explicit tunnelling actions,
given by imaginary parts of integrals of the frozen momentum $p$,
defined by $2\alpha\cos p=E-V$. The actions $A_L$ and $A_R$ are
integrals along the real position interval, to the left and to the right
of the site at which the horizontal line through $E$ meets the potential
curve. The actions $U_+$ and $U_-$ are integrals along curves of real
momentum leading to the complex turning points above and below that
interval. The support of $\mu$ lies on the curves where two of the four
actions tie, and the density of $\mu$ along these curves is explicit
(Corollary~\ref{cor:curves}). In particular, $\mu$ is a mother body of
$\M$ (Corollary~\ref{cor:MB}), so all three parts of
Conjecture~\ref{conj:LT} hold for this class.

\emph{Sharp asymptotics.} The Main Theorem follows from an asymptotic
expansion of the Dirichlet determinant $\det(E-H_N)$, valid for $E$ in
the interior of $\Sweep$ outside a closed set of zero area
(Theorem~\ref{thm:sharp}).
The expansion has four terms, each with relative error $O(1/N)$. Two of
them come from the two frozen WKB solutions, and two from switches
between these solutions at the turning points above and below the real
interval. To obtain them we connect solutions normalized at the two
Dirichlet ends through a periodic Wronskian, and we separate its Fourier
modes by continuing above and below the turning points. The analytic
input is the complex WKB method for difference equations
\cite{BF95,FK19,KO26}.

\emph{Non-monotone potentials.} When $\Im V$ is not monotone, a
horizontal line of energies can meet the potential curve several times.
We prove an exact recursion for the Dirichlet determinant with one
switching coefficient per crossing (Theorem~\ref{thm:recursion}). Its
expansion attaches a cost to every history of dominant and recessive
roots. Under a holomorphic extension hypothesis, if the switching
coefficients obey the natural asymptotics \textup{(Sw)} and distinct
histories have no persistent ties, then $\mu_N$ converges, and
$\Pot{\M}-\Pot{\mu}$ is twice the minimal cost almost everywhere on
$\Sweep$ (Theorem~\ref{thm:conditional}). A geometric condition, Stokes
regularity, implies \textup{(Sw)} (Theorem~\ref{thm:stokesregular}). We
verify it almost everywhere for the circular arc
$V(x)=-\ii e^{3\pi\ii x/2}$ (Theorem~\ref{thm:circle} and
Corollary~\ref{cor:circle}). The general case is formulated as
Conjecture~\ref{conj:general}, and we record partial results towards it.

\emph{Special cases and numerics.} For purely imaginary potentials a
reflection symmetry simplifies the formula, and affine potentials
connect our results with the Lommel polynomials of \cite{BlSt}
(Section~\ref{sec:real}). We compare all predictions with numerically
computed spectra in Sections~\ref{sec:numerics}
and~\ref{sec:numgeneral}.

\subsection{Organization of the paper}
Section~\ref{sec:proofLT} contains the proof of
Theorem~\ref{thm:LTfacts}. Section~\ref{sec:operator} introduces the
Dirichlet discrete Schr\"odinger operator, states the Main Theorem and
Theorem~\ref{thm:sharp}, and gives an overview of the case of several
crossings. Sections~\ref{sec:setting}--\ref{sec:WKB} collect general
facts on Dirichlet determinants and the standard measure, the geometry
of the complex momentum and the four actions, and the complex WKB tools.
Section~\ref{sec:proof} proves Theorem~\ref{thm:sharp} and the Main
Theorem. Section~\ref{sec:real} treats purely imaginary and affine
potentials, and Section~\ref{sec:numerics} compares the results with
numerical spectra. Section~\ref{sec:general} deals with potentials whose
imaginary part is not monotone, and Section~\ref{sec:outlook} lists open
questions. 

\FloatBarrier
\section{Proof of Theorem~\ref{thm:LTfacts}}\label{sec:proofLT}

In this section $a$ is a banded symbol as in Theorem~\ref{thm:LTfacts}.
We put $P_N(E)=\det(E-T_N(a))$, so that
$\Pot{\mu_N(a)}(E)=\frac1N\log|P_N(E)|$. Part~(1) follows from the
row-sum bound. For part~(2) we compare $P_N$ with a product of frozen
transfer matrices, whose growth is controlled by the
Schmidt--Spitzer--Hirschman theorem. Parts~(3) and~(4) are formal
consequences of~(2), and part~(5) is an explicit computation.

\subsection{Compactness}
By the row-sum bound, every eigenvalue of $T_N(a)$ lies in the disk
$|E|\le R=\sum_j\|a_j\|_\infty$. Probability measures on a compact disk
form a weakly compact set. This proves part~(1).

\subsection{Frozen potentials and transfer matrices}
For fixed $x$ and $E$, let $W(x,E)$ be $a_d(x)$ times the $d$-th
compound of the companion matrix of the recurrence
$\sum_{j=-d}^{d}(a_j(x)-E\delta_{j0})u_{n+j}=0$. Its entries are
polynomials in $E$ and the $a_j(x)$, so $W$ is defined and continuous
also where $a_d(x)=0$. If $a_d(x)\ne0$, its spectral radius is
$\rho(W(x,E))=|a_d(x)|\prod_{j=1}^{d}|z_j(x,E)|$, where
$|z_1|\ge\dots\ge|z_{2d}|$ are the roots of $z^d(a(x,z)-E)$. By the
theorem of Schmidt--Spitzer and Hirschman \cite{SS60,Hi67,BG05}, applied
to the active band at each $x$,
\[
 \Pot{\nu_x}(E)=\log\rho(W(x,E))
\]
for every $x\in[0,1]$ and almost every $E$. In particular
$x\mapsto\nu_x$ is weakly continuous, and integration in $x$ gives
$\Pot{\M(a)}(E)=\int_0^1\log\rho(W(x,E))\,dx$ for almost every $E$.

Expanding $P_N(E)$ along the band shows that, up to sign, it is a fixed
corner entry of the ordered product $W_{N,N-1}(E)\cdots W_{N,0}(E)$.
Here $W_{N,i}(E)$ is defined like $W(x,E)$, with the coefficients in the
window of row $i$ taken at their sampling points $(i+j)/(2N)$. By
uniform continuity of the coefficients,
$\max_i\|W_{N,i}(E)-W(i/N,E)\|\to0$.

\subsection{The upper bound}
Fix $E$ and $\varepsilon>0$. Since $x\mapsto W(x,E)$ and the spectral
radius are continuous, $[0,1]$ splits into finitely many intervals with
a norm on each, adapted to $W$ at one point of the interval, such that
$\|W(y,E)\|\le\rho(W(y,E))+\varepsilon$ for $y$ in the interval.
Changing norms at the finitely many junctions costs a bounded factor,
so
\[
 \limsup_{N\to\infty}\frac1N\log|P_N(E)|\le
 \int_0^1\log\bigl(\rho(W(x,E))+\varepsilon\bigr)\,dx .
\]
Letting $\varepsilon\to0$ gives
$\limsup_N\Pot{\mu_N(a)}(E)\le\Pot{\M(a)}(E)$ for almost every $E$.

Now let $\nu$ be the weak limit of a subsequence $\mu_{N_k}(a)$. All
these measures are supported in $|E|\le R$, so
$\Pot{\mu_{N_k}(a)}\to\Pot{\nu}$ in $L^1_{\mathrm{loc}}$ \cite{Ra95},
hence almost everywhere along a further subsequence. Therefore
$\Pot{\nu}\le\Pot{\M(a)}$ almost everywhere. The potential of a
compactly supported finite measure is, at every point, the limit of its
averages over small disks about that point, so the inequality holds
everywhere.

\subsection{The exterior}
Let $O=\C\setminus\supp\M(a)$. By the continuity of $x\mapsto\nu_x$,
$\supp\M(a)$ is the closure of $\bigcup_x\supp\nu_x$, so each $E\in O$
lies outside every $\supp\nu_x$. Off $\supp\nu_x$ the Schmidt--Spitzer
theorem gives $|z_d(x,E)|>|z_{d+1}(x,E)|$, so $W(x,E)$ has a simple
dominant eigenvalue for every $x\in[0,1]$. Outside a further set of zero
area, the dominant eigenvector has a nonzero component along the
boundary coordinate of the corner entry for every $x$. For such $E$,
conjugating each factor by a continuous choice of dominant eigenvector
shows that the corner entry of the product has modulus
$\exp\bigl(N\int_0^1\log\rho(W(x,E))\,dx+o(N)\bigr)$. Hence
$\Pot{\mu_N(a)}(E)\to\Pot{\M(a)}(E)$ for almost every $E\in O$, and
$\Pot{\nu}=\Pot{\M(a)}$ almost everywhere on $O$. Since $\Pot{\M(a)}$ is
harmonic on $O$, the distributional Laplacian $2\pi\nu$ of $\Pot{\nu}$
vanishes there. Thus $\nu(O)=0$, that is, $\supp\nu\subseteq\supp\M(a)$.
Then $\Pot{\nu}$ is harmonic on $O$ as well, and the equality holds
everywhere on $O$. This completes the proof of part~(2).

\subsection{Mother bodies and the zero-area case}
For part~(3), both measures have mass one, and
$\supp\nu\subseteq\supp\M(a)\subseteq\overline U$ by~(2). Also
$\C\setminus\overline U\subseteq O$, so (2) gives $\Pot{\nu}=\Pot{\M(a)}$
off $\overline U$ and $\Pot{\nu}\le\Pot{\M(a)}$ on $\C$. These are
properties (i)--(iii) of Definition~\ref{def:MB}.

For part~(4), suppose that $\supp\M(a)$ has zero area. Then (2) gives
$\Pot{\nu}=\Pot{\M(a)}$ almost everywhere. Applying the Laplacian, we
get $\nu=\M(a)$. Thus every weak subsequential limit equals $\M(a)$,
and by (1) the whole sequence converges.

\subsection{A symbol without a limit}
For part~(5), only even offsets occur, so $T_N(a)$ is the direct sum of
its restrictions to the even and to the odd indices. For $N=2m$ the
restriction to the odd indices has characteristic polynomial $X^m$.
Hence $\mu_{2m}(a)(\{0\})\ge\frac12$, and by the Portmanteau theorem
every weak limit along even sizes has an atom of mass at least
$\frac12$ at the origin. For $N=2m+1$ the determinants of both
restrictions at $E=0$ can be evaluated in closed form: with
$s=m+\frac12$ their moduli are $s^{-n}\prod_{j=1}^{n}|s-j|$ for $n=m$
and $n=m+1$. Their product is at least $(6e)^{-(2m+1)}$, so
$\int\log|w|\,d\mu_{2m+1}(a)(w)\ge-\log(6e)$. Since $\log|w|$ is upper
semicontinuous and all eigenvalues lie in a fixed disk, every weak limit
$\nu$ along odd sizes satisfies $\int\log|w|\,d\nu(w)\ge-\log(6e)$, and
in particular $\nu(\{0\})=0$. So the limits along even and odd sizes
differ. If $\supp\M(a)$ had zero area, (4) would give convergence, so
its area is positive. This completes the proof of
Theorem~\ref{thm:LTfacts}.

\section{The Dirichlet discrete Schr\"odinger operator}\label{sec:operator}
Fix $\alpha>0$ and let $V:[0,1]\to\C$ be continuous. Take $c(x)=\alpha$
and $b(x)=V(x)$ in the tridiagonal symbol of Section~\ref{sec:examples}.
The resulting $N\times N$ matrix acts on the sites $0,\ldots,N-1$ by
\begin{equation}\label{eq:TN}
 (H_Nu)_k=\alpha(u_{k-1}+u_{k+1})+V(k/N)u_k,\qquad u_{-1}=u_N=0 .
\end{equation}
Thus $V$ occupies the main diagonal, both adjacent diagonals equal
$\alpha$, and the endpoint conditions are Dirichlet. The off-diagonal part
follows the adjacency convention for the free discrete Schr\"odinger operator
\cite{Teschl,IbSt}. With $\Delta_{\mathrm d}u_k=u_{k+1}-2u_k+u_{k-1}$,
the unitary change $u_k\mapsto(-1)^ku_k$ gives the equivalent form
$-\alpha\Delta_{\mathrm d}+V-2\alpha I$; the same change shows why we may
take the common real off-diagonal coefficient to be positive.

We count eigenvalues with algebraic multiplicity and write
\[
 p_N(E)=\det(EI-H_N),\qquad
 \mu_N=N^{-1}\sum_{p_N(E)=0}\delta_E .
\]
Here the frozen segments all have length $4\alpha$ and point in the
same direction, and their centres trace $V([0,1])$. We write $\M$ for the
standard measure and $\Sweep$ for the region swept out by the segments, so
\begin{equation}\label{eq:physical-sweep}
 \M=\int_0^1\nu_x\,dx,\qquad
 \Sweep=\bigcup_{0\le x\le1}\bigl(V(x)+[-2\alpha,2\alpha]\bigr).
\end{equation}
For real $V$, the matrices are self-adjoint and $\mu_N$ converges to
$\M$. Adding a constant imaginary part merely translates this picture.
When the imaginary part varies from site to site, $H_N$ is non-normal.
Indeed, its Hermitian real part commutes with
$\operatorname{diag}(\Im V(k/N))$ exactly when neighbouring sampled
imaginary parts agree. It is in this non-normal case that the support of
$\mu_N$ can develop the patterns seen in Section~\ref{sec:intro}.

Two solvable examples put the problem in context. An affine purely
imaginary potential leads to the Lommel polynomials of Blaschke and
\v{S}tampach \cite{BlSt}. Their Charlier-polynomial example
\cite{BlStC} has a varying product of adjacent coefficients and so
belongs to the broader tridiagonal setting. Boundary conditions also
matter: the periodic cosine models in \cite{KO26,BO26} lead to
a different characteristic equation, whereas here it is the two Dirichlet
endpoints that enter the actions in the Main Theorem.

\subsection{Monotone imaginary potential}
The geometry is simplest when the imaginary part of $V$ orders the frozen
segments. For a complex position $z$ in the domain of $V$, write
\[
 f(z)=-\Im V(z),\qquad r(z)=\Re V(z).
\]
The minus sign makes $f$ decrease along the lattice in the orientation
used below. We write energies as $E=\tau-\ii\sigma$, so $\tau=\Re E$ and
$\sigma=-\Im E$. The spectral plots use the ordinary axes $(\Re E,\Im E)$.
When we say that a \emph{turning point} lies in the upper or lower
half-plane, we always refer to the position variable $z$, not to the
energy $E$.

\begin{enumerate}
\item[(H1)] $V$ is complex-valued and real-analytic on $[0,1]$, with
$(\Im V)'(x)>0$ there.
\end{enumerate}
Set $\sigma_0=f(0)>\sigma_1=f(1)$. For
$\sigma\in[\sigma_1,\sigma_0]$ let $x_c=x_c(\sigma)$ solve $f(x_c)=\sigma$
and set $g(\sigma)=r(x_c)$. Then
\begin{equation}\label{eq:sweep}
 \Sweep=\Sweep_{2\alpha},\qquad
 \Sweep_r=\{\tau-\ii\sigma:\ \sigma_1\le\sigma\le\sigma_0,\quad
                |\tau-g(\sigma)|\le r\}.
\end{equation}
Thus $g(\sigma)$ is the real coordinate of the centre of the corresponding
horizontal frozen segment. The next assumption gives us room to continue the
momentum to both turning points.
\begin{enumerate}
\item[(H2)] $V$ extends holomorphically to a convex domain $\Omega$ containing
$[0,1]$, with $\Im V'>0$ there and $\Sweep_{4\alpha}\subset V(\Omega)$.
Writing $Y=\sup\{|\Im z|:z\in\Omega,\ V(z)\in\Sweep_{4\alpha}\}$, we
require $[0,1]+\ii[-Y,Y]\subset\Omega$.
\end{enumerate}
Under (H2) the function $V$ is injective: for distinct $z,w\in\Omega$ the
quotient $(V(z)-V(w))/(z-w)$ is the average of $V'$ over the segment joining
them and therefore has positive imaginary part. Write $\beta=V^{-1}$. Condition (H2) holds for sufficiently
small $\alpha$ under (H1), and for every $\alpha$ for affine potentials
satisfying (H1); see Remark~\ref{rem:H2}. A potential with
$(\Im V)'<0$ reduces to this orientation by the unitary identity
$-H_N(V)=D H_N(-V)D^{-1}$, $D_{kk}=(-1)^k$, together with the point
reflection $E\mapsto-E$ of the energy. A real constant may also be added to
both $V$ and $E$ without changing the analysis.

For $E\in\Sweep^\circ$, the frozen momentum on the real interval is the
continuous solution of
\begin{equation}\label{eq:momentum}
 2\alpha\cos p(x,E)=E-V(x),\qquad 0<\Re p(x,E)<\pi.
\end{equation}
The two frozen roots have the same modulus at $x_c$. Define
\begin{equation}\label{eq:actions}
\begin{gathered}
 A_L=\int_0^{x_c}|\Im p(x,E)|\,dx,\qquad
 A_R=\int_{x_c}^1|\Im p(x,E)|\,dx,\\
 \begin{aligned}
 U_+&=\int_{\tau-2\alpha}^{g(\sigma)}
       \arccos\frac{\tau-t}{2\alpha}\,|\Im\beta'(t-\ii\sigma)|\,dt,\\
 U_-&=\int_{g(\sigma)}^{\tau+2\alpha}
       \arccos\frac{t-\tau}{2\alpha}\,|\Im\beta'(t-\ii\sigma)|\,dt .
 \end{aligned}
\end{gathered}
\end{equation}
Here and below, the energy is omitted from the notation for the actions
when it is fixed. The turning points
$x_+=\beta(E-2\alpha)$ and $x_-=\beta(E+2\alpha)$ lie above and below the
real position interval, respectively. The inverse images of the two
horizontal energy segments from $V(x_c)$ to $E\mp2\alpha$ are curves on
which the momentum is real, and $U_+$ and $U_-$ are the integrals of
$p\,|dy|$ and $(\pi-p)\,|dy|$ along them. All four actions are strictly positive in
$\Sweep^\circ$ and extend continuously to its boundary.

\begin{maintheorem}
Assume \textup{(H1)} and \textup{(H2)}.  Then $\mu_N$ converges weakly to a
probability measure $\mu$ whose logarithmic potential is
\begin{equation}\label{eq:mainformula}
 \Pot{\mu}(E)=
 \begin{cases}
 \Pot{\M}(E)-2\min\{A_L,A_R,U_+,U_-\},
 &E\in\Sweep,\\[2pt]
 \Pot{\M}(E), &E\notin\Sweep .
 \end{cases}
\end{equation}
Equivalently,
$\Pot{\mu}=\log\alpha+\max\{|A_L-A_R|,\,A_L+A_R-2U_+,\,A_L+A_R-2U_-\}$ on
$\Sweep$.  The support of $\mu$ is a finite union of real-analytic arcs, contained in the closure of the set of
$E\in\Sweep^\circ$ at which the minimum in \eqref{eq:mainformula} is
attained at least twice.  The densities on these arcs are given in
Corollary~\ref{cor:curves}.
\end{maintheorem}

\begin{corollary}\label{cor:MB}
Under \textup{(H1)} and \textup{(H2)}, $\mu$ is a mother body of $\M$.  Thus,
under these explicit hypotheses, this family realizes the three predictions of
Conjecture~\ref{conj:LT}.
\end{corollary}

\begin{proof}
The measure $\M$ is supported on $\Sweep$, whose boundary consists of the curves
$V(x)\pm2\alpha$ and two horizontal segments, and its density
\eqref{eq:Mdensity} is real-analytic in $\Sweep^\circ$.  Property (i) is part
of the Main Theorem.  Properties (ii) and (iii) follow from
\eqref{eq:mainformula}, since the four actions are non-negative.
\end{proof}

The formula also explains why curves appear in the numerical spectra.
Where one action is strictly smallest, the limiting potential is the
real part of a single analytic phase and hence is harmonic. Since its
Laplacian is $2\pi\mu$, mass can occur in the interior only where at least
two actions attain the minimum. On a regular arc where two phases
exchange dominance, the jump in the normal derivative gives the density
with respect to arclength (Corollary~\ref{cor:curves}).

The deficit $\Pot{\M}-\Pot{\mu}$, which measures the non-normality of
$H_N$ (Remark~\ref{rem:deficit}), is twice the smallest of the four
tunnelling actions. It is positive in $\Sweep^\circ$ and vanishes on $\partial\Sweep$: $A_L=0$ on
the edge $\sigma=\sigma_0$, $A_R=0$ on $\sigma=\sigma_1$, $U_+=0$ on the right
boundary $\tau=g+2\alpha$, and $U_-=0$ on the left boundary $\tau=g-2\alpha$.

The Main Theorem follows from sharp asymptotics of $p_N$.  For
$E\in\Sweep^\circ$ put
\begin{equation}\label{eq:IJ}
 I=\int_0^1p(x)\,dx,\qquad
 J_+=\int_0^{x_+}p(z)\,dz,\qquad J_-=\int_0^{x_-}\bigl(p(z)-\pi\bigr)\,dz ,
\end{equation}
where $p$ is continued to the complex domain $D$ of Lemma~\ref{lem:geometry}
and the integrals are taken along any paths in $D$ ending at $x_\pm$; they
converge, since $p$ is continuous up to $x_\pm$.  Fix a branch of
$(\sin p)^{1/2}$ on $D$ and put
\[
 A=\frac{e^{\ii p(0)}}{2\ii(\sin p(0))^{1/2}},\qquad
 B=-\frac{e^{-\ii p(0)}}{2\ii(\sin p(0))^{1/2}} .
\]

\begin{theorem}\label{thm:sharp}
Assume \textup{(H1)} and \textup{(H2)}.  There is a relatively closed set
$\mathcal E\subset\Sweep^\circ$ of zero area such that for
$E\in\Sweep^\circ\setminus \mathcal E$, as $N\to\infty$,
\begin{multline}\label{eq:sharp}
 \frac{(\sin p(1))^{1/2}\,p_N(E)}{\alpha^N}
 =Ae^{\ii NI}\bigl(1+O(\tfrac1N)\bigr)
 +Be^{-\ii NI}\bigl(1+O(\tfrac1N)\bigr)\\
 +\ii A\,e^{\ii N(2J_+-I)}\bigl(1+O(\tfrac1N)\bigr)
 -\ii A\,e^{\ii N(2J_--I)}\bigl(1+O(\tfrac1N)\bigr)
 +O\bigl(e^{-cN}M_N\bigr),
\end{multline}
where $M_N$ is the largest modulus of the four exponentials and
$c=c(E)>0$.
\end{theorem}

The error terms in \eqref{eq:sharp} are relative to the individual terms. At
a point where a single term dominates, the determinant is asymptotic to that
term, with relative error $O(N^{-1})$. On a tie curve two terms can cancel,
so there the expansion gives no relative estimate for the determinant as a
whole.

There are four possible leading phases. Without a switch, the two roots give
$\Theta_{1,2}=\log\alpha\pm\ii I$. A switch at the upper or lower turning
point gives, respectively,
\begin{equation}\label{eq:phases0}
 \Theta_3=\log\alpha+\ii(2J_+-I),\qquad
 \Theta_4=\log\alpha+\ii(2J_--I).
\end{equation}
The subtraction of $\pi$ in $J_-$ keeps track of the branch at which
$p=\pi$; it cannot be dropped from the complex phase. The real parts of
$\Theta_1,\ldots,\Theta_4$ are $\Pot{\M}$ minus $2A_R,2A_L,2U_+,2U_-$,
respectively, and their maximum gives \eqref{eq:mainformula}.

For purely imaginary $V$, reflection across the imaginary energy axis
interchanges the two turning-point costs. The smaller cost suffices, and
a tie between the two can occur only on that axis
(Section~\ref{sec:real}). For a general complex potential the curve on
which the two costs are equal can bend away from this axis; Figures~\ref{fig:eigs} and~\ref{fig:complex}
show both cases.

\subsection{Several crossings}
When $\Im V$ is not monotone, a horizontal energy segment can meet the
potential curve at more than one site. The corresponding transfer problem
still admits an exact recursion, with one switching coefficient per crossing
(Theorem~\ref{thm:recursion}). Expanding this recursion suggests a cost for
each sequence of dominant and recessive roots. When the switching asymptotics
\textup{(Sw)} hold and distinct histories have no persistent equality of
costs, the minimum of those costs determines the limiting potential
(Theorem~\ref{thm:conditional}). A geometric condition called Stokes
regularity implies \textup{(Sw)}; we verify it almost everywhere for
$V(x)=-\ii e^{3\pi\ii x/2}$, although non-resonance remains a hypothesis in
that example. The spiral in Section~\ref{sec:numgeneral} illustrates the
open problem beyond these hypotheses.

Sections~\ref{sec:setting} and~\ref{sec:actions} develop the potential
theory and the four actions. The analytic input is the canonical-domain WKB
estimate together with its horizontal continuation
\cite[Appendix~A]{KO26}; Section~\ref{sec:WKB} gives the versions needed
here. The Dirichlet determinant requires connecting solutions
normalized at opposite ends of the interval. Section~\ref{sec:proof}
expresses this connection as a periodic Wronskian and separates its Fourier
modes by continuing above and below the turning points.

\section{Setting and general facts}\label{sec:setting}
\subsection{Dirichlet determinants and frozen roots}
For fixed $N$, let $P_k$ be the leading principal minor of order $k$ of
$EI-H_N$, a monic polynomial in the energy. Expansion along the last row gives
\begin{equation}\label{eq:recurrence}
 P_{k+1}=(E-V(k/N))P_k-\alpha^2P_{k-1},\qquad
 P_{-1}=0,\quad P_0=1,\quad p_N=P_N.
\end{equation}
The minus sign comes from multiplying the product of the two real
off-diagonal entries by the cofactor sign. At a fixed site $x$, the frozen
roots, that is, the roots of $\rho^2-(E-V(x))\rho+\alpha^2=0$, are
\begin{equation}\label{eq:roots}
 \rho_\pm=\alpha e^{\pm\ii p},\qquad
 |\rho_\pm|=\alpha e^{\mp\Im p},\qquad \rho_+\rho_-=\alpha^2 .
\end{equation}
They coincide precisely when $V(x)=E\mp2\alpha$.

\begin{lemma}\label{lem:momentum}
Assume \textup{(H1)} and let $E=\tau-\ii\sigma\in\Sweep^\circ$.  Then
there are no turning points on $[0,1]$, and there is a unique continuous
solution $p$ of \eqref{eq:momentum} on $[0,1]$ with $0<\Re p<\pi$.  It is
real-analytic in $x$, holomorphic in $E$, and
\[
 \Im p(x)<0\ \ (x<x_c),\qquad \Im p(x_c)=0,\qquad \Im p(x)>0\ \ (x>x_c).
\]
\end{lemma}

\begin{proof}
By \eqref{eq:momentum}, $\cos p=\mathfrak c(x):=\frac{\tau-r(x)}{2\alpha}
-\ii\frac{\sigma-f(x)}{2\alpha}$.  The cosine maps the strip
$\{0<\Re p<\pi\}$ biholomorphically onto
$\C\setminus\bigl((-\infty,-1]\cup[1,\infty)\bigr)$.  The value $\mathfrak c(x)$
is real only at $x=x_c$, where $\mathfrak c(x_c)=(\tau-g(\sigma))/(2\alpha)\in(-1,1)$
because $E\in\Sweep^\circ$.  Hence $\mathfrak c([0,1])$ avoids the two
half-lines, which gives existence, uniqueness and regularity, and shows that
$\mathfrak c(x)\ne\pm1$, so that there are no turning points.  Writing
$p=u+\ii v$, we have $\Im\cos p=-\sin u\sinh v$ with $\sin u>0$, so $v$ has the
sign of $\sigma-f(x)$, which is the sign of $x-x_c$.
\end{proof}

By the lemma, $\rho_+$ is the dominant root on $[0,x_c)$ and $\rho_-$ the
dominant root on $(x_c,1]$; the two have equal moduli exactly at the crossing
point $x_c$.

\subsection{The standard measure}
At each $x$ the frozen Dirichlet matrix has eigenvectors proportional to
$\sin((k+1)\pi j/(N+1))$ and eigenvalues
$V(x)+2\alpha\cos(\pi j/(N+1))$, $1\le j\le N$. Consequently, $\nu_x$
is the image of $d\theta/\pi$ on $[0,\pi]$ under
$\theta\mapsto V(x)+2\alpha\cos\theta$. Integrating the arcsine density
and changing variables $x\mapsto\sigma=f(x)$ in \eqref{eq:physical-sweep}
gives, under \textup{(H1)},
\begin{equation}\label{eq:Mdensity}
 \frac{d\M}{dA}(\tau-\ii\sigma)=
 \frac{1}{\pi\,|f'(x_c(\sigma))|\,
 \sqrt{4\alpha^2-(\tau-g(\sigma))^2}},\qquad E\in\Sweep^\circ .
\end{equation}
There is no mass on the boundary.

\begin{lemma}\label{lem:PotM}
For every $E\in\C$,
\[
 \begin{aligned}
 \Pot{\nu_x}(E)&=\log\max(|\rho_+|,|\rho_-|),\\
 \Pot{\M}(E)&=\int_0^1\log\max(|\rho_+|,|\rho_-|)\,dx .
 \end{aligned}
\]
In particular, under \textup{(H1)},
\begin{equation}\label{eq:PotM}
 \Pot{\M}(E)=\log\alpha+A_L(E)+A_R(E),\qquad E\in\Sweep^\circ .
\end{equation}
Off $V(x)+[-2\alpha,2\alpha]$ the Cauchy transform
$C_{\nu_x}(E)=\int d\nu_x(w)/(E-w)=2\partial_E\Pot{\nu_x}(E)$ is
$(\rho_{\rm dom}-\rho_{\rm rec})^{-1}$, where the dominant and recessive
roots are the roots of larger and smaller modulus, respectively.
\end{lemma}
\begin{proof}
Let $\rho$ be either root. Then $E-V(x)=\rho+\alpha^2/\rho$ and
\[
 E-V(x)-2\alpha\cos\theta=
 \frac{(\rho-\alpha e^{\ii\theta})(\rho-\alpha e^{-\ii\theta})}{\rho}.
\]
Averaging the logarithm of the modulus of this expression over the circle
gives $2\log\max(|\rho|,\alpha)-\log|\rho|$, which is the logarithm of the
modulus of the larger root, also in the equal-modulus case. This average is
computed from the circular mean of $\log|\rho-\alpha e^{\ii\theta}|$; in the
equal-modulus case its logarithmic singularity is integrable. Integrating over
$x$ and using \eqref{eq:roots} now gives \eqref{eq:PotM}. Off the frozen segment,
differentiating the quadratic yields
$2\partial_E\Pot{\nu_x}=\rho'_{\rm dom}/\rho_{\rm dom}
=(\rho_{\rm dom}-\rho_{\rm rec})^{-1}$.
\end{proof}

\subsection{Exponential sums}
The following principle turns pointwise asymptotics of polynomials into weak
convergence of their zero-counting measures.

\begin{proposition}\label{prop:expsum}
Let $q_N$ be monic polynomials of degree $N\to\infty$ whose zeros lie in a
fixed disk.  Put $v_N=N^{-1}\log|q_N|$ and let $\nu_N$ be the normalized
zero-counting measure of $q_N$.  Suppose that for almost every $E\in\C$
there are $J\in\N$, complex numbers $\Phi_1(E),\dots,\Phi_J(E)$
with pairwise distinct real parts, and complex numbers $a_{j,N}(E)$ with
$0<\liminf_N|a_{j,N}(E)|\le\limsup_N|a_{j,N}(E)|<\infty$, such
that
\begin{equation}\label{eq:expsum}
 q_N(E)=\sum_{j=1}^Je^{N\Phi_j(E)}\bigl(a_{j,N}(E)+o(1)\bigr)
 \qquad(N\to\infty).
\end{equation}
Then $v_N(E)\to v(E):=\max_j\Re\Phi_j(E)$ for almost every
$E$.  There is a subharmonic function $v^*$ with $v^*=v$ almost
everywhere such that $v_N\to v^*$ in $L^1_{\rm loc}(\C)$ and
$\nu_N\Rightarrow\nu:=(2\pi)^{-1}\Delta v^*$.  Moreover, $v^*=\Pot{\nu}$.
\end{proposition}

\begin{proof}
Fix such an energy $E$. For $J=1$ the pointwise estimate below is
immediate from the single term. For $J>1$, let $j_0$ be the index of the phase
with the largest real part, and put
$\eta=\Re\Phi_{j_0}-\max_{j\ne j_0}\Re\Phi_j>0$.  The other terms are
$O(e^{N(\Re\Phi_{j_0}-\eta)})$, so $q_N(E)=e^{N\Phi_{j_0}}(a_{j_0,N}+o(1))$,
and the last factor is eventually bounded above and below in modulus.  Hence
$v_N(E)=\Re\Phi_{j_0}(E)+O(N^{-1})$.

The functions $v_N$ are subharmonic and uniformly bounded above on compact
sets, and $v_N(E)=\log|E|+O(|E|^{-1})$ uniformly in $N$ for
large $|E|$.  By the compactness theorem for subharmonic functions
\cite[Thm.~4.1.9]{Ho83}, every subsequence has a further subsequence that
converges in $L^1_{\rm loc}$ to a subharmonic function $w$; convergence to
$-\infty$ is excluded by the behaviour at infinity.  A further subsequence
converges almost everywhere, so $w=v$ almost everywhere.  Two subharmonic
functions that agree almost everywhere agree everywhere, so all limit points
coincide; call the common limit $v^*$.  Thus $v_N\to v^*$ in $L^1_{\rm loc}$.
Since $\nu_N=(2\pi)^{-1}\Delta v_N$ in the sense of distributions and all
$\nu_N$ are supported in a fixed disk, $\nu_N\Rightarrow(2\pi)^{-1}\Delta v^*$.
Finally, $v_N=\Pot{\nu_N}$, so for $\phi\in C_c(\C)$ we have $\int\phi v_N\,dA
=\int(\log|\cdot|*\phi)\,d\nu_N\to\int(\log|\cdot|*\phi)\,d\nu=\int\phi\,
\Pot{\nu}\,dA$, because $\log|\cdot|*\phi$ is continuous.  Hence $v^*=\Pot{\nu}$
almost everywhere, and everywhere, as both are subharmonic.
\end{proof}

\begin{corollary}\label{cor:charts}
Assume the hypotheses of Proposition~\ref{prop:expsum}.  Let $F\subset\C$ be
a closed set of zero area, and suppose that on every component $W$ of
$\C\setminus F$ there are finitely many holomorphic functions
$\Phi_{W,1},\dots,\Phi_{W,J_W}$ with $\Re(\Phi_{W,j}-\Phi_{W,k})\not\equiv0$
for $j\ne k$, such that $v=\max_j\Re\Phi_{W,j}$ almost everywhere on $W$;
this holds, for instance, if the phases in \eqref{eq:expsum} can be taken to be
$\Phi_{W,1}(E),\dots,\Phi_{W,J_W}(E)$ for almost every
$E\in W$.  Then
$v^*=\max_j\Re\Phi_{W,j}$ on every $W$, and
\[
 \supp\nu\subseteq F\cup\bigcup_W\{E\in W:\text{the maximum is
 attained by at least two indices}\}.
\]
The intersection of the right-hand side with each $W$ is contained in a
real-analytic subset of $W$ of dimension at most one.  On an arc in $W$ where
exactly two phases $\Phi_{W,j},\Phi_{W,k}$ attain the maximum and
$\Phi_{W,j}'\ne\Phi_{W,k}'$, the measure $\nu$ has the density
\[
 \frac1{2\pi}\bigl|\Phi_{W,j}'-\Phi_{W,k}'\bigr|
\]
with respect to arclength.
\end{corollary}

\begin{proof}
On $W$ the function $\max_j\Re\Phi_{W,j}$ is continuous and subharmonic, and
it equals $v^*$ almost everywhere; hence it equals $v^*$ on $W$.  Off the tie
set, where the maximum is attained only once, it is locally harmonic, so $\nu$
vanishes there.  The tie set is contained
in $\bigcup_{j<k}\{\Re(\Phi_{W,j}-\Phi_{W,k})=0\}$, a finite union of zero sets
of harmonic functions that do not vanish identically.  Such sets are
real-analytic of dimension at most one.  Near an arc where exactly the phases
$j$ and $k$ attain the maximum, put $u=\Re(\Phi_{W,j}-\Phi_{W,k})$.  Then
$v^*=\Re\Phi_{W,k}+\max\{u,0\}$, and integration by parts on the two sides of
$\{u=0\}$ gives $\Delta\max\{u,0\}=|\nabla u|\,ds$.  By the Cauchy--Riemann
equations $|\nabla u|=|\Phi_{W,j}'-\Phi_{W,k}'|$.  These steps are the
elementary part of the theory of subharmonic configurations of
Bj\"ork, Borcea and B\o gvad \cite[Lemma~2.1, Prop.~2.3]{BBB}, applied
to the configuration $\max_j\Re\Phi_{W,j}$.
\end{proof}

The zeros of $p_N$ lie in the disk $|E|\le\|V\|_\infty+2\alpha$, since
$\|H_N\|\le\|V\|_\infty+2\alpha$.  Hence Proposition~\ref{prop:expsum} applies
to $q_N=p_N$ as soon as \eqref{eq:expsum} is known almost everywhere.  Outside
$\Sweep$ this follows from the absence of a crossing
(Proposition~\ref{prop:exterior}); the whole difficulty lies inside $\Sweep$.

\section{The momentum in the complex domain and the four actions}
\label{sec:actions}

\subsection{Geometry}
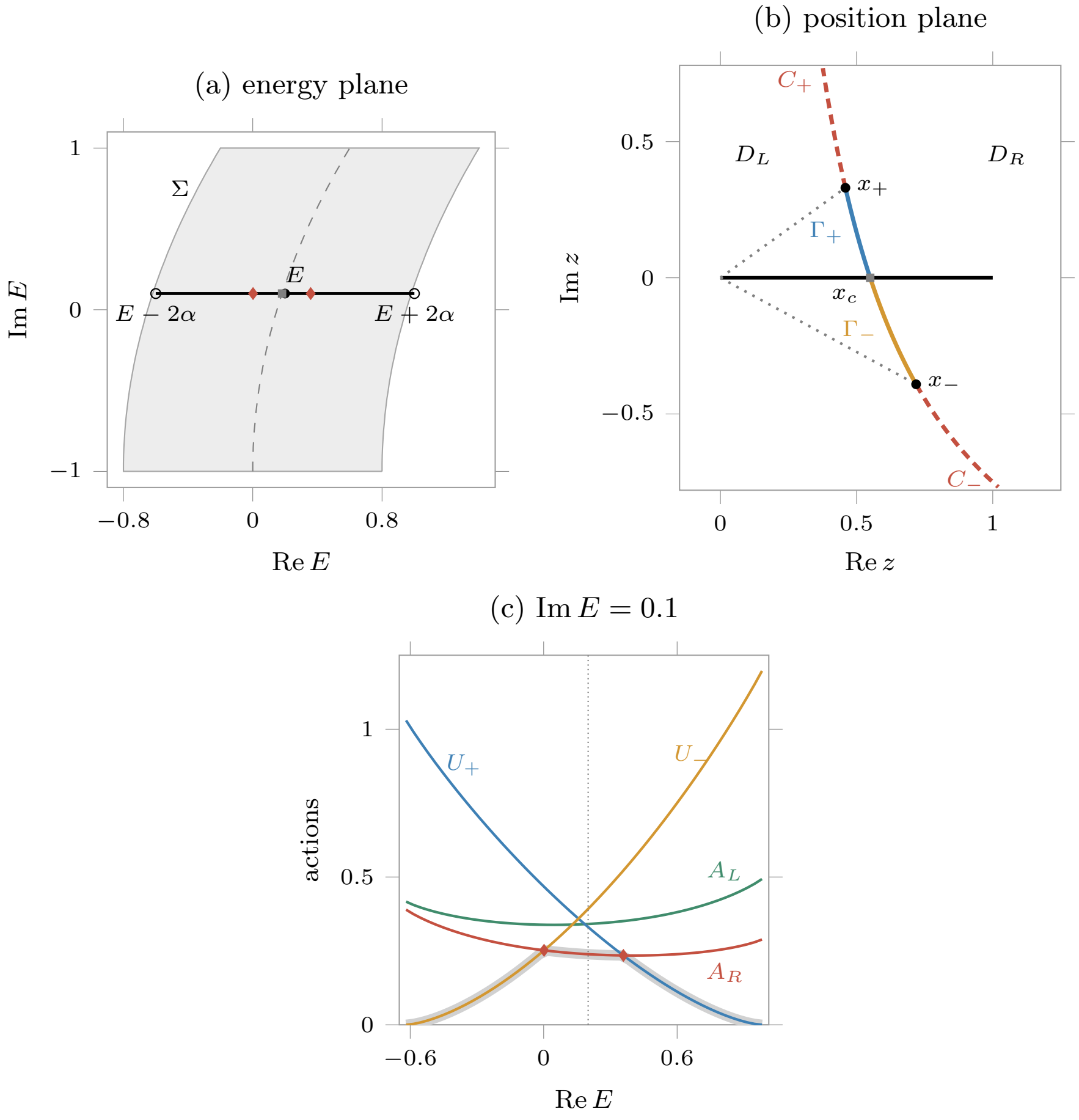
\begin{figure}[tbp]
\centering
\begin{tikzpicture}
\begin{axis}[paperplot,width=4.2cm,height=4.6cm,xmin=-.9,xmax=1.5,ymin=-1.1,ymax=1.1,xtick={-.8,0,.8},ytick={-1,0,1},title={(a) energy plane}]
\addplot[fill=gray!14,draw=gray!70,thin] coordinates {(0.800,-1.000) (0.800,-0.950) (0.801,-0.900) (0.803,-0.850) (0.806,-0.800) (0.809,-0.750) (0.814,-0.700) (0.818,-0.650) (0.824,-0.600) (0.830,-0.550) (0.838,-0.500) (0.845,-0.450) (0.854,-0.400) (0.863,-0.350) (0.874,-0.300) (0.884,-0.250) (0.896,-0.200) (0.908,-0.150) (0.922,-0.100) (0.935,-0.050) (0.950,0.000) (0.965,0.050) (0.982,0.100) (0.998,0.150) (1.016,0.200) (1.034,0.250) (1.054,0.300) (1.073,0.350) (1.094,0.400) (1.115,0.450) (1.137,0.500) (1.160,0.550) (1.184,0.600) (1.208,0.650) (1.234,0.700) (1.259,0.750) (1.286,0.800) (1.313,0.850) (1.342,0.900) (1.370,0.950) (1.400,1.000) (-0.200,1.000) (-0.230,0.950) (-0.258,0.900) (-0.287,0.850) (-0.314,0.800) (-0.341,0.750) (-0.366,0.700) (-0.392,0.650) (-0.416,0.600) (-0.440,0.550) (-0.463,0.500) (-0.485,0.450) (-0.506,0.400) (-0.527,0.350) (-0.546,0.300) (-0.566,0.250) (-0.584,0.200) (-0.602,0.150) (-0.619,0.100) (-0.635,0.050) (-0.650,0.000) (-0.665,-0.050) (-0.678,-0.100) (-0.692,-0.150) (-0.704,-0.200) (-0.716,-0.250) (-0.727,-0.300) (-0.737,-0.350) (-0.746,-0.400) (-0.755,-0.450) (-0.763,-0.500) (-0.770,-0.550) (-0.776,-0.600) (-0.782,-0.650) (-0.787,-0.700) (-0.791,-0.750) (-0.794,-0.800) (-0.797,-0.850) (-0.799,-0.900) (-0.800,-0.950) (-0.800,-1.000) (0.800,-1.000)};
\addplot[gray,dashed,thin] coordinates {(0.000,-1.000) (0.000,-0.950) (0.002,-0.900) (0.003,-0.850) (0.006,-0.800) (0.009,-0.750) (0.014,-0.700) (0.018,-0.650) (0.024,-0.600) (0.030,-0.550) (0.037,-0.500) (0.045,-0.450) (0.054,-0.400) (0.063,-0.350) (0.074,-0.300) (0.084,-0.250) (0.096,-0.200) (0.108,-0.150) (0.122,-0.100) (0.135,-0.050) (0.150,0.000) (0.165,0.050) (0.182,0.100) (0.198,0.150) (0.216,0.200) (0.234,0.250) (0.254,0.300) (0.273,0.350) (0.294,0.400) (0.315,0.450) (0.337,0.500) (0.360,0.550) (0.384,0.600) (0.408,0.650) (0.434,0.700) (0.459,0.750) (0.486,0.800) (0.513,0.850) (0.542,0.900) (0.570,0.950) (0.600,1.000)};
\addplot[black,line width=.9pt] coordinates {(-0.6000,0.1000) (1.0000,0.1000)};
\addplot[only marks,mark=*,mark size=1.3pt,black] coordinates {(0.2000,0.1000)};
\addplot[only marks,mark=o,mark size=1.5pt,black] coordinates {(-0.6000,0.1000) (1.0000,0.1000)};
\addplot[only marks,mark=square*,mark size=1.1pt,gray] coordinates {(0.1815,0.1000)};
\addplot[only marks,mark=diamond*,mark size=1.7pt,actionred] coordinates {(0.0018,0.1000) (0.3585,0.1000)};
\node[font=\scriptsize,above] at (axis cs:0.260,0.100) {$E$};
\node[font=\scriptsize,below] at (axis cs:-0.600,0.100) {$E-2\alpha$};
\node[font=\scriptsize,below] at (axis cs:1.000,0.100) {$E+2\alpha$};
\node[font=\scriptsize] at (axis cs:-.45,.75) {$\Sweep$};
\end{axis}
\end{tikzpicture}
\hspace{.2cm}
\begin{tikzpicture}
\begin{axis}[paperplot,xlabel={$\Re z$},ylabel={$\Im z$},width=4.2cm,height=4.6cm,xmin=-.15,xmax=1.25,ymin=-.78,ymax=.78,xtick={0,.5,1},ytick={-.5,0,.5},title={(b) position plane}]
\addplot[black,line width=1.1pt] coordinates {(0,0) (1,0)};
\addplot[gray,dotted,thick] coordinates {(0,0) (0.4590,0.3304)};
\addplot[gray,dotted,thick] coordinates {(0,0) (0.7186,-0.3910)};
\addplot[actionblue,line width=1.3pt] coordinates {(0.4590,0.3304) (0.4608,0.3225) (0.4627,0.3145) (0.4646,0.3065) (0.4665,0.2985) (0.4684,0.2904) (0.4703,0.2823) (0.4723,0.2742) (0.4743,0.2661) (0.4763,0.2579) (0.4783,0.2497) (0.4804,0.2415) (0.4825,0.2333) (0.4846,0.2250) (0.4867,0.2167) (0.4889,0.2084) (0.4911,0.2000) (0.4933,0.1917) (0.4955,0.1832) (0.4978,0.1748) (0.5001,0.1663) (0.5024,0.1578) (0.5048,0.1493) (0.5072,0.1408) (0.5096,0.1322) (0.5120,0.1236) (0.5145,0.1150) (0.5170,0.1063) (0.5196,0.0976) (0.5222,0.0889) (0.5248,0.0801) (0.5274,0.0713) (0.5301,0.0625) (0.5328,0.0537) (0.5356,0.0448) (0.5384,0.0359) (0.5412,0.0270) (0.5441,0.0180) (0.5470,0.0090) (0.5500,0.0000)};
\addplot[actiongold,line width=1.3pt] coordinates {(0.5500,0.0000) (0.5531,-0.0095) (0.5563,-0.0190) (0.5596,-0.0285) (0.5629,-0.0381) (0.5662,-0.0477) (0.5696,-0.0574) (0.5730,-0.0670) (0.5765,-0.0767) (0.5801,-0.0865) (0.5837,-0.0963) (0.5874,-0.1061) (0.5911,-0.1159) (0.5949,-0.1257) (0.5987,-0.1356) (0.6026,-0.1456) (0.6066,-0.1555) (0.6106,-0.1655) (0.6147,-0.1755) (0.6189,-0.1855) (0.6231,-0.1956) (0.6274,-0.2057) (0.6318,-0.2158) (0.6363,-0.2260) (0.6408,-0.2361) (0.6454,-0.2463) (0.6501,-0.2566) (0.6548,-0.2668) (0.6597,-0.2771) (0.6646,-0.2873) (0.6696,-0.2976) (0.6747,-0.3080) (0.6798,-0.3183) (0.6851,-0.3286) (0.6904,-0.3390) (0.6959,-0.3494) (0.7014,-0.3598) (0.7070,-0.3702) (0.7127,-0.3806) (0.7186,-0.3910)};
\addplot[actionred,line width=1.3pt,densely dashed] coordinates {(0.4590,0.3304) (0.4501,0.3701) (0.4416,0.4091) (0.4336,0.4476) (0.4259,0.4854) (0.4187,0.5227) (0.4118,0.5595) (0.4052,0.5957) (0.3989,0.6315) (0.3928,0.6667) (0.3871,0.7015) (0.3815,0.7359) (0.3762,0.7698)};
\addplot[actionred,line width=1.3pt,densely dashed] coordinates {(0.7186,-0.3910) (0.7329,-0.4160) (0.7479,-0.4410) (0.7634,-0.4659) (0.7796,-0.4908) (0.7963,-0.5155) (0.8137,-0.5401) (0.8317,-0.5646) (0.8504,-0.5887) (0.8696,-0.6126) (0.8895,-0.6361) (0.9100,-0.6593) (0.9310,-0.6821) (0.9526,-0.7044) (0.9746,-0.7262) (0.9972,-0.7474) (1.0202,-0.7681)};
\addplot[only marks,mark=*,mark size=1.3pt,black] coordinates {(0.4590,0.3304) (0.7186,-0.3910)};
\addplot[only marks,mark=square*,mark size=1.1pt,gray] coordinates {(0.5500,0)};
\node[font=\scriptsize,right] at (axis cs:0.459,0.330) {$x_+$};
\node[font=\scriptsize,right] at (axis cs:0.719,-0.391) {$x_-$};
\node[font=\scriptsize,below left] at (axis cs:0.550,0) {$x_c$};
\node[font=\scriptsize,actionblue] at (axis cs:0.390,0.166) {$\Gamma_+$};
\node[font=\scriptsize,actiongold] at (axis cs:0.513,-0.196) {$\Gamma_-$};
\node[font=\scriptsize,actionred] at (axis cs:0.276,0.720) {$C_+$};
\node[font=\scriptsize,actionred] at (axis cs:0.900,-0.748) {$C_-$};
\node[font=\scriptsize] at (axis cs:.12,.45) {$D_L$};
\node[font=\scriptsize] at (axis cs:1.05,.45) {$D_R$};
\end{axis}
\end{tikzpicture}
\hspace{.2cm}
\begin{tikzpicture}
\begin{axis}[scale only axis,width=4.0cm,height=4.0cm,xmin=-0.649,xmax=1.012,ymin=0,ymax=1.25,xlabel={$\Re E$},ylabel={actions},label style={font=\footnotesize},tick label style={font=\scriptsize},title style={font=\small},tick align=outside,axis line style={gray!75},scaled ticks=false,xtick={-.6,0,.6},ytick={0,.5,1},title={(c) $\Im E=0.1$}]
\addplot[black!18,line width=3.2pt] coordinates {(-0.6185,0.00000) (-0.5852,0.00290) (-0.5518,0.00825) (-0.5185,0.01522) (-0.4852,0.02354) (-0.4518,0.03304) (-0.4185,0.04363) (-0.3852,0.05522) (-0.3518,0.06776) (-0.3185,0.08120) (-0.2852,0.09552) (-0.2518,0.11069) (-0.2185,0.12668) (-0.1852,0.14347) (-0.1518,0.16105) (-0.1185,0.17940) (-0.0852,0.19851) (-0.0518,0.21838) (-0.0185,0.23898) (0.0148,0.25073) (0.0482,0.24794) (0.0815,0.24543) (0.1148,0.24320) (0.1482,0.24123) (0.1815,0.23952) (0.2148,0.23805) (0.2482,0.23683) (0.2815,0.23586) (0.3148,0.23512) (0.3482,0.23462) (0.3815,0.22155) (0.4148,0.20325) (0.4482,0.18550) (0.4815,0.16831) (0.5148,0.15171) (0.5482,0.13570) (0.5815,0.12032) (0.6148,0.10557) (0.6482,0.09149) (0.6815,0.07810) (0.7148,0.06545) (0.7482,0.05356) (0.7815,0.04250) (0.8148,0.03233) (0.8482,0.02313) (0.8815,0.01503) (0.9148,0.00818) (0.9482,0.00289) (0.9815,0.00000)};
\addplot[actiongreen,thick] coordinates {(-0.6185,0.41660) (-0.5852,0.40548) (-0.5518,0.39641) (-0.5185,0.38856) (-0.4852,0.38164) (-0.4518,0.37549) (-0.4185,0.37000) (-0.3852,0.36509) (-0.3518,0.36071) (-0.3185,0.35680) (-0.2852,0.35333) (-0.2518,0.35027) (-0.2185,0.34759) (-0.1852,0.34527) (-0.1518,0.34329) (-0.1185,0.34164) (-0.0852,0.34031) (-0.0518,0.33929) (-0.0185,0.33857) (0.0148,0.33814) (0.0482,0.33801) (0.0815,0.33816) (0.1148,0.33860) (0.1482,0.33933) (0.1815,0.34034) (0.2148,0.34165) (0.2482,0.34326) (0.2815,0.34517) (0.3148,0.34739) (0.3482,0.34994) (0.3815,0.35281) (0.4148,0.35603) (0.4482,0.35961) (0.4815,0.36356) (0.5148,0.36791) (0.5482,0.37267) (0.5815,0.37788) (0.6148,0.38354) (0.6482,0.38971) (0.6815,0.39641) (0.7148,0.40369) (0.7482,0.41160) (0.7815,0.42019) (0.8148,0.42953) (0.8482,0.43973) (0.8815,0.45090) (0.9148,0.46322) (0.9482,0.47701) (0.9815,0.49310)};
\addplot[actionred,thick] coordinates {(-0.6185,0.38943) (-0.5852,0.37431) (-0.5518,0.36134) (-0.5185,0.34972) (-0.4852,0.33914) (-0.4518,0.32944) (-0.4185,0.32051) (-0.3852,0.31226) (-0.3518,0.30463) (-0.3185,0.29756) (-0.2852,0.29102) (-0.2518,0.28496) (-0.2185,0.27935) (-0.1852,0.27417) (-0.1518,0.26938) (-0.1185,0.26497) (-0.0852,0.26092) (-0.0518,0.25721) (-0.0185,0.25382) (0.0148,0.25073) (0.0482,0.24794) (0.0815,0.24543) (0.1148,0.24320) (0.1482,0.24123) (0.1815,0.23952) (0.2148,0.23805) (0.2482,0.23683) (0.2815,0.23586) (0.3148,0.23512) (0.3482,0.23462) (0.3815,0.23436) (0.4148,0.23433) (0.4482,0.23455) (0.4815,0.23502) (0.5148,0.23575) (0.5482,0.23674) (0.5815,0.23801) (0.6148,0.23957) (0.6482,0.24144) (0.6815,0.24365) (0.7148,0.24621) (0.7482,0.24918) (0.7815,0.25258) (0.8148,0.25649) (0.8482,0.26097) (0.8815,0.26614) (0.9148,0.27216) (0.9482,0.27932) (0.9815,0.28846)};
\addplot[actionblue,thick] coordinates {(-0.6185,1.02975) (-0.5852,0.99084) (-0.5518,0.95416) (-0.5185,0.91888) (-0.4852,0.88473) (-0.4518,0.85154) (-0.4185,0.81920) (-0.3852,0.78765) (-0.3518,0.75683) (-0.3185,0.72669) (-0.2852,0.69719) (-0.2518,0.66831) (-0.2185,0.64003) (-0.1852,0.61232) (-0.1518,0.58518) (-0.1185,0.55858) (-0.0852,0.53251) (-0.0518,0.50698) (-0.0185,0.48196) (0.0148,0.45746) (0.0482,0.43346) (0.0815,0.40998) (0.1148,0.38700) (0.1482,0.36452) (0.1815,0.34255) (0.2148,0.32109) (0.2482,0.30014) (0.2815,0.27971) (0.3148,0.25979) (0.3482,0.24041) (0.3815,0.22155) (0.4148,0.20325) (0.4482,0.18550) (0.4815,0.16831) (0.5148,0.15171) (0.5482,0.13570) (0.5815,0.12032) (0.6148,0.10557) (0.6482,0.09149) (0.6815,0.07810) (0.7148,0.06545) (0.7482,0.05356) (0.7815,0.04250) (0.8148,0.03233) (0.8482,0.02313) (0.8815,0.01503) (0.9148,0.00818) (0.9482,0.00289) (0.9815,0.00000)};
\addplot[actiongold,thick] coordinates {(-0.6185,0.00000) (-0.5852,0.00290) (-0.5518,0.00825) (-0.5185,0.01522) (-0.4852,0.02354) (-0.4518,0.03304) (-0.4185,0.04363) (-0.3852,0.05522) (-0.3518,0.06776) (-0.3185,0.08120) (-0.2852,0.09552) (-0.2518,0.11069) (-0.2185,0.12668) (-0.1852,0.14347) (-0.1518,0.16105) (-0.1185,0.17940) (-0.0852,0.19851) (-0.0518,0.21838) (-0.0185,0.23898) (0.0148,0.26032) (0.0482,0.28239) (0.0815,0.30518) (0.1148,0.32868) (0.1482,0.35290) (0.1815,0.37783) (0.2148,0.40346) (0.2482,0.42980) (0.2815,0.45683) (0.3148,0.48457) (0.3482,0.51300) (0.3815,0.54212) (0.4148,0.57194) (0.4482,0.60246) (0.4815,0.63368) (0.5148,0.66560) (0.5482,0.69822) (0.5815,0.73155) (0.6148,0.76559) (0.6482,0.80037) (0.6815,0.83589) (0.7148,0.87217) (0.7482,0.90923) (0.7815,0.94712) (0.8148,0.98587) (0.8482,1.02555) (0.8815,1.06627) (0.9148,1.10817) (0.9482,1.15153) (0.9815,1.19717)};
\addplot[only marks,mark=diamond*,mark size=1.7pt,actionred] coordinates {(0.0018,0.25190) (0.3585,0.23451)};
\draw[gray,densely dotted] (axis cs:0.2000,0) -- (axis cs:0.2000,1.25);
\node[font=\scriptsize,actiongreen] at (axis cs:0.815,0.520) {$A_L$};
\node[font=\scriptsize,actionred] at (axis cs:0.815,0.176) {$A_R$};
\node[font=\scriptsize,actionblue] at (axis cs:-0.358,0.879) {$U_+$};
\node[font=\scriptsize,actiongold] at (axis cs:0.668,0.909) {$U_-$};
\end{axis}
\end{tikzpicture}
\caption{The geometry of Section~\ref{sec:actions} for $V(z)=0.6z^2-\ii(1-2z)$, $\alpha=0.4$ and $E=0.2+0.1\ii$. (a) The tube $\Sweep$ (shaded), the centre curve $V([0,1])$ (dashed) and the horizontal segment $[E-2\alpha,E+2\alpha]$; the square is $V(x_c)$. (b) The inverse image of this segment under $V$: the real-momentum curves $\Gamma_+$ (blue) and $\Gamma_-$ (gold) join $x_c$ to the turning points $x_\pm$, and the cuts $C_\pm$ (red, dashed) continue the level curve $L$ beyond them. The dotted segments $[0,x_\pm]$ are the paths for $J_\pm$. (c) The four actions along the line $\Im E=0.1$ across $\Sweep$; the grey band is their minimum and the dotted line marks $\Re E=0.2$. $U_-$ vanishes on the left edge and $U_+$ on the right edge. The two diamonds, where $A_R=U_-$ and $A_R=U_+$, mark the points at which this line meets the support of $\mu$; they are also marked in (a).}\label{fig:geometry}
\end{figure}

Fix $E=\tau-\ii\sigma\in\Sweep^\circ$, let $g=g(\sigma)$, and put
$\mathfrak c(z)=(E-V(z))/(2\alpha)$ for $z\in\Omega$, so that
\eqref{eq:momentum} reads $\cos p=\mathfrak c$.  Let
\[
 L=\{z\in\Omega:f(z)=\sigma\},\qquad
 D_L=\{z\in\Omega:f(z)>\sigma\},
\]
$D_R=\{z\in\Omega:f(z)<\sigma\}$,
let $x_\pm=\beta(E\mp2\alpha)$ be the turning points, and let
\[
 \Gamma_+=\beta\bigl([\tau-2\alpha,g]-\ii\sigma\bigr),\qquad
 \Gamma_-=\beta\bigl([g,\tau+2\alpha]-\ii\sigma\bigr),
\]
\[
 \begin{aligned}
 C_+&=\{z\in\Omega:V(z)\in(-\infty,\tau-2\alpha]-\ii\sigma\},\\
 C_-&=\{z\in\Omega:V(z)\in[\tau+2\alpha,\infty)-\ii\sigma\}.
 \end{aligned}
\]
Let $D$ be the component of $\Omega\setminus(C_+\cup C_-)$ that contains
$[0,1]$.  The curves $\Gamma_\pm$ lie in $\Omega$, because
$[\tau-2\alpha,\tau+2\alpha]-\ii\sigma\subset\Sweep_{4\alpha}\subseteq
V(\Omega)$.  Since $V(x_c)=g-\ii\sigma$, both curves start at $x_c$, and they
end at $x_+$ and $x_-$, respectively.  Figure~\ref{fig:geometry}(a),(b) shows
these objects for the complex quadratic potential of Section~\ref{sec:numerics}.

\begin{lemma}\label{lem:geometry}
Assume \textup{(H1)}, \textup{(H2)} and $E\in\Sweep^\circ$.
\begin{enumerate}[(i)]
\item On every horizontal line, $x\mapsto f(x+\ii y)$ is strictly
decreasing.  Hence $L$ meets every horizontal line at most once, $D_L$ lies to
the left of $L$ and $D_R$ to its right.  Along $L$ the function $r$ is a
strictly decreasing function of $\Im z$.
\item $\Gamma_+$ is a graph over $0\le\Im z\le\Im x_+$, and $\Gamma_-$ a graph
over $\Im x_-\le\Im z\le0$; in particular $\Im x_+>0>\Im x_-$.  The cuts
satisfy $C_+\cup C_-\subset\{\Im z\ge\Im x_+\}\cup\{\Im z\le\Im x_-\}$, and
the points of $L$ just above $x_+$ (just below $x_-$) belong to $C_+$ (to
$C_-$).  The turning points in $\Omega$, that is, the zeros of $1-\mathfrak c^2$, are
exactly $x_+$ and $x_-$.
\item $D$ is simply connected and contains
$(\Gamma_+\cup\Gamma_-)\setminus\{x_+,x_-\}$ and the half-open segments
$[0,x_\pm)$.  On $D$ there is a unique
holomorphic $p$ with $\cos p=\mathfrak c$ and $0<\Re p<\pi$.  It satisfies $\Im p<0$ on
$D\cap D_L$, $\Im p>0$ on $D\cap D_R$, and on $\Gamma_+\cup\Gamma_-$
\[
 p\bigl(\beta(t-\ii\sigma)\bigr)=\arccos\frac{\tau-t}{2\alpha}.
\]
As $z\to x_+$, $p(z)=\kappa_+(z-x_+)^{1/2}(1+O(z-x_+))$, and as $z\to x_-$,
$\pi-p(z)=\kappa_-(z-x_-)^{1/2}(1+O(z-x_-))$, with $\kappa_\pm\ne0$.  In
particular $\sin p$ has branch points of square-root type at $x_\pm$.
\item There is an $\varepsilon$, with $0<\varepsilon<\min\{Y-\Im x_+,Y+\Im x_-\}$,
that has the following property.  Let $y\in(\Im x_+,\Im x_++\varepsilon]$ be such
that $\ii y\in D_L$ and $1+\ii y\in D_R$; when $0<\Re x_+<1$ this holds for
every such $y$ once $\varepsilon$ is small.  Then the segment $[\ii y,1+\ii y]$ meets $L$ at one point,
which lies on $C_+$.  The analytic continuation $\tilde p$ of $p$
along the segment from its left end equals $p$ on $D_L$ and $-p$ on $D_R$, and
$\Im\tilde p<0$ on the segment.  The same holds for
$y\in[\Im x_--\varepsilon,\Im x_-)$, with $C_-$ in place of $C_+$ and
$\tilde p=2\pi-p$ on $D_R$.
\end{enumerate}
\end{lemma}

\begin{proof}
(i) We have $f_x=-\Im V'<0$, which proves the claims about horizontal
lines. Next, if $z_1,z_2\in L$ are distinct, then $d=(V(z_2)-
V(z_1))/(z_2-z_1)$ has $\Im d>0$ by convexity. Thus
$\Im(z_2-z_1)=(r(z_2)-r(z_1))\Im(1/d)$, which has the sign opposite to that
of $r(z_2)-r(z_1)$. This proves the strict monotonicity along $L$.

(ii) The first statements follow from (i), since $r$ takes the values $g$ at
$x_c$, $\tau-2\alpha<g$ at $x_+$ and $\tau+2\alpha>g$ at $x_-$.  The heights strictly
between $\Im x_-$ and $\Im x_+$ are occupied by $\Gamma_\pm$, on which
$\tau-2\alpha<r<\tau+2\alpha$; since $L$ meets each horizontal line at most
once, no point of $C_\pm$ lies at these heights.  The turning points are the
solutions of $V=E\mp2\alpha$, which are unique by univalence.

(iii) The cosine maps $\{0<\Re p<\pi\}$ biholomorphically onto
$\C\setminus((-\infty,-1]\cup[1,\infty))$, and $\mathfrak c(z)$ lies in
$(-\infty,-1]\cup[1,\infty)$ exactly when $z\in C_+\cup C_-$.  This gives $p$.
Since $\Im\mathfrak c=(f-\sigma)/(2\alpha)$ and $\Im\cos p=-\sin(\Re p)\sinh(\Im p)$,
the sign of $\Im p$ is that of $\sigma-f$.  On $\Gamma_\pm$ one has
$V=t-\ii\sigma$ and $\mathfrak c=(\tau-t)/(2\alpha)\in[-1,1]$.  Every component of
$C_+\cup C_-$ is an arc of a level curve of the harmonic function $f$, and
these level curves contain no critical points ($V'\ne0$).  In the $V$-plane each
component of the removed set is a portion of one of two closed rays.
Every such portion approaches the boundary of $V(\Omega)$ at one end at
least, or is unbounded; only an end at $E\mp2\alpha$ can lie in this
domain. The inverse image of each portion is therefore relatively closed in
$\Omega$ and leaves every compact subset of $\Omega$ at one end at least;
hence, on the Riemann sphere, it attaches to the connected set
$\hat\C\setminus\Omega$.  This holds for each of the possibly infinitely many
components.  Hence $\hat\C\setminus D$ is
connected, and $D$ is simply connected.

By (ii) and convexity, the segments
$[0,x_\pm)$ lie in $D$.  The expansions at $x_\pm$ follow from
$1-\cos p\sim p^2/2$ near $p=0$ and $1+\cos p\sim(\pi-p)^2/2$ near $p=\pi$,
together with $\mathfrak c'(x_\pm)=-V'(x_\pm)/(2\alpha)\ne0$.

(iv) The first statement follows from (i) and (ii), since near $x_+$ the
curve $L$ continues upwards in $C_+$.  On $C_+$ one has
$\mathfrak c\in[1,\infty)$, and the boundary values of the principal arccosine are
$\pm\ii\operatorname{arccosh}\mathfrak c$.  Since $\Im p<0$ on $D_L$, the value from the
left is $-\ii\operatorname{arccosh}\mathfrak c$ and the value from the right is
$+\ii\operatorname{arccosh}\mathfrak c$.  Hence the continuation across $C_+$ is $-p$,
which has negative imaginary part on $D_R$.  At the intersection with $C_+$,
$\Im\tilde p=-\operatorname{arccosh}\mathfrak c<0$, because the segment passes above
$x_+$, where $\mathfrak c>1$.  On $C_-$ one has $\mathfrak c\in(-\infty,-1]$, and the boundary
values are $\pi\mp\ii\operatorname{arccosh}|\mathfrak c|$ from the left and from the
right.  The continuation is therefore $2\pi-p$, whose imaginary part is $-\Im p<0$
on $D_R$.
\end{proof}

\begin{remark}\label{rem:H2}
For $V(z)=V_0+\gamma z$ with $\Im\gamma>0$ one may take $\Omega=\C$.
If $V=-\ii B$ for a real quadratic $B$ with $B'<0$ on $[0,1]$, take the
half-plane that is bounded by the vertical line through the real critical point
$z_0$ of $B$ and contains $[0,1]$. There $\Im V'>0$, and the image of this
half-plane under $B$ is the plane cut along a real ray that starts at
$B(z_0)$ and points away from $[B(1),B(0)]$. Thus $V(\Omega)$ contains $\Sweep_{4\alpha}$
for every $\alpha$, and $\Omega$ contains every required rectangle.

For the complex quadratic $V(z)=0.6z^2-\ii(1-2z)$, use
$\Omega=\{\Im z>-5/3\}$. Completing the square gives
$V(z)=\tfrac35(z+5\ii/3)^2+5/3-\ii$, so the missing ray is
$[5/3,\infty)-\ii$. The tube $\Sweep_{4\alpha}$ misses this ray for $\alpha\le0.4$: on
$\Im E=-1$ its right endpoint is $4\alpha\le1.6$. The inverse image of the
tube is compact in $\Omega$, and in fact $|\Im z|\le4/3$ on it. Indeed, let
$z=x+\ii y\in\Omega$ with $V(z)=V(s)+t$, $s\in[0,1]$, $|t|\le1.6$. The
imaginary parts give $x(1.2y+2)=2s$. If $y>1$, then $0\le x\le s\le1$ and
$\Re V(z)=0.6(x^2-y^2)-2y<-2$, outside the tube. If $y<-4/3$, then
$d=1.2y+2\in(0,0.4)$, so $x^2=4s^2/d^2\ge s^2$, and the real parts give
$t\ge-0.6y^2-2y>1.6$, because $(3y+4)(y+2)<0$. Thus $Y\le4/3<5/3$,
which proves the rectangle condition.

In general, choose a fixed thin convex neighbourhood of $[0,1]$ on which
$\Im V'>0$. The image of this neighbourhood contains a neighbourhood of the
compact curve $V([0,1])$, and $\beta'$ is bounded in a smaller neighbourhood.
Every point of $\Sweep_{4\alpha}$ lies within horizontal distance $4\alpha$ of
that curve, and integrating $\beta'$ along the horizontal segment joining them
gives $Y=O(\alpha)$. This proves (H2) for small $\alpha$.
\end{remark}

\subsection{The actions}
Let $I$ and $J_\pm$ be as in \eqref{eq:IJ}.  The integrals defining $J_\pm$
converge and, since $D$ is simply connected, do not depend on the path in it;
by Lemma~\ref{lem:geometry}(iii) one may integrate along the straight segments
$[0,x_\pm]$.
Let $\Theta_1,\dots,\Theta_4$ be the four phases introduced with \eqref{eq:phases0}, with
$\log\alpha\in\R$.

\begin{lemma}\label{lem:actions}
For $E\in\Sweep^\circ$,
\[
 \Im I=A_R-A_L,\qquad \Im J_\pm=U_\pm-A_L,\qquad \Im(J_\pm-I)=U_\pm-A_R,
\]
and consequently
\begin{equation}\label{eq:realparts}
\begin{aligned}
 \Re\Theta_1&=\Pot{\M}-2A_R,&\qquad \Re\Theta_2&=\Pot{\M}-2A_L,\\
 \Re\Theta_3&=\Pot{\M}-2U_+,&\qquad \Re\Theta_4&=\Pot{\M}-2U_- .
\end{aligned}
\end{equation}
\end{lemma}

\begin{proof}
The first identity follows from Lemma~\ref{lem:momentum}.  For $J_+$, integrate
along $[0,x_c]$ and then along $\Gamma_+$.  The first part contributes $-A_L$.  On
$\Gamma_+$ we have $z=\beta(t-\ii\sigma)$, with $t$ decreasing from $g$ to
$\tau-2\alpha$; there the momentum is real, and $dy=\Im\beta'\,dt$ with
$\Im\beta'=-\Im V'(\beta)/|V'(\beta)|^2<0$.  So $\Im(p\,dz)=p\,dy=p\,|\Im\beta'|\,|dt|$, and by
Lemma~\ref{lem:geometry}(iii) this part contributes $U_+$.  For $J_-$,
integrate $p-\pi$ along $[0,x_c]$ and $\Gamma_-$; the first part again gives
$-A_L$, because $\pi$ is real.  On $\Gamma_-$, $t$ increases from $g$ to
$\tau+2\alpha$ and $dy<0$, so
$\Im((p-\pi)dz)=(\pi-p)|dy|=\arccos\frac{t-\tau}{2\alpha}|\Im\beta'|\,dt$,
and this part contributes $U_-$.  The identities for $J_\pm-I$ follow by
subtraction. Finally, $\Re\Theta_3=\log\alpha+\Im I-2\Im J_+$, and similarly
for $\Theta_4$; now \eqref{eq:PotM} gives \eqref{eq:realparts}.
\end{proof}

Figure~\ref{fig:geometry}(c) shows the four actions along one horizontal
line across $\Sweep$.  Each of $U_\pm$ vanishes on one edge, and the two
changes of the minimizing action are the points where this line meets the
support of $\mu$.

\subsection{Derivatives and densities}
Differentiating \eqref{eq:momentum} gives
$\partial_E p=-1/(2\alpha\sin p)=-\ii/(\rho_+-\rho_-)$, hence
\begin{equation}\label{eq:derivatives}
 I'(E)=-\ii\int_0^1\frac{dx}{\rho_+(x)-\rho_-(x)},\qquad
 J_\pm'(E)=-\ii\int_0^{x_\pm}\frac{dz}{\rho_+(z)-\rho_-(z)} ,
\end{equation}
the second pair of integrals being taken along paths in $D$.  The moving
endpoints contribute
nothing, because the integrands of $J_+$ and $J_-$ vanish at $x_+$ and $x_-$,
respectively; the integrands in \eqref{eq:derivatives} have integrable
singularities $(z-x_\pm)^{-1/2}$.  In particular, $\Theta_1,\dots,\Theta_4$ are
holomorphic on $\Sweep^\circ$.  Note also that
$\Theta_1'=\int_0^1\frac{dx}{\rho_+-\rho_-}$ is the Cauchy transform of $\M$,
continued analytically into $\Sweep$ from the region $\sigma<\sigma_1$; this
follows from Lemma~\ref{lem:PotM}.

\begin{corollary}\label{cor:curves}
Granting the Main Theorem, on an arc where exactly two actions attain the
minimum, $\mu$ has density $(2\pi)^{-1}|\Theta_i'-\Theta_j'|$ with respect to
arclength, where $\Theta_i,\Theta_j$ are the corresponding phases.  Explicitly:
\[
\begin{array}{lll}
 A_L=A_R:\ \tfrac1\pi|I'|, &
 A_R=U_\pm:\ \tfrac1\pi|I'-J_\pm'|, &
 A_L=U_\pm:\ \tfrac1\pi|J_\pm'|,\\[3pt]
 U_+=U_-:\ \tfrac1\pi|J_+'-J_-'|. & &
\end{array}
\]
\end{corollary}

\begin{proof}
The densities follow from Corollary~\ref{cor:charts} and the identities
$\Theta_1-\Theta_2=2\ii I$,
$\Theta_1-\Theta_{3,4}=2\ii(I-J_\pm)$, $\Theta_2-\Theta_{3,4}=-2\ii J_\pm$ and
$\Theta_3-\Theta_4=2\ii(J_+-J_-)$.  Each $\Theta_i'-\Theta_j'$ is holomorphic
and not identically zero, since otherwise $\Re(\Theta_i-\Theta_j)$ would be
constant, hence zero on the arc, contradicting the non-resonance shown in the
proof of the Main Theorem.  So its zeros are isolated; they carry no mass,
because $\Pot{\mu}$ is continuous.
\end{proof}

The next lemma shows that the support is a finite union of arcs, even near
the corners of $\Sweep$, where two actions vanish simultaneously.

\begin{lemma}\label{lem:finite}
Assume \textup{(H1)} and \textup{(H2)}, and let $T$ be the set of
$E\in\Sweep^\circ$ at which the minimum of the four actions is attained
at least twice.  Then $T$ is a finite union of real-analytic arcs and points.
Near the corner $V(0)+2\alpha$, $T$ coincides with $\{A_L=U_+\}$ and
is a single analytic arc ending at the corner; the same holds at the other
three corners $V(0)-2\alpha$, $V(1)\pm2\alpha$, with the pairs $(A_L,U_-)$,
$(A_R,U_\pm)$, respectively.
\end{lemma}

\begin{proof}
\emph{Harmonicity.}  By Lemma~\ref{lem:actions} the six differences of the
actions are, up to sign, $\Im I$, $\Im J_\pm$, $\Im(J_\pm-I)$ and
$\Im(J_+-J_-)$.  These are harmonic on $\Sweep^\circ$, and none of them is
identically zero, by the non-resonance argument in the proof of the Main
Theorem.

\emph{Edges.}  The actions are continuous on $\Sweep$.  Moreover, $A_L>0$ where
$\sigma<\sigma_0$, $A_R>0$ where $\sigma>\sigma_1$, $U_+>0$ where
$\tau<g+2\alpha$, and $U_->0$ where $\tau>g-2\alpha$.  Uniformly on $\Sweep$,
\begin{itemize}
\item $A_L\le x_c\max|\Im p|\to0$ as $\sigma\to\sigma_0$, and likewise
$A_R\to0$ as $\sigma\to\sigma_1$;
\item $U_+=O((g+2\alpha-\tau)^{3/2})$ and $U_-=O((\tau-g+2\alpha)^{3/2})$ at
the right and left boundaries.
\end{itemize}
Hence, for every $\varepsilon>0$ there is $\delta>0$ with the following
property: at the points of $\Sweep$ within $\delta$ of $\partial\Sweep$ and at
distance at least $\varepsilon$ from the four corners, the minimum is attained
only once.

\emph{Corners.}  Consider the corner $E_0=V(0)+2\alpha$.  Near it
$A_R$ and $U_-$ are bounded below while $A_L,U_+\to0$, so $T$ lies in
$\{\Im J_+=0\}$ there.  Put $\xi=x_+(E)$.  The map $E\mapsto\xi$
is conformal near $E_0$, since $\beta'(V(0))\ne0$, and $\xi(E_0)=0$.
It maps $\Sweep$ near $E_0$ onto a region with a corner at $0$ of some
opening $\theta\in(0,\pi)$.  Since $V(\xi)=E-2\alpha$, the
substitution $z=\xi t$ gives
\[
 J_+=\xi\int_0^1p(\xi t)\,dt,\qquad
 \cos p(\xi t)=1+\xi(1-t)G(\xi,t),
\]
where $G(\xi,t)=(2\alpha)^{-1}\int_0^1V'\bigl(\xi(t+s(1-t))\bigr)ds$ and
$G(0,t)=V'(0)/(2\alpha)\ne0$; the segment $[0,\xi)$ lies in $D$ by
Lemma~\ref{lem:geometry}(iii).  Write
$\arccos(1-w)=(2w)^{1/2}\varphi(w)$, with $\varphi$ holomorphic near $0$ and
$\varphi(0)=1$.  For suitable branches this gives
\[
 J_+=\xi^{3/2}F(\xi),\qquad
 F(\xi)=\int_0^1(1-t)^{1/2}\bigl(-2G(\xi,t)\bigr)^{1/2}
 \varphi\bigl(-\xi(1-t)G(\xi,t)\bigr)\,dt ,
\]
with $F$ holomorphic near $0$ and $F(0)=\frac23(-2G(0,0))^{1/2}\ne0$.  Hence
the function $\eta=\xi^{1/2}F(\xi)^{1/3}$, with a branch of $\xi^{1/2}$ on the
corner region, maps $\Sweep$ near $E_0$ conformally onto a region with a
corner of opening $\theta/2<\pi/2$ at $0$, bounded by two analytic arcs, and
$J_+=\eta^3$.  The
set $\{\Im\eta^3=0\}$ consists of six rays from $0$ at mutual angles $\pi/3$.
We determine which of them enter the region.  Write $a=\arg V'(0)\in(0,\pi)$.
On the edge $\sigma=\sigma_0$ one has $x_c=0$, $J_+=\int_{\Gamma_+}p\,dz$ with
$p\ge0$, and $\Gamma_+$ leaves $0$ in the direction of $-\beta'(V(0))$, so
$\arg J_+\to\pi-a$ at the corner.  On the side $\tau=g+2\alpha$ one has
$x_+=x_c\in[0,1]$ and $p(x)\sim(-V'(0)(x_c-x)/\alpha)^{1/2}$ with $\Im p<0$
for $x<x_c$, so $\arg J_+\to-(\pi-a)/2$.  Thus $\theta=\pi-a$, and $3\arg\eta$
ranges over an interval of length $3\theta/2<3\pi/2$ whose end points
$-(\pi-a)/2$ and $\pi-a$ lie in $(-\pi,0)$ and $(0,\pi)$, respectively, so
the interval contains $0$ but not $\pm\pi$.  Hence exactly one of the six rays enters the
region; it stays away from the two boundary arcs near $0$ and leaves the
corner at the angle $\theta/3$ from the side $V([0,1])+2\alpha$.  This is
consistent with the signs on the edges: $\Im J_+=U_+-A_L$ is positive on
$\sigma=\sigma_0$, where $A_L=0<U_+$, and negative on $\tau=g+2\alpha$, where
$U_+=0<A_L$.  Since $A_R,U_-$ are bounded below near $E_0$, $T$ coincides
there with $\{\Im J_+=0\}$, which is a single analytic arc ending at $E_0$.
At the other three corners the same computation applies to
$J_-=\int_0^{x_-}(p-\pi)\,dz$, to
$J_+-I=-\int_{x_+}^1p\,dz$ and to $J_--I=-\int_{x_-}^1(p-\pi)\,dz-\pi$, whose
integrands vanish at $x_\pm$, with
$\xi=x_-$, $1-x_\pm$, respectively.

\emph{Conclusion.}  The closure of $T$ meets $\partial\Sweep$ only at the
four corners.  Removing small open neighbourhoods of the corners from it
therefore leaves a compact subset of $\Sweep^\circ$.  At an interior tie point,
a local nodal chart gives finitely many analytic half-branches.  Each other action difference restricted to a
half-branch is real-analytic at its initial point, so it either vanishes
identically or has constant strict sign sufficiently near that point.
Consequently each sufficiently short half-branch lies wholly in $T$ or is
disjoint from $T$.  A finite subcover of the compact part, together with the
four shortened corner arcs, expresses $T$ as a finite union of analytic arcs
and points.  The parametrizations may be chosen regular and injective, with
analytic extensions across their endpoints.
\end{proof}

\section{Complex WKB tools}\label{sec:WKB}

\subsection{The equation}
Put $h=1/N$ and $P_k=\alpha^k\varphi_k$.  Then \eqref{eq:recurrence}
becomes
\begin{equation}\label{eq:schrodinger}
 \varphi_{k+1}+\varphi_{k-1}+v(kh)\varphi_k=0,\qquad
 v(z)=-\frac{E-V(z)}{\alpha},\qquad \varphi_0=1,\ \varphi_{-1}=0,
\end{equation}
and $p_N(E)=\alpha^N\varphi_N$.  This recurrence is the lattice restriction of
the analytic difference equation
\begin{equation}\label{eq:analytic}
 \psi(z+h)+\psi(z-h)+v(z)\psi(z)=0 ,
\end{equation}
whose complex momentum, defined by $2\cos p+v=0$, is exactly the frozen
momentum of \eqref{eq:momentum} \cite{BF95,FS15,FK19,FS22,Fed23}.  A solution
that is holomorphic on a vertical band of width greater than $2h$ extends, by
the recurrence, to the horizontal strip through that band, as far as $v$
remains holomorphic.

\subsection{Canonical curves and WKB solutions}
Let $q$ be holomorphic near a $C^1$ curve $\gamma(y)=x(y)+\ii y$,
$y_1\le y\le y_2$, and suppose that $\sin q\ne0$ there.  The curve is \emph{canonical} for
$q$ with margin $\eta>0$ if
\[
 \eta\le\Re q+x'\Im q\le\pi-\eta ,
\]
that is, $\Im(q\dot\gamma)\ge\eta$ and $\Im((q-\pi)\dot\gamma)\le-\eta$
\cite[\S5]{FK19}, \cite[Def.~2.3]{KO26}.

\begin{theorem}[{\cite[Thm.~6]{KO26}}; see also \cite{FS15}, {\cite[Thm.~1.2]{FK19}}]
\label{thm:WKB}
Let $K_0\Subset K_1$ be bounded simply connected domains whose
intersections with horizontal lines are intervals.  Let $q$ be a
branch of the momentum holomorphic near $\overline{K_1}$, with
$|\sin q|\ge c_0>0$ there, and fix a branch of $(\sin q)^{1/2}$.  Suppose that
$K_1$ is swept by a compact $C^1$ family of curves with two common endpoints,
all canonical with margin $\eta$.  Let $z_0\in K_0$.  Then for small $h$ there
are solutions $\psi_\pm$ of \eqref{eq:analytic}, holomorphic in $K_1$, with
\[
 \psi_\pm(z)=(\sin q(z))^{-1/2}\exp\Bigl(\pm\frac\ii h\int_{z_0}^zq\Bigr)
 \bigl(1+O(h)\bigr)
\]
uniformly on $K_0$, and $W[\psi_+,\psi_-]=2\ii(1+O(h))$ on $K_0$.  The
estimates are locally uniform in $E$ when the domains, the branch bounds and
the canonical margins are uniform.
\end{theorem}

\begin{lemma}\label{lem:segment}
Let $q$ be holomorphic near a vertical segment
$\{x_0+\ii y:y_1\le y\le y_2\}$, with $y_1<y_2$, on which $0<\Re q<\pi$.  Then for every
$\varepsilon>0$ there is $w>0$ such that Theorem~\ref{thm:WKB} applies with
$K_0\supseteq[x_0-w,x_0+w]+\ii[y_1+\varepsilon,y_2-\varepsilon]$.
\end{lemma}

\begin{proof}
For $|s|\le1$ let
$\gamma_s(y)=x_0+sw_1(y-y_1)(y_2-y)+\ii y$.  For small $w_1$
these curves stay in the region where $q$ is holomorphic, and by compactness
they are canonical with a uniform margin, since $\Re q+x'\Im q=\Re q+O(w_1)$.
They have common endpoints and sweep a bounded convex lens.  There is an open
rectangle that has compact closure in this lens and contains the required
trimmed rectangle; this rectangle and the lens give $K_0$ and $K_1$.
\end{proof}

We call the solutions obtained in this way the \emph{WKB solutions normalized
at $z_0$ along the segment}; they depend on the choice of the segment.

\subsection{Horizontal continuation}

\begin{theorem}[{\cite[Thm.~7]{KO26}}; see also {\cite[Thm.~3.1]{FK21}}, \cite{FK19}]\label{thm:horizontal}
Let $q$ and a branch of $(\sin q)^{1/2}$ be holomorphic on a neighbourhood of
a horizontal segment $[w_0,w_1]$, with $\Re w_1>\Re w_0$ and $\Im q\le-\eta<0$
there.  Suppose that a solution $\psi$ of \eqref{eq:analytic} satisfies
$\psi=(\sin q)^{-1/2}e^{\frac\ii hS}(1+O(h))$, with $S'=q$, on a neighbourhood
of $w_0$.  Then the same representation holds, uniformly, on a neighbourhood
of $[w_0,w_1]$.
\end{theorem}

Applied to $q$ or to $-q$, and after reflection, the theorem says that a
representation $(\sin q)^{-1/2}e^{\frac\ii hS}$ with $S'=q$ can be continued
in the direction in which $|e^{\frac\ii hS}|$ grows, that is, to the right
where $\Im q<0$ and to the left where $\Im q>0$.

Only fixed-energy versions of Theorems~\ref{thm:WKB} and~\ref{thm:horizontal}
are used below; all constants may depend on $E$.  For fixed $E$,
Theorem~\ref{thm:WKB} is the scalar case of \cite[Thm.~2.2 and (2.14)]{FS22}
(see also \cite{FS15}), and Theorem~\ref{thm:horizontal} with $1+o(1)$ in
place of $1+O(h)$ is the continuation principle \cite[Thm.~3.1]{FK21}.
These published versions suffice for the Main Theorem and for the
weak-limit statements of Section~\ref{sec:general}.  Only the rates
$O(1/N)$, as in Theorem~\ref{thm:sharp} and Proposition~\ref{prop:exterior},
require the $O(h)$ form of the continuation proved in \cite{KO26}.

\subsection{Wronskians}
For solutions $f,g$ of \eqref{eq:analytic} put
$W[f,g](z)=f(z+h)g(z)-f(z)g(z+h)$.  This Wronskian is $h$-periodic in $z$,
and hence constant on the lattice $h\Z$.  If $W[f,g]\ne0$ on $h\Z$, then every solution
$\varphi$ of \eqref{eq:schrodinger} satisfies
$\varphi=\frac{W[\varphi,g]}{W[f,g]}f+\frac{W[f,\varphi]}{W[f,g]}g$ on $h\Z$.

\begin{lemma}\label{lem:wronskian}
Suppose that near $[z,z+h]$
\[
 f=F\,e^{\frac\ii hS_f}(1+O(h)),\qquad g=G\,e^{\frac\ii hS_g}(1+O(h)),
\]
with $F,G,S_f,S_g$ holomorphic on a fixed neighbourhood, $F,G$ nonvanishing there, and $S_f'=q_f$,
$S_g'=q_g$.  Then
\[
 W[g,f](z)=F(z)G(z)\,e^{\frac\ii h(S_f(z)+S_g(z))}
 \bigl(e^{\ii q_g(z)}-e^{\ii q_f(z)}+O(h)\bigr).
\]
In particular, if $e^{\ii q_f}=e^{-\ii q_g}$ and $\sin q_g\ne0$, the bracket is
$2\ii\sin q_g\,(1+O(h))$, and if $q_g=q_f$ it is $O(h)$.
\end{lemma}

\begin{proof}
By Taylor expansion,
$e^{\frac\ii hS(z+h)}=e^{\frac\ii hS(z)}e^{\ii S'(z)}(1+O(h))$ and
$F(z+h)=F(z)(1+O(h))$.
\end{proof}

\begin{lemma}\label{lem:lattice}
Let $h>0$ and let $F$ be holomorphic and $h$-periodic on $\{y_-\le\Im z\le y_+\}$, with
$y_-<0<y_+$, and let $F_n=h^{-1}\int_0^hF(x+\ii y)e^{-2\pi\ii n(x+\ii y)/h}dx$
be its Fourier coefficients; they do not depend on $y$.  On $h\Z$ one has
$F=F^{\le0}+F^{\ge1}$ with $F^{\le0}=\sum_{n\le0}F_n$ and
$F^{\ge1}=\sum_{n\ge1}F_n$.  Let $0<\delta\le1$ and $h\le\min\{y_+,|y_-|\}$.
\begin{enumerate}[(i)]
\item If $|F-F_+^*|\le\delta|F_+^*|$ on the line $\Im z=y_+$ for a constant
$F_+^*$, then $|F^{\le0}-F_+^*|\le3\delta|F_+^*|$.
\item If $|F|\le M_+$ on the line $\Im z=y_+$, then $|F^{\le0}|\le2M_+$.
\item If $|F-F_-^*e^{2\pi\ii z/h}|\le\delta|F_-^*e^{2\pi\ii z/h}|$ on the line
$\Im z=y_-$ for a constant $F_-^*$, then $|F^{\ge1}-F_-^*|\le3\delta|F_-^*|$.
\item If $|F|\le M_-$ on the line $\Im z=y_-$, then
$|F^{\ge1}|\le2M_-e^{-2\pi|y_-|/h}$.
\end{enumerate}
\end{lemma}

\begin{proof}
(i) Computing the coefficients on $\Im z=y_+$, we obtain
$|F_0-F_+^*|\le\delta|F_+^*|$ and, for $n<0$,
$|F_n|\le\delta|F_+^*|e^{-2\pi|n|y_+/h}$, since the constant $F_+^*$ has no
other Fourier coefficients.  The proof of (iii) is the same, on $\Im z=y_-$,
for $F_1$ and $F_n$, $n\ge2$, after multiplication by $e^{-2\pi\ii z/h}$.
Parts (ii) and (iv) follow from $|F_n|\le M_\pm e^{2\pi ny_\pm/h}$.  In each
case it remains to sum the geometric series.  Finally, continuity and
periodicity give finite bounds on both boundary lines, so the same
coefficient estimates show that the Fourier series converges absolutely on
the real axis, and Fourier inversion gives the asserted identity on $h\Z$.
\end{proof}

\section{Proofs}\label{sec:proof}

Fix $E\in\Sweep^\circ$ and let $h=1/N$; constants may depend on
$E$.  We use the domain $D$, the branch $p$ and the curves
$\Gamma_\pm$ of Lemma~\ref{lem:geometry}, together with a fixed branch of
$(\sin p)^{1/2}$ on $D$, which exists since $\sin p\ne0$ on $D$ and $D$ is
simply connected.  Put
\[
 S_0(z)=\int_0^zp,\qquad S_1(z)=\int_1^zp,\qquad S_c(z)=\int_{x_c}^zp ,
\]
with integration paths in $D$.

\subsection{Three solutions}
Short vertical segments through $0$, $x_c$ and $1$ lie in $D$ by
Lemma~\ref{lem:geometry}(ii), since $C_\pm\subset\{\Im z\ge\Im x_+\}\cup
\{\Im z\le\Im x_-\}$; moreover $0<\Re p<\pi$ on them.  Let $Y_0^\pm>0$ be heights
such that the vertical segments
\[
 V_0=\{\ii y:-Y_0^-\le y\le Y_0^+\},\qquad V_1=1+V_0
\]
lie in $D$; these heights will be fixed in the proof of
Proposition~\ref{prop:switch}.  Let
\begin{itemize}
\item $\psi_0$ be the plus WKB solution normalized at $0$ along $V_0$;
\item $\psi_1$ be the plus WKB solution normalized at $1$ along $V_1$;
\item $u$ be the minus WKB solution normalized at $x_c$ along a short vertical
segment through $x_c$.
\end{itemize}
Thus, near the respective points,
\[
 \psi_0=\frac{e^{\frac\ii hS_0}}{(\sin p)^{1/2}}(1+O(h)),\quad
 \psi_1=\frac{e^{\frac\ii hS_1}}{(\sin p)^{1/2}}(1+O(h)),\quad
 u=\frac{e^{-\frac\ii hS_c}}{(\sin p)^{1/2}}(1+O(h)).
\]

\begin{lemma}\label{lem:realaxis}
For every small $\delta>0$ the representation of $\psi_0$ holds on a
neighbourhood of $[0,x_c-\delta]$, that of $u$ on neighbourhoods of $[0,x_c]$
and of $[x_c,1]$, and that of $\psi_1$ on a neighbourhood of $[x_c+\delta,1]$.
On $h\Z$,
\[
 W[\psi_0,u]=2\ii e^{\frac\ii hS_0(x_c)}(1+O(h)),\qquad
 W[\psi_1,u]=2\ii e^{\frac\ii hS_1(x_c)}(1+O(h)).
\]
\end{lemma}

\begin{proof}
By Lemma~\ref{lem:momentum}, $\Im p<0$ on $[0,x_c)$ and $\Im p>0$ on
$(x_c,1]$.  By Theorem~\ref{thm:horizontal}, $\psi_0$ continues to the right
on $[0,x_c-\delta]$, and $\psi_1$ continues to the left on $[x_c+\delta,1]$.
The representation of $u$ holds on a box $[x_c-w,x_c+w]+\ii[-w,w]$ by
Lemma~\ref{lem:segment}, and from $x_c\mp w/2$ it continues to the left on
$[0,x_c]$ and to the right on $[x_c,1]$.  Since the Wronskians are constant
on $h\Z$, we may evaluate them by Lemma~\ref{lem:wronskian} at lattice points
of $[x_c-2\delta,x_c-\delta]$ and of $[x_c+\delta,x_c+2\delta]$,
respectively, noting that $S_0-S_c=S_0(x_c)$ and $S_1-S_c=S_1(x_c)$.
\end{proof}

\subsection{Reduction to one Wronskian}
By Lemma~\ref{lem:realaxis}, $\psi_1$ and $u$ form a basis of solutions on
$h\Z$.  Write
$\psi_0=a\psi_1+cu$ on $h\Z$, and call
\[
 \W=W[\psi_1,\psi_0]\quad\text{(a constant on $h\Z$)}
\]
the \emph{switching Wronskian}.  Then $a=W[\psi_0,u]/W[\psi_1,u]$ and
$c=\W/W[\psi_1,u]$.

\begin{proposition}\label{prop:reduction}
With $A,B$ as in Theorem~\ref{thm:sharp},
\begin{equation}\label{eq:reduction}
 \frac{(\sin p(1))^{1/2}p_N(E)}{\alpha^N}
 =Ae^{\frac\ii hI}(1+O(h))+Be^{-\frac\ii hI}(1+O(h))
 +\frac{A}{2\ii}\,\W\,(1+O(h)).
\end{equation}
\end{proposition}

\begin{proof}
On $h\Z$ write $\varphi=A'\psi_0+B'u$.  The conditions $\varphi_0=1$,
$\varphi_{-1}=0$ form a linear system, whose coefficients are given by the WKB
forms at $0$ and $-h$: $\psi_0(kh)=e^{\ii kp(0)}(\sin p(0))^{-1/2}(1+O(h))$
and $u(kh)=\Lambda_0e^{-\ii kp(0)}(\sin p(0))^{-1/2}(1+O(h))$ for $k=0,-1$,
where $\Lambda_0=e^{\frac\ii hS_0(x_c)}$.  Solving gives $A'=A(1+O(h))$ and
$B'=B\Lambda_0^{-1}(1+O(h))$.  Then
\[
 \varphi_N=A'a\,\psi_1(1)+(A'c+B')\,u(1).
\]
By Lemma~\ref{lem:realaxis}, $a=e^{\frac\ii h(S_0(x_c)-S_1(x_c))}(1+O(h))
=e^{\frac\ii hI}(1+O(h))$ and $\psi_1(1)=(\sin p(1))^{-1/2}(1+O(h))$.
Moreover $u(1)=e^{-\frac\ii hS_c(1)}(\sin p(1))^{-1/2}(1+O(h))$ and
$S_1(x_c)=-S_c(1)$.  Hence $u(1)/W[\psi_1,u]=(2\ii)^{-1}(\sin p(1))^{-1/2}
(1+O(h))$ and $B'u(1)=Be^{-\frac\ii hI}(\sin p(1))^{-1/2}(1+O(h))$.  It
remains to use $p_N=\alpha^N\varphi_N$.
\end{proof}

The first two terms of \eqref{eq:reduction} are already in final form, and
everything has been reduced to the single number $\W$.  The factor $1+O(h)$
multiplying $\W$ is one common number, coming from $A'$ and
$u(1)/W[\psi_1,u]$.  Hence any splitting of $\W$, such as
$\W=\W^{\le0}+\W^{\ge1}$ below, yields independent relative errors for its
parts.  Both $\psi_0$ and $\psi_1$ are plus solutions, so even if their
representations held at a common point near $x_c$,
Lemma~\ref{lem:wronskian} would only give the upper bound
$|\W|\le Ch|e^{\frac\ii h(S_0+S_1)}|$.  The true size is found in the
complex domain: the Fourier modes $n\le0$ of $\W$ are determined above the
real axis, and the modes $n\ge1$ below it.

\subsection{The switching Wronskian}
Let $\mathcal E\subset\Sweep^\circ$ be the set of $E$ with
$\Re x_+(E)\in\{0,1\}$ or $\Re x_-(E)\in\{0,1\}$.  It is a union
of level sets of the non-constant harmonic functions
$\Re\beta(E\mp2\alpha)$, and hence a closed set of zero area.  For
$E\notin \mathcal E$, exactly one of the following holds for each sign:
\begin{enumerate}
\item[(a$_\pm$)] $\Gamma_\pm$ lies in the strip $\{0<\Re z<1\}$;
\item[(b$_\pm$)] $\Gamma_\pm$ touches or crosses $\{\Re z=0\}\cup\{\Re z=1\}$,
and does so for the first time (starting from $x_c$) at a height $\pm y_*^\pm$ with
$0<y_*^\pm<|\Im x_\pm|$.
\end{enumerate}

\begin{proposition}\label{prop:switch}
Let $E\in\Sweep^\circ\setminus \mathcal E$, and let $\W=\W^{\le0}+\W^{\ge1}$ be
the splitting of Lemma~\ref{lem:lattice}.
\begin{enumerate}
\item[(a$_+$)] In case \textup{(a$_+$)},
$\W^{\le0}=-2e^{\frac\ii h(2J_+-I)}(1+O(h))$.
\item[(b$_+$)] In case \textup{(b$_+$)} there is $c>0$ such that
$|\W^{\le0}|\le e^{-c/h}\max\{|e^{\frac\ii hI}|,|e^{-\frac\ii hI}|\}$ and
$U_+\ge\min\{A_L,A_R\}+c$.
\item[(a$_-$)] In case \textup{(a$_-$)},
$\W^{\ge1}=2e^{\frac\ii h(2J_--I)}(1+O(h))$.
\item[(b$_-$)] In case \textup{(b$_-$)} there is $c>0$ such that
$|\W^{\ge1}|\le e^{-c/h}\max\{|e^{\frac\ii hI}|,|e^{-\frac\ii hI}|\}$ and
$U_-\ge\min\{A_L,A_R\}+c$.
\end{enumerate}
\end{proposition}

\begin{proof}
\emph{Choice of the segments.}  In case (a$_+$), every horizontal line
$\Im z=y$ with $0\le y\le\Im x_+$ meets $L$ at a point of the open strip, by
Lemma~\ref{lem:geometry}(i),(ii); by continuity the same holds for
$0\le y\le Y_1^+$, for some $Y_1^+\in(\Im x_+,\Im x_++\varepsilon]$, with
$\varepsilon$ from Lemma~\ref{lem:geometry}(iv).  Put $Y_0^+=Y_1^+$.  In case
(b$_+$) put $Y_0^+=y_*^++\varepsilon$ with small $\varepsilon>0$.  Define
$Y_0^-$ analogously.  In all cases the segments $V_0$ and $V_1$ meet $L$ at
most in points of $\Gamma_\pm$, so they lie in $D$, and $0<\Re p<\pi$ on
them.  By Lemma~\ref{lem:segment} the representations of $\psi_0$ and $\psi_1$
hold on $[-w,w]+\ii[-Y_2^-,Y_2^+]$ and on $[1-w,1+w]+\ii[-Y_2^-,Y_2^+]$ for
some $w>0$ and $Y_2^\pm<Y_0^\pm$ as close to $Y_0^\pm$ as we wish.  By the
recurrence, $\psi_0$ and $\psi_1$ are holomorphic on
$[0,1]+\ii[-Y_2^-,Y_2^+]$, so $\W$ extends to an $h$-periodic holomorphic
function on the strip $-Y_2^-\le\Im z\le Y_2^+$.

\emph{Case (a$_+$): above the turning point.}  Fix
$Y'\in(\Im x_+,Y_2^+)$.  The segment $[\ii Y',1+\ii Y']$ has its left end in
$D_L$ and its right end in $D_R$.  By Lemma~\ref{lem:geometry}(iv), the
continuation $\tilde p$ of $p$ along it has $\Im\tilde p<0$.  This
continuation is regular, because $x_+$ is the only zero of $\sin p$ in the
upper half of $\Omega$ and $\Im x_+<Y'$.  Theorem~\ref{thm:horizontal}
therefore continues $\psi_0$ from $\ii Y'$ to $1+\ii Y'$: near $1+\ii Y'$,
\[
 \psi_0(z)=(\sin\tilde p(z))^{-1/2}e^{\frac\ii h\tilde S_0(z)}(1+O(h)),
 \qquad \tilde S_0(z)=\int_{[0,\ii Y']\cup[\ii Y',z]}\tilde p ,
\]
where $(\sin\tilde p)^{1/2}$ is continued along the same path.

We now compare this with $\psi_1$.  Choose the principal square root of
$\sin p$ on $D$, which is holomorphic since $\Re\sin p>0$.
Across $C_+$ its reciprocal continues as
$\sqrt{\ii}/\sqrt{\ii\sin\tilde p}$, with principal square roots: the
continued sine avoids the positive imaginary ray, so the denominator is
holomorphic and nonzero.  The continued amplitude equals $(\sin p)^{-1/2}$ on
$D_L$ and $\ii(\sin p)^{-1/2}$ on $D_R$.  Thus, near $1+\ii Y'$,
\[
 \tilde p=-p,\qquad (\sin\tilde p)^{1/2}=-\ii(\sin p)^{1/2},
\]
where $p$ and $(\sin p)^{1/2}$ are the branches on $D$, which are the ones
used in $\psi_1$.  Deforming the path of $\tilde S_0$ into one that runs in
$D$ to $x_+$ and then returns in $D$ on the other side of $C_+$, we obtain
$\tilde S_0(z)=J_+-\int_{x_+}^zp=2J_+-S_0(z)$, and since $S_1=S_0-I$,
\[
 \tilde S_0(z)+S_1(z)=2J_+-I .
\]
Lemma~\ref{lem:wronskian}, with $g=\psi_1$, $q_g=p$, $f=\psi_0$ and
$q_f=\tilde p=-p$, gives, near $1+\ii Y'$,
\[
 \W(z)=\frac{e^{\frac\ii h(2J_+-I)}\,2\ii\sin p\,(1+O(h))}
 {(\sin p)^{1/2}\cdot(-\ii)(\sin p)^{1/2}}
 =-2e^{\frac\ii h(2J_+-I)}(1+O(h)) .
\]
By periodicity this holds on the whole line $\Im z=Y'$, and
Lemma~\ref{lem:lattice}(i) gives (a$_+$).

\emph{Case (a$_-$): below the turning point.}  The argument is the same, with
the following changes.  Fix $Y''\in(|\Im x_-|,Y_2^-)$ and continue $\psi_0$
along $[-\ii Y'',1-\ii Y'']$.  By Lemma~\ref{lem:geometry}(iv) the continued
momentum is $\tilde p=2\pi-p$ on $D_R$, with $\Im\tilde p<0$.  Now the
continued sine avoids the negative imaginary ray, and the reciprocal square
root continues as $\sqrt{-\ii}/\sqrt{-\ii\sin\tilde p}$, which equals
$(\sin p)^{-1/2}$ on $D_L$ and $-\ii(\sin p)^{-1/2}$ on $D_R$.
Consequently, near $1-\ii Y''$,
\[
 \tilde p=2\pi-p,\qquad (\sin\tilde p)^{1/2}=\ii(\sin p)^{1/2}.
\]
Deforming the path through $x_-$ gives
$\tilde S_0(z)=\int_0^{x_-}p+\int_{x_-}^z(2\pi-p)=2J_-+2\pi z-S_0(z)$, hence
$\tilde S_0+S_1=2J_--I+2\pi z$.  Since $e^{\ii\tilde p}=e^{-\ii p}$,
Lemma~\ref{lem:wronskian} gives, near $1-\ii Y''$,
\[
 \W(z)=\frac{e^{\frac\ii h(2J_--I)}e^{2\pi\ii z/h}\,2\ii\sin p\,(1+O(h))}
 {(\sin p)^{1/2}\cdot\ii(\sin p)^{1/2}}
 =2e^{\frac\ii h(2J_--I)}e^{2\pi\ii z/h}(1+O(h)) ,
\]
and Lemma~\ref{lem:lattice}(iii) gives (a$_-$).

\emph{Case (b$_+$).}  Suppose that $\Gamma_+$ first meets the boundary of the
strip at $\ii y_*$, where $y_*=y_*^+$; the case $1+\ii y_*$ is symmetric under
exchanging $0$ and $1$, $\psi_0$ and $\psi_1$, and $A_L$ and $A_R$.  Then
$f(\ii y_*)=\sigma$, so $(0,1]+\ii y_*\subset D_R$.  By
Theorem~\ref{thm:horizontal} ($\Im p>0$ on $D_R$), $\psi_1$ continues to the
left along $[\delta+\ii y_*,1+\ii y_*]$ for any $\delta\in(0,w)$.  On the
interval $(\delta,w)+\ii y_*$ both $\psi_0$ and $\psi_1$ are plus solutions
with the same branch, so Lemma~\ref{lem:wronskian} gives
\[
 |\W(z)|\le Ch\,\bigl|e^{\frac\ii h(S_0(z)+S_1(z))}\bigr|
 =Ch\exp\Bigl(\frac1h\bigl(\Im I-2\Im S_0(z)\bigr)\Bigr).
\]
For $z=x+\ii y_*$ with $0<x<w$,
\[
 \Im S_0(z)=\int_0^{y_*}\Re p(\ii t)\,dt+\int_0^x\Im p(t+\ii y_*)\,dt
 \ge c:=\int_0^{y_*}\Re p(\ii t)\,dt>0 .
\]
Since $w-\delta>h$ for small $h$, the interval $(\delta,w)+\ii y_*$ contains a
full period, and $y_*<Y_2^+$ because $Y_0^+=y_*+\varepsilon$ and $Y_2^+$ is as
close to $Y_0^+$ as we wish.  By periodicity and Lemma~\ref{lem:lattice}(ii),
$|\W^{\le0}|\le Ch\,e^{-2c/h}e^{\Im I/h}$, which is at most
$e^{-c/h}e^{\Im I/h}$ for small $h$, and
$e^{\Im I/h}=|e^{-\frac\ii hI}|$.  Finally, $U_+-A_L=\Im J_+$ by
Lemma~\ref{lem:actions}, and integrating along $[0,\ii y_*]$ and then along the
part of $\Gamma_+$ from $\ii y_*$ to $x_+$, on which $p$ is real and $y$
increases, gives $\Im J_+\ge c$.  If the exit point is $1+\ii y_*$, the same
argument gives $|\W^{\le0}|\le Ch\,e^{-2c/h}|e^{\frac\ii hI}|$ and
$U_+-A_R=\Im(J_+-I)=\Im\int_1^{x_+}p\ge c$, with
$c=\int_0^{y_*}\Re p(1+\ii t)\,dt$.

\emph{Case (b$_-$).}  Now $\Gamma_-$ first meets the boundary at $-\ii y_*$ or
$1-\ii y_*$, where $y_*=y_*^-$; consider $-\ii y_*$.  As before,
$(0,1]-\ii y_*\subset D_R$, $\psi_1$ continues to the left along this line, and
on $(\delta,w)-\ii y_*$
\[
 |\W(z)|\le Ch\exp\Bigl(\frac1h\bigl(\Im I-2\Im S_0(z)\bigr)\Bigr),\qquad
 \Im S_0(z)\ge-\int_0^{y_*}\Re p(-\ii t)\,dt .
\]
By Lemma~\ref{lem:lattice}(iv),
$|\W^{\ge1}|\le Ch\exp\bigl(\frac1h(\Im I-2c)\bigr)$ with
$c=\int_0^{y_*}(\pi-\Re p(-\ii t))\,dt>0$.  Moreover $U_--A_L=\Im J_-$, and
integrating $p-\pi$ along $[0,-\ii y_*]$ and then along $\Gamma_-$ gives
$\Im J_-\ge c$.  The exit point $1-\ii y_*$ is treated symmetrically.
\end{proof}

\subsection{Exterior asymptotics}

\begin{proposition}\label{prop:exterior}
Let $V$ be real-analytic on $[0,1]$, $\alpha>0$, and let $E\notin\Sweep$
with $\sigma\notin\{f(0),f(1)\}$.  Then
\[
 p_N(E)=He^{N\Theta_\infty}\bigl(1+O(N^{-1})\bigr),\qquad H\ne0,
\]
where $\Theta_\infty=\int_0^1\log\rho_{\rm dom}(x)\,dx$ for a continuous branch
of the logarithm of the dominant frozen root, so that
$\Re\Theta_\infty=\Pot{\M}(E)$.
\end{proposition}

\begin{proof}
Since $E\notin\Sweep$, the number $(E-V(x))/(2\alpha)$ is never in
$[-1,1]$.  Hence there is no turning point on $[0,1]$, and there is a
continuous solution $p$ of $2\cos p=(E-V)/\alpha$ on $[0,1]$ with
$0<\Re p(0)<\pi$; moreover $\Im p\ne0$ on $[0,1]$, since $\cos p\in[-1,1]$ for real
$p$.  Fix a branch of $(\sin p)^{1/2}$ continuous on a neighbourhood of
$[0,1]$.  Suppose that $\Im p<0$.  The case $\Im p>0$ is the same, with the
plus and minus solutions exchanged; it gives $H=B(\sin p(1))^{-1/2}$ and
$\Theta_\infty=\log\alpha-\ii\int_0^1p$.  Because $\sigma\ne f(0)$, a short
vertical segment through $0$ satisfies the hypothesis of
Lemma~\ref{lem:segment} for $p$.  Let $\psi_0$ be the plus WKB solution
normalized at $0$ along it.  At $1$ we apply Lemma~\ref{lem:segment} to the
branch $q=\pm p+2\pi l$ with $0<\Re q(1)<\pi$; this branch exists because
$\sigma\ne f(1)$.  Take the WKB solution for $q$ at $1$ whose phase is
$-\frac\ii h\int_1^zp$ up to a multiple of $2\pi z/h$: the minus solution if
$q=p+2\pi l$, and the plus solution if $q=-p+2\pi l$.  Multiply it by
$e^{\pm2\pi\ii lz/h}$, which equals $1$ on $h\Z$, and by a constant (needed
when $\sin q=-\sin p$), so that
$\psi_1=(\sin p)^{-1/2}e^{-\frac\ii h\int_1^zp}(1+O(h))$ near $1$.  Since
$\Im p<0$, Theorem~\ref{thm:horizontal} continues $\psi_0$ to the right and
$\psi_1$ to the left along $[0,1]$.  Lemma~\ref{lem:wronskian} gives
$W[\psi_1,\psi_0]=-2\ii\Lambda(1+O(h))\ne0$ with $\Lambda=e^{\frac\ii h\int_0^1p}$.

As in the proof of Proposition~\ref{prop:reduction}, write $\varphi=A'\psi_0
+B'\psi_1$ on $h\Z$.  The values at $k=0$ and $k=-1$ are
\[
 \psi_0(kh)=\frac{e^{\ii kp(0)}(1+O(h))}{(\sin p(0))^{1/2}},\qquad
 \psi_1(kh)=\frac{\Lambda e^{-\ii kp(0)}(1+O(h))}{(\sin p(0))^{1/2}} .
\]
Hence
$A'=A(1+O(h))$ and $B'=B\Lambda^{-1}(1+O(h))$, and
\[
 \varphi_N=\frac{A\Lambda(1+O(h))+B\Lambda^{-1}(1+O(h))}{(\sin p(1))^{1/2}} .
\]
Since $|\Lambda|=e^{\frac1h\int_0^1|\Im p|}$, the second term is exponentially
smaller.  Hence $p_N=\alpha^N\varphi_N=He^{N\Theta_\infty}(1+O(h))$, where
\[
 H=A(\sin p(1))^{-1/2},\qquad
 \Theta_\infty=\log\alpha+\ii\int_0^1p=\int_0^1\log\rho_+ ,
\]
and $\rho_+$ is the dominant frozen root.  Finally, by Lemma~\ref{lem:PotM},
$\Re\Theta_\infty=\log\alpha+\int_0^1|\Im p|=\Pot{\M}(E)$.
\end{proof}

\subsection{Proofs of Theorem~\ref{thm:sharp} and of the Main Theorem}

\begin{proof}[Proof of Theorem~\ref{thm:sharp}]
Insert Proposition~\ref{prop:switch} into \eqref{eq:reduction}, noting that
$\frac A{2\ii}\cdot(-2)=\ii A$ and $\frac A{2\ii}\cdot2=-\ii A$.  In case
(b$_+$) the corresponding part of $\W$ is $O(e^{-c/h}M_N)$, and so is the term
$\ii Ae^{\ii N(2J_+-I)}$ of \eqref{eq:sharp}, since by
Lemma~\ref{lem:actions}
$|e^{\ii N(2J_+-I)}|=e^{N(\Im I+2A_L-2U_+)}\le e^{-2cN}|e^{-\ii NI}|$ if the
exit point is on $\Re z=0$, and similarly with $A_R$ and $e^{\ii NI}$
otherwise.  Both are therefore absorbed in the error term.  Case (b$_-$) is
the same, with $J_-$ and $U_-$.
\end{proof}

\begin{proof}[Proof of the Main Theorem]
We apply Proposition~\ref{prop:expsum} and Corollary~\ref{cor:charts} to
$q_N=p_N$, with the closed null set $F$ formed by $\partial\Sweep$ and the two
horizontal lines $\sigma=\sigma_0$ and $\sigma=\sigma_1$, which contain
the horizontal edges of $\Sweep$.  The components of $\C\setminus F$ are
$\Sweep^\circ$ and the components of
$\C\setminus(\Sweep\cup\{\sigma=\sigma_0\}\cup\{\sigma=\sigma_1\})$.  Off the two
lines, $\mathfrak c(0)$ and $\mathfrak c(1)$ are non-real, so the branches used
in Proposition~\ref{prop:exterior} depend holomorphically on $E$ in each of
the latter components, and $\Theta_\infty$ is single-valued and holomorphic
there; this is why the lines are included in $F$.  On
$\Sweep^\circ\setminus \mathcal E$, Theorem~\ref{thm:sharp} gives
$p_N=\sum_{j=1}^4e^{N\Theta_j}(c_j+o(1))$ with $c_1=A(\sin p(1))^{-1/2}$,
$c_2=B(\sin p(1))^{-1/2}$ and $c_{3,4}=\pm\ii c_1$, all non-zero; the error
$O(e^{-cN}M_N)$ is absorbed in the $o(1)$ of a term of maximal modulus.
The phases $\Theta_j$ are holomorphic on the whole of $\Sweep^\circ$.  Outside
$\Sweep$, Proposition~\ref{prop:exterior} gives
$p_N=He^{N\Theta_\infty}(1+O(N^{-1}))$ with $H\ne0$ and
$\Re\Theta_\infty=\Pot{\M}$.

\emph{Non-resonance.}  The six differences $\Re(\Theta_i-\Theta_j)$ are, by
\eqref{eq:realparts}, twice the differences of the actions.  They are harmonic
on the connected set $\Sweep^\circ$, and each of them changes sign:
\begin{itemize}
\item $A_L\to0$ at the edge $\sigma=\sigma_0$, $A_R\to0$ at the edge
$\sigma=\sigma_1$;
\item $U_+\to0$ at the right boundary $\tau=g+2\alpha$, and $U_-\to0$ at the
left boundary $\tau=g-2\alpha$;
\item in each case the other actions remain positive there.
\end{itemize}
Hence none of these differences vanishes identically.

The exponential-sum principle therefore gives $\mu_N\Rightarrow\mu$, and Corollary~\ref{cor:charts}
gives $\Pot{\mu}=\max_j\Re\Theta_j$ on $\Sweep^\circ$ and $\Pot{\mu}=\Pot{\M}$
on $\C\setminus\Sweep$.  By Lemma~\ref{lem:actions}, $\max_j\Re\Theta_j=
\Pot{\M}-2\min\{A_L,A_R,U_+,U_-\}$.  The right-hand side of
\eqref{eq:mainformula} is continuous on $\C$.  Since a subharmonic function
that agrees almost everywhere with a continuous one agrees with it everywhere,
\eqref{eq:mainformula} holds on $\C$.  As $\Pot{\mu}=\Pot{\M}$ is harmonic on
$\C\setminus\Sweep$, $\mu$ has no mass there, and on $\Sweep^\circ$ its support
lies in the tie set $T$ by Corollary~\ref{cor:charts}.  The densities are given
explicitly in Corollary~\ref{cor:curves}, and the number of arcs is finite by
Lemma~\ref{lem:finite}.

Finally, we show that $\mu(\partial\Sweep)=0$.  Put $m=\Pot{\M}-\Pot{\mu}\ge0$;
it vanishes outside $\Sweep$.  In a compact edge chart use Cartesian
coordinates $(u,v)$ in which the edge is given by $v=g(u)$, and put
$\delta=v-g(u)$.  We claim that $m=o(|\delta|)$ uniformly in the chart.  On the
horizontal edges this follows by factoring $A_L=x_c\widetilde A_L$ and
$A_R=(1-x_c)\widetilde A_R$: the normalized averages are continuous on the
closed crossing rectangle and vanish on their respective edges, while the
positive lower bound for $\Im V'$ bounds the crossing gap by a multiple of
$|\delta|$.  On the curved edges the turning gap is $|\delta|$, and
\[
 U_\pm\le M|\delta|\arccos\bigl(1-|\delta|/(2\alpha)\bigr)
 =o(|\delta|),
\]
with one constant $M$ for the whole sweep.

Take $\psi\in C_c^\infty$ supported in the chart and
$\chi\in C_c^\infty((-1,1))$, $\chi(0)=1$, and set
$\phi_\varepsilon=\psi(u)\chi(\delta/\varepsilon)$.
The Laplacian of this function is $O(\varepsilon^{-2})$ and is supported in a
graph strip of area $O(\varepsilon)$, on which $m=o(\varepsilon)$ uniformly.
Thus $\int m\,\Delta\phi_\varepsilon\,dA=o(1)$.  Since
$\Delta m=2\pi(\M-\mu)$ and $\M(\partial\Sweep)=0$, dominated convergence gives
\[
 \int_{\partial\Sweep}\psi\,d\mu
 =\lim_{\varepsilon\to0}\int\phi_\varepsilon\,d\mu
 =\lim_{\varepsilon\to0}\Bigl(\int\phi_\varepsilon\,d\M
 -\frac1{2\pi}\int m\,\Delta\phi_\varepsilon\,dA\Bigr)=0.
\]
So $\mu$ has no mass on the edges, and none at the corners, because
$\Pot{\mu}$ is continuous.  Hence $\supp\mu\subseteq\overline T$.
Conversely, take a point of $T$ and suppose that $\Pot{\mu}$ is harmonic on a
disk around it.  The difference between this potential and the real part of
either of the two maximizing phases is then harmonic and nonnegative on the
disk and vanishes at its centre, so by the minimum principle it vanishes
identically.  The real parts of the two phases would therefore coincide on the
disk, contrary to non-resonance.  Thus $T\subset\supp\mu$, so
$\supp\mu=\overline T$, and Lemma~\ref{lem:finite} now shows that the support
is the asserted finite union of arcs.
\end{proof}

\begin{remark}[Non-normality]\label{rem:deficit}
Let $\mathcal N_N$ be the strictly upper triangular part of a Schur form of
$H_N$.  Under \textup{(H1)} and \textup{(H2)},
\[
 \frac2\pi\int_\C\bigl(\Pot{\M}-\Pot{\mu}\bigr)\,dA
 =\lim_{N\to\infty}\frac{\|\mathcal N_N\|_F^2}N
 =\int_0^1|V|^2\,dx+2\alpha^2-\int_\C|E|^2\,d\mu .
\]
Thus the integrated deficit measures the departure of $H_N$ from normality.
Indeed,
$\|H_N\|_F^2=\sum_j|E_j|^2+\|\mathcal N_N\|_F^2$, and
$N^{-1}\|H_N\|_F^2\to\int_0^1|V|^2\,dx+2\alpha^2$.  Moreover
$N^{-1}\sum_j|E_j|^2=\int|E|^2\,d\mu_N\to\int|E|^2\,d\mu$,
because all supports lie in a fixed disk.  This gives the second equality.
For the first, $u=\Pot{\M}-\Pot{\mu}$ vanishes outside $\Sweep$ by the Main
Theorem, and $\Delta u=2\pi(\M-\mu)$.  Green's formula with $\Delta|E|^2=4$
gives $\int|E|^2\,d(\M-\mu)=\frac2\pi\int u\,dA$.  It remains to note that
$\int|w|^2\,d\nu_x(w)=|V(x)|^2+2\alpha^2$.
\end{remark}

\section{Purely imaginary and affine potentials}\label{sec:real}

\subsection{Reflection symmetry}
Let $V=-\ii B$, where $B$ is real on $[0,1]$ and $B'<0$. Then $g=0$ and
$\Sweep=[-2\alpha,2\alpha]+\ii[-\sigma_0,-\sigma_1]$. The anti-linear
symmetry $D\overline{H_N}D^{-1}=-H_N$, $D_{kk}=(-1)^k$, yields
$p_N(-\bar E)=(-1)^N\overline{p_N(E)}$. The reflection $E\mapsto-\bar E$
fixes $A_L,A_R$ and interchanges $U_+,U_-$; indeed,
$V(\bar z)=-\overline{V(z)}$ wherever both sides are defined, so the
inverse branches obtained by continuation from $x_c$ are reflected as well.

\begin{lemma}[Mirror inequality]\label{lem:mirror}
If $0<\tau=\Re E<2\alpha$ and $E\in\Sweep^\circ$, then $U_+<U_-$.
\end{lemma}
\begin{proof}
Put $b_\sigma(s)=|\Im\beta'(s-\ii\sigma)|$. The reflection symmetry of the
inverse makes this weight even in $s$, and (H2) makes it strictly positive.
Changing $t$ to $-s$ in $U_+$, and restricting the integral defining $U_-$ to
the same interval, we obtain
\[
 U_- - U_+\ \ge\!
 \int_0^{2\alpha-\tau}\!
 \left(\arccos\frac{s-\tau}{2\alpha}
       -\arccos\frac{s+\tau}{2\alpha}\right)b_\sigma(s)\,ds>0 .
\]
\end{proof}
\begin{corollary}\label{cor:real}
For a purely imaginary potential satisfying \textup{(H1)} and \textup{(H2)},
$\Pot{\mu}=\Pot{\M}-2\min(A_L,A_R,U)$ on $\Sweep$, where
\[
 U=\int_0^{2\alpha-|\tau|}\arccos\frac{|\tau|+s}{2\alpha}\,
            |\Im\beta'(s-\ii\sigma)|\,ds .
\]
The two turning-point costs can tie in the interior only on the imaginary
energy axis.
\end{corollary}
\begin{proof}
By reflection symmetry, the inverse-derivative weight is even. For $\tau\ge0$,
changing $t$ to $-s$ in $U_+$ gives the displayed integral; for $\tau\le0$,
the same integral is $U_-$. Lemma~\ref{lem:mirror} and reflection show that
this integral is $\min(U_+,U_-)$, and that in the interior the two costs are
equal exactly when $\tau=0$. The Main Theorem therefore gives the three-action
formula in the interior. On the edges $\sigma=\sigma_0,\sigma_1$, one endpoint
action vanishes; on the edges $|\tau|=2\alpha$, the displayed integral
vanishes. Since all the actions are nonnegative, the boundary formula follows
from the boundary values in the Main Theorem.
\end{proof}

\begin{remark}[A smaller extension domain]\label{rem:realH2}
For this symmetric potential the three-action \emph{potential} formula
needs only $\Sweep_r\subset V(\Omega)$ in (H2) for some $r>2\alpha$, with
$Y$ defined using $\Sweep_r$.  The margin $r-2\alpha$ ensures that the
vertical segments and the rectangle on which $\W$ is periodic reach
strictly below $\Im z_*$ (see below), and that the lens of
Lemma~\ref{lem:segment} fits inside $\Omega$. We now explain why the distant
turning point need not be reached.
Take $\tau>0$ and follow the upper arc to $x_+=\beta(E-2\alpha)$; both this arc
and its reflection in the real axis of the position plane are available, by
continuation of the inverse from $V(x_c)$. The reflected arc has real momentum.
Indeed, $f(\bar z)=f(z)$ and $r(\bar z)=-r(z)$, and $g=0$, so the reflected arc
is the initial part of $\Gamma_-$, on which $r$ runs over
$[0,2\alpha-\tau]\subset[0,\tau+2\alpha]$. If the
upper arc leaves the position strip before reaching $x_+$, use the first exit
point and its reflection in the two boundary estimates of
Proposition~\ref{prop:switch}; both groups of modes are then exponentially
smaller than an unswitched term. Otherwise the upper calculation is as
in case (a$_+$). On the lower side it suffices to use the horizontal line
through $z_* = \overline{x_+}$. Along the reflected arc from $x_c$ to
$z_*$ let $U^*=\int(\pi-p)|dy|$. The integral in the proof of
Lemma~\ref{lem:mirror} gives $U^*>U_+$, using only $0\le s\le2\alpha-\tau$.

In more detail, on the line through $z_*$ the endpoint solutions can be
continued to within $\delta$ of $z_*$ from their respective growing sides,
where $\Im p$ has a strict sign. Across the remaining distance, which is at most
$2\delta$, the difference equation bounds each two-component solution
vector by $K^{2\delta/h+2}$ times its incoming norm; here $K$ is independent
of small $\delta$, since the coefficients are bounded there. Thus on one
period near $z_*$ the bound for the Wronskian is
$O(\exp((\Im I-2\Im S_0(z_*)+O(\delta))/h))$. The constant in
$O(\delta)$ also accounts for moving the incoming actions to $z_*$.
Lemma~\ref{lem:lattice}(iv) contributes the factor
$e^{-2\pi|\Im z_*|/h}$; the resulting exponent is
$A_L+A_R-2U^*+O(\delta)$. Choosing $\delta$ after fixing $E$ makes
this strictly smaller than $A_L+A_R-2U_+$. Reflection handles $\tau<0$;
the omitted energy axis has zero area, and the formula extends to it by
continuity. This argument does not give an expansion involving a
turning point outside the domain.
\end{remark}

\subsection{Affine potentials and Lommel polynomials}
For $V(x)=V_0+\gamma x$ with $\Im\gamma>0$, $\beta'=1/\gamma$ and
$\Sweep$ is a parallelogram. Put $r_0=|\Im(1/\gamma)|$ and
$u=(\tau-g(\sigma))/(2\alpha)$. Direct integration gives
\[
 U_+=2\alpha r_0(\sqrt{1-u^2}-u\arccos u),\qquad
 U_-=2\alpha r_0(\sqrt{1-u^2}+u\arccos(-u)).
\]
In particular $U_--U_+=2\pi\alpha r_0u$. For nonzero $\gamma$, the
substitution $P_k=\gamma^k Q_k$ gives the unit-slope recurrence with energy
$(E-V_0)/\gamma$ and squared off-diagonal $\alpha^2/\gamma^2$.

Now take $V(x)=-\ii(1-2x)$. To compare with the order variable in
\cite{BlSt}, put $w=\ii E$, $\nu=N(w-1)/2$ and $z=\ii\alpha N$. By
\eqref{eq:recurrence}, $\varphi_k=P_k(E)/\alpha^k$ satisfies
\[
 \varphi_{k+1}=\frac{2(k+\nu)}z\,\varphi_k-\varphi_{k-1},
 \qquad \varphi_{-1}=0,\quad \varphi_0=1 .
\]
Hence $p_N(E)=\alpha^N R_{N,\nu}(\ii\alpha N)$, where $R_{N,\nu}$ is the
Lommel polynomial. In \cite[Eq.~(4)]{BlSt} the centred polynomial is
$Q_N^{(\alpha)}(\zeta)=R_{N,(1-N-N\zeta)/2}(\ii\alpha N)$; thus our
argument corresponds exactly to $\zeta=N^{-1}-w$. The zeros are symmetric under
$\zeta\mapsto-\zeta$ for every $N$: with the reversal $R$ of the sites and
$D=\operatorname{diag}((-1)^k)$ one has $DRH_NRD^{-1}=-H_N-(2\ii/N)I$, so
the spectrum is invariant under $w\mapsto2/N-w$, which is
$\zeta\mapsto-\zeta$.
For $0<\alpha<\alpha_0$, where $\alpha_0\simeq0.33137$ is the threshold
in \cite[Eq.~(8) and Theorem~1]{BlSt}, the central segment of the
two-turning-point tie lies on the imaginary $E$-axis. In the notation
of that paper, put
\[
 \chi_\alpha(\zeta)=\tfrac12\sqrt{4\alpha^2+(1+\zeta)^2}
 +\tfrac{1+\zeta}{2}\log\frac{2\alpha}{1+\zeta+\sqrt{4\alpha^2+(1+\zeta)^2}} .
\]
\begin{lemma}\label{lem:chi}
For $V(x)=-\ii(1-2x)$ and $E\in\Sweep^\circ$,
$U_+-A_L=\Im J_+=\Re\chi_\alpha(-w)-\frac\pi4\Im w$, $w=\ii E$.
\end{lemma}
\begin{proof}
Substitute $v=(w-1+2z)/(2\ii\alpha)$ in $J_+$. The new variable runs from
$v_0=\ii(1-w)/(2\alpha)$ to $1$ along a segment whose interior lies in the
upper half-plane. There the principal arccosine gives
\[
 J_+=\ii\alpha\int_{v_0}^1\arccos v\,dv
     =\ii\alpha\bigl(\sqrt{1-v_0^2}-v_0\arccos v_0\bigr)
     =\ii\chi_\alpha(-w)+\tfrac\pi4(1-w).
\]
For the last equality write $v_0=\ii t$, $\Re t>0$, and use
$\arccos(\ii t)=\pi/2-\ii\log(t+\sqrt{1+t^2})$.
The first equality in the lemma is Lemma~\ref{lem:actions}.
\end{proof}
Thus the curve $A_L=U_+$ is given by precisely the equation of \cite{BlSt};
reflection and the symmetry $V(1-x)=-V(x)$ yield the other branches.
The Charlier limit in \cite{BlStC} also admits a description by curves and
densities, but its recurrence is different: the monic Charlier polynomials
with parameter $a$ satisfy
$C_{k+1}(t)=(t-k-a)C_k(t)-akC_{k-1}(t)$.
In the scaling of \cite{BlStC} the parameter is negative and proportional
to $N$. The resulting product of adjacent off-diagonal coefficients varies
with $k$, whereas the Dirichlet recurrence \eqref{eq:recurrence} has the fixed
product $\alpha^2$.

\subsection{Changes of normalization}
More generally, a symmetric tridiagonal matrix $J_N$ with potential $B$
and common nonzero off-diagonal coefficient $c$ satisfies
$J_N=\omega H_N(B/\omega,|c|)$, where $\omega=c/|c|$ and
$H_N(V,\alpha)$ denotes \eqref{eq:TN}. In particular, the convention
with $c=\ii\alpha$ has $V=-\ii B$ and $E=-\ii w$, and
\[
 \det(EI-H_N(V,\alpha))=(-\ii)^N\det(wI-J_N),\qquad w=\ii E .
\]
This is a rotation of the energy plane, not a change in the spectral
problem at fixed $B$ and $w$: counting measures are pushed forward by
$w\mapsto-\ii w$, and logarithmic potentials rotate without an additive
constant. Distinct constant off-diagonals $c_-,c_+$ with
$c_-c_+=\alpha^2$ give the same Dirichlet determinants, by the continuant
recurrence; for such matrices we assert only that the eigenvalues coincide.
Throughout the paper, $H_N$ denotes the discrete Schr\"odinger operator with
equal real off-diagonals in \eqref{eq:TN}.

\section{Numerical comparison in the monotone case}\label{sec:numerics}
Figures~\ref{fig:eigs} and~\ref{fig:complex} compare the Dirichlet
eigenvalues for $N=400$ with the curves on which two of the four tunnelling
actions tie for the minimum. The eigenvalues were computed by Aberth
iteration, using the determinant recurrence and its derivative, in
$340$-digit arithmetic, and the plotted points are rotated
into the physical $E$-plane as in Section~\ref{sec:real}. The tie curves
were located by evaluating the actions on a grid and refining by bisection. For
$V(x)=-\ii(1-2x+0.6x^2)$ and $\alpha=0.15$, $142$ computed eigenvalues lie on
the imaginary axis, between $-0.5196\ii$ and $0.0100\ii$; the predicted
segment is $[-0.5190\ii,0.0119\ii]$.

For the genuinely complex quadratic potential, Table~\ref{tab:complex} lists
finite-size values computed directly from the Dirichlet determinant via \eqref{eq:recurrence}. The column
headed $L$ is the right-hand side of \eqref{eq:mainformula}, computed
independently by quadrature of \eqref{eq:actions}. Let $S_N$ be the sum of
the four leading terms on the right-hand side of \eqref{eq:sharp}, and $D_N$ the
largest of their moduli. The last column gives
\[
 \varepsilon_N=
 \big|(\sin p(1))^{1/2}p_N(E)/\alpha^N-S_N\big|/D_N
 \quad\text{at }N=800 .
\]
The recurrence was
evaluated with $\lfloor 0.7N\rfloor+80$ decimal digits, and the action and
phase integrals with 53-digit working precision. Each of
the four actions is minimal in at least one row. In the second
row $U_-$ and $A_L$ are close, and the finite-size values converge
more slowly than a single dominant exponential would suggest: at
$N=800$ the second largest of the four terms has relative size about
$0.02$ in that row, but less than $10^{-46}$ in all other rows.  Thus
$\varepsilon_{800}$ mainly tests the dominant term; the switched terms are
tested by the rows in which they dominate.
\begin{table}[tbp]
\centering\small
\caption{Complex quadratic potential $V(x)=0.6x^2-\ii(1-2x)$, $\alpha=0.4$.
Here $L_N=N^{-1}\log|p_N(E)|$.}\label{tab:complex}

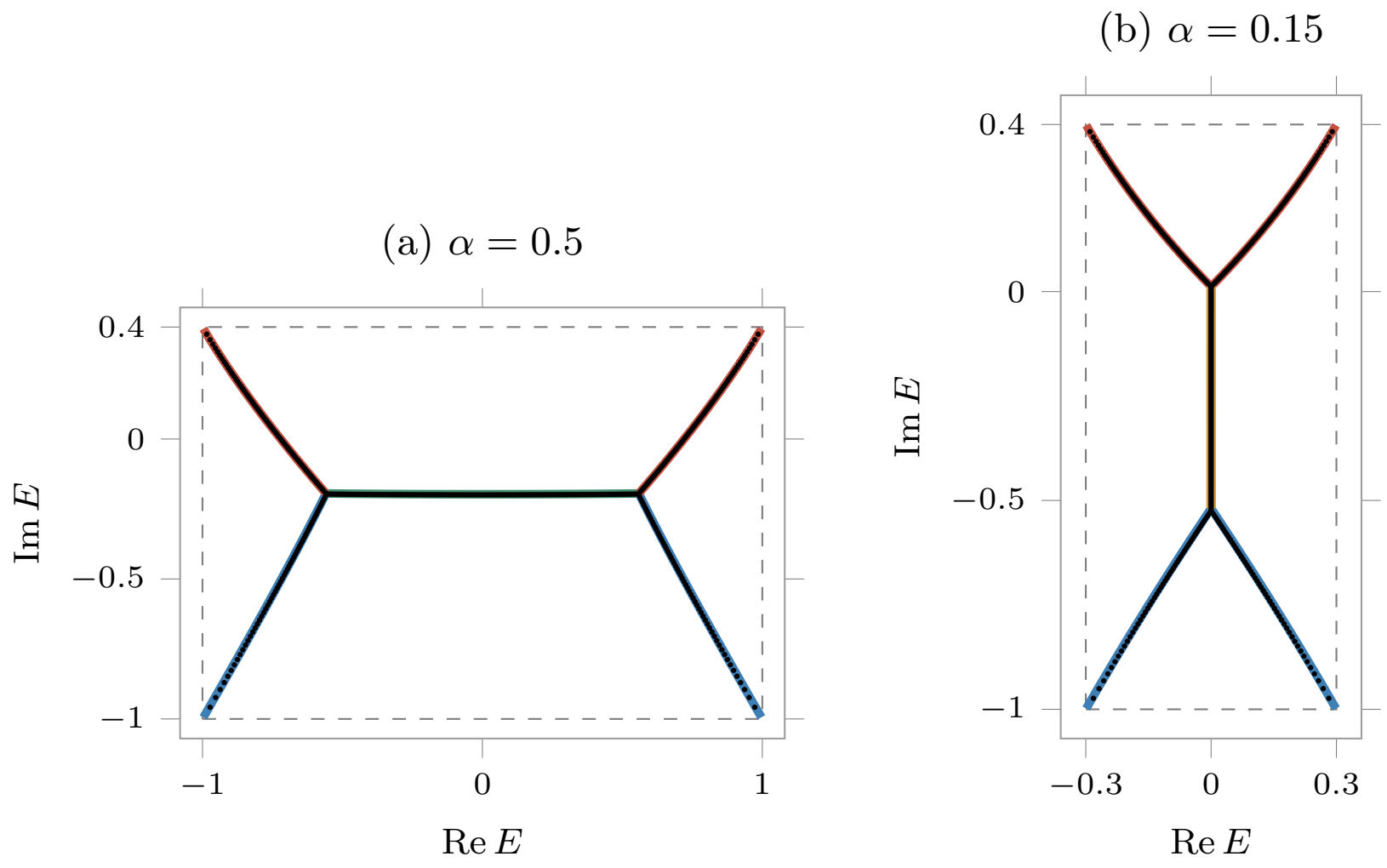
\caption{Purely imaginary potential $V(x)=-\ii(1-2x+0.6x^2)$, $N=400$.
Dots are eigenvalues; the dashed rectangle is the frozen sweep. Coloured
lines mark ties for the minimum in \eqref{eq:mainformula}: $A_L=A_R$
(green), $A_L=U_\pm$ (blue), $A_R=U_\pm$ (red), $U_+=U_-$ (gold).
The gold branch is visible in (b), on the imaginary axis.}\label{fig:eigs}
\end{figure}

\begin{figure}[tbp]
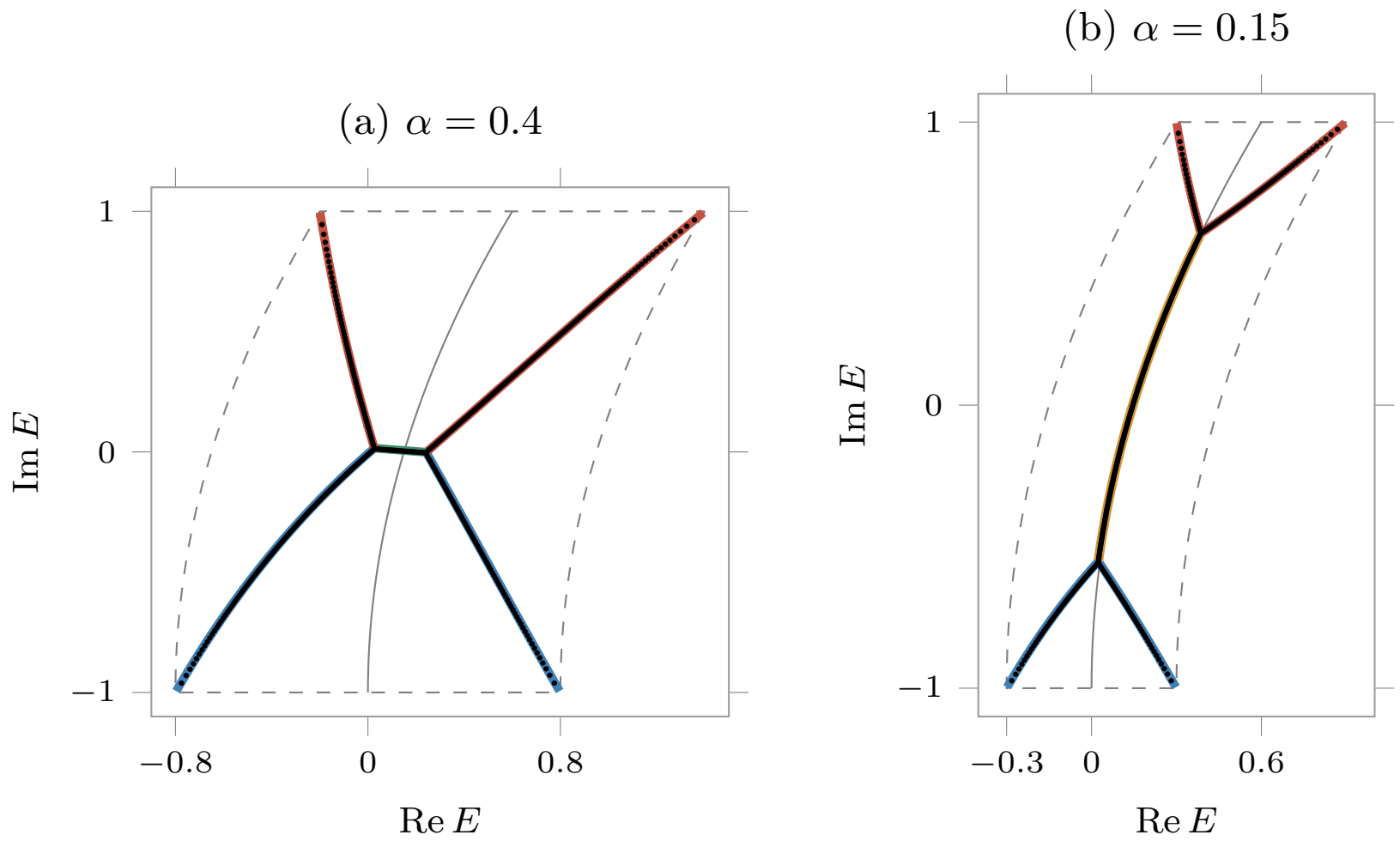

\centering
%
\caption{Complex quadratic potential $V(x)=0.6x^2-\ii(1-2x)$, $N=400$.
The grey curve is $V([0,1])$; its horizontal translates (dashed) bound the sweep.
Colours are as in Figure~\ref{fig:eigs}. In (b) the gold
two-turning-point branch bends away from the imaginary axis.}\label{fig:complex}
\end{figure}

\section{Beyond monotonicity}\label{sec:general}

In this section $V\colon[0,1]\to\C$ is an arbitrary real-analytic function;
$f$ need not be monotone, and the curve $V([0,1])$ may have horizontal
tangents.  We show that the structure found above persists: the problem
reduces exactly to switching coefficients at the crossing points, and the
natural candidate for the limit is governed by the minimum of a finite family of
tunnelling costs.  We prove existence of the limit, and the formula, under an asymptotic
hypothesis on the switching coefficients that holds in the monotone case.
For univalent $V$ we prove this hypothesis under a geometric condition on the
curves of real momentum; for a circular arc the condition holds almost everywhere.
We also give numerical comparisons outside the monotone case.

\subsection{Crossing points and layers}
Let $E=\tau-\ii\sigma$.  A point $x\in(0,1)$ with $f(x)=\sigma$ and
$|\tau-r(x)|<2\alpha$ is a \emph{crossing point} of $E$; the number
$m(E)$ of crossing points equals the number of frozen segments that contain $E$.  Let $\mathcal E_0$
be the set of $E\in\Sweep$ for which some solution of $f(x)=\sigma$
in $[0,1]$ is an endpoint or satisfies $f'(x)=0$, or for which
$E\mp2\alpha\in V([0,1])$.  This set is contained in finitely many analytic
curves, so it has zero area.  The function $m$ is locally constant on
$\Sweep^\circ\setminus \mathcal E_0$; we call $m(E)$ the number of \emph{layers}
at $E$, and the arcs of the curves $V([0,1])\pm2\alpha$ across which the number of layers jumps
the \emph{layer boundaries}.

Fix $E\in\Sweep^\circ\setminus \mathcal E_0$, with crossing points
$0<x_1<\dots<x_m<1$, and put $x_0=0$, $x_{m+1}=1$, $I_k=(x_k,x_{k+1})$.  Since there
is no turning point on $[0,1]$, the equation
$2\cos p=(E-V)/\alpha$ has a continuous solution $p$ on $[0,1]$ with
$0<\Re p(0)<\pi$; in general, $\Re p$ leaves $(0,\pi)$ at the real
solutions of $f=\sigma$ that are not crossing points.  Since $\cos p$ is real only
where $f=\sigma$, and lies in $(-1,1)$ exactly at the crossing points,
$\Im p$ vanishes exactly at $x_1,\dots,x_m$, and it changes sign there because
$f'(x_k)\ne0$.  Thus one frozen root dominates on each $I_k$, and the
dominant root changes at every crossing point.  Put
\[
 A_k=\int_{I_k}|\Im p(x)|\,dx\qquad(0\le k\le m),\qquad
 \Pot{\M}=\log\alpha+\sum_{k=0}^mA_k
\]
(Lemma~\ref{lem:PotM}).  For monotone $f$ we have $m=1$, $A_0=A_L$ and
$A_1=A_R$.

\subsection{Real-momentum arcs and switching costs}
Let $S=[E-2\alpha,E+2\alpha]$ and $w_\pm=E\mp2\alpha$.
Every crossing point satisfies $V(x_k)\in S$.  We assume the following.
\begin{enumerate}
\item[(H2$'$)] $V$ extends holomorphically to a domain $\Omega\supset[0,1]$
with $V'\ne0$ on $\Omega$, and for $E\in\Sweep^\circ\setminus \mathcal E_0$ the
local inverse of $V$ at $V(x_k)$ continues analytically in $\Omega$ along
both segments of $S$ from $V(x_k)$ to $w_\pm$.
\end{enumerate}
This holds, for instance, if $V$ is univalent on $\Omega$ and $V(\Omega)$
contains every segment $V(x)+[-4\alpha,4\alpha]$.  Denote by $\Gamma_{k,\pm}$ the arc traced by the continued inverse; it runs from
$x_k$ to a turning point $x_{k,\pm}$ with $V(x_{k,\pm})=w_\pm$.  The
continued momentum is real on $\Gamma_{k,\pm}$ and moves monotonically from
$p(x_k)$ to a value $p_{k,\pm}\in\pi\Z$.  Let $s_k=1$ if $\rho_+$ dominates on
$I_{k-1}$, and $s_k=-1$ otherwise; in the notation of
Section~\ref{sec:exact}, $s_k=\varsigma_{k-1}$.

\begin{definition}\label{def:admissible}
The \emph{switch} at $x_k$ through $x_{k,t}$, $t\in\{+,-\}$, is
\emph{admissible} if $\Gamma_{k,t}$ meets the real axis only at $x_k$.  Its
\emph{cost} is
\[
 U_{k,t}=\Im\int_{\Gamma_{k,t}}s_k\,(p-p_{k,t})\,dz .
\]
\end{definition}

If the switch at $x_k$ through $x_{k,t}$ is admissible, we call $t$ an
\emph{admissible end} at $x_k$. We also put
\begin{equation}\label{eq:calU}
 \mathcal U_{k,t}=-\ii s_k\Bigl(\int_{\Gamma_{k,t}}p\,dz-p_{k,t}\,x_{k,t}\Bigr).
\end{equation}
Since $\Im(p_{k,t}x_{k,t})=\Im\int_{\Gamma_{k,t}}p_{k,t}\,dz$, we have
$\Re\mathcal U_{k,t}=U_{k,t}$.

The quantity $-2U_{k,t}$ is the real part of the extra phase acquired by
following the dominant root from $x_k$ to $x_{k,t}$ and returning along the other
root.

\begin{lemma}\label{lem:costpositive}
For an admissible switch,
\[
 U_{k,t}=\int_{\Gamma_{k,t}}|\Im z|\,|dp|>0 .
\]
If $V$ is univalent on a domain containing all the arcs, then
$\Gamma_{k,+}\cup\Gamma_{k,-}=\beta(S)$ is the same arc $\Gamma$ for every $k$,
and this arc passes through all crossing points. The switch through $x_\pm$ is
admissible only for the crossing point that $\Gamma$ meets first when
traversed from $x_\pm$, and not even for that one if $\Gamma$ meets
$\Omega\cap\R$ outside $[0,1]$ before reaching it.  In particular, at most two crossing points have
admissible switches.
\end{lemma}

\begin{proof}
On $\Gamma_{k,t}$ the momentum is real, so $\Im((p-p_t)dz)=(p-p_t)\,dy$.
Integrating by parts, and using $\Im x_k=0$ and $p(x_{k,t})=p_t$, we get
$U_{k,t}=-s_k\int_{\Gamma_{k,t}}y\,dp$.  For an admissible switch, $y$ has a
constant sign $\varepsilon$ on $\Gamma_{k,t}\setminus\{x_k\}$, and $dp$ has a
constant sign $\delta$.  It remains to show that $s_k\varepsilon\delta=-1$; this
is a local computation at $x_k$.  From
$2\cos p=(E-V)/\alpha$ we get $p'=V'/(2\alpha\sin p)$ on the
real axis.  Hence the sign of $\Im p$ on $I_{k-1}$ near $x_k$ is
$\eta_k\,\mathrm{sgn}(f_x(x_k))$, where $\eta_k=\mathrm{sgn}\sin p(x_k)$, and therefore
$s_k=-\eta_k\,\mathrm{sgn}(f_x(x_k))$.  Along $\Gamma_{k,t}$ we have
$z=\beta(t'-\ii\sigma)$, so $dy=\Im(1/V')\,dt'$, and
$p=\eta_k\arccos\frac{\tau-t'}{2\alpha}+2\pi j$ for a fixed $j\in\Z$.  Moving
towards $w_+$ means $dt'<0$; then $\varepsilon=-\mathrm{sgn}(f_x(x_k))$ and
$\delta=-\eta_k$.  Moving towards $w_-$ gives $\varepsilon=\mathrm{sgn}(f_x(x_k))$ and $\delta=\eta_k$.  In both cases $s_k\varepsilon\delta=-1$.  The statements about
the univalent case follow because $\beta$ is single-valued.
\end{proof}

Under \textup{(H1)}, the lemma gives a second formula for the switching actions of
the Main Theorem: $U_\pm=\int_{\Gamma_\pm}|\Im z|\,|dp|$.

\subsection{Histories and the predicted potential}

\begin{definition}\label{def:history}
A \emph{history} $H$ is a word $w_0w_1\cdots w_m$ in the letters $D,R$ with
no two consecutive letters $R$, together with an admissible end $t_k$ for
every $k$ with $w_{k-1}=w_k=D$.  Its \emph{cost} is
\[
 c(H)=\sum_{k:\,w_k=R}A_k+\sum_{k:\,w_{k-1}=w_k=D}U_{k,t_k}.
\]
Let $C(E)$ be the minimal cost of a history; it is defined for
$E\in\Sweep^\circ\setminus \mathcal E_0$.  Put $F=\Pot{\M}-2C$ there and
$F=\Pot{\M}$ outside $\Sweep$.
\end{definition}

The letter $w_k$ records whether the solution follows the dominant ($D$) or
the recessive ($R$) root on $I_k$.  Since the dominant root changes at every
crossing point, following one root continuously through a
crossing point turns $D$ into $R$ and $R$ into $D$ at no cost.  Keeping $D$
on both sides requires a switch through a turning point, and $R$ cannot
follow $R$.  The two alternating words are always histories, so
$C(E)<\infty$; they correspond to following $\rho_+$ and
$\rho_-$, respectively, along all of $[0,1]$.  For monotone $f$ the histories are
$DD$ (twice, once with each end), $DR$ and $RD$, and
$F=\Pot{\M}-2\min\{A_L,A_R,U_+,U_-\}$ is the formula of the Main Theorem.

\subsection{An exact reduction}\label{sec:exact}
Let $\varsigma_k$ be the sign for which $e^{\frac\ii h\varsigma_k\int p}$ grows to
the right on $I_k$ ($0\le k\le m$), and put $\varsigma_{m+1}=-\varsigma_m$.
These signs alternate. At a normalization site one may choose a local
branch $q=\pm p+2\pi l$ with $0<\Re q<\pi$. Applying
Lemma~\ref{lem:segment} to this branch, and multiplying by a constant and
an $h$-periodic exponential, gives the phase and normalization stated below, expressed in
terms of $p$. In particular, the construction also covers a branch continued past
one of the principal arccosine cuts.  A \emph{WKB solution of type $\varsigma$} is a solution
with the representation $(\sin p)^{-1/2}e^{\frac\ii h\varsigma\int p}(1+O(h))$;
it is a plus solution if $\varsigma=1$ and a minus solution if $\varsigma=-1$.
Define solutions of \eqref{eq:analytic} as follows:
\begin{itemize}
\item $u_0$ is the WKB solution of type $\varsigma_0$ normalized at $0$;
\item $u_k$, $1\le k\le m$, is the WKB solution of type $\varsigma_k$
normalized at $x_k$ along a short vertical segment;
\item $u_{m+1}$ is the WKB solution of type $\varsigma_{m+1}$ normalized at $1$.
\end{itemize}
Each $u_k$ with $1\le k\le m$ grows away from $x_k$ in both directions.  As
in Lemma~\ref{lem:realaxis}, for every small $\delta>0$ its WKB representation holds
on neighbourhoods of $[x_{k-1}+\delta,x_k]$ and of $[x_k,x_{k+1}-\delta]$, with
$\delta=0$ at $x_0=0$ and $x_{m+1}=1$.

\begin{theorem}\label{thm:recursion}
Let $E\in\Sweep^\circ\setminus \mathcal E_0$.  On $h\Z$,
\[
 u_{k-1}=\gamma_ku_{k+1}+\alpha_ku_k\qquad(1\le k\le m),\qquad
 \gamma_k=e^{\frac\ii h\varsigma_{k-1}\int_{x_{k-1}}^{x_{k+1}}p}\,(1+O(h)),
\]
with \emph{switching coefficients}
$\alpha_k=W[u_{k+1},u_{k-1}]/W[u_{k+1},u_k]$.  Write
$\varphi=\mathcal P_ku_k+\mathcal Q_ku_{k+1}$ on $h\Z$, and let $(a_0,a_1)=(A,B)$ if
$\varsigma_0=1$ and $(a_0,a_1)=(B,A)$ if $\varsigma_0=-1$.  Then
$\mathcal P_0=a_0(1+O(h))$, $\mathcal Q_0=a_1\,e^{\frac\ii h\varsigma_1\int_0^{x_1}p}(1+O(h))$,
and
\[
 \mathcal P_k=\alpha_k\mathcal P_{k-1}+\mathcal Q_{k-1},\qquad
 \mathcal Q_k=\gamma_k\mathcal P_{k-1},
\]
\[
 p_N(E)=\alpha^N\bigl(\mathcal P_mu_m(1)+\mathcal Q_mu_{m+1}(1)\bigr).
\]
Consequently, $\alpha^{-N}(\sin p(1))^{1/2}p_N(E)$ is the sum of the terms below over all words in $D,R$ with no two consecutive letters $R$; admissible ends are not yet chosen at this stage. The terms are
\[
 T_w=a_w\,e^{\frac\ii h\Theta_w}(1+O(h))\prod_{k:\,w_{k-1}=w_k=D}\tilde\alpha_k,
\]
where $\Theta_w=\sum_{j=0}^{m}\epsilon_j\int_{x_j}^{x_{j+1}}p$,
with $\epsilon_j=\varsigma_j$ for $w_j=D$ and
$\epsilon_j=-\varsigma_j$ for $w_j=R$. Here $\tilde\alpha_k=\alpha_k\,e^{-\frac\ii h\varsigma_{k-1}\int_{x_{k-1}}^{x_k}p}$,
$a_w\in\{A,B\}$, and $\Theta_w$ is real-analytic, and holomorphic
if $w$ contains no pair $DD$ (see Lemma~\ref{lem:phases}); moreover,
\[
 \bigl|e^{\frac\ii h\Theta_w}\bigr|=
 \exp\Bigl(\frac1h\Bigl(\Pot{\M}-\log\alpha-2\sum_{k:\,w_k=R}A_k\Bigr)\Bigr).
\]
\end{theorem}

\begin{proof}
As in Lemma~\ref{lem:realaxis}, $W[u_k,u_{k-1}]$ is evaluated on $I_{k-1}$ and
$W[u_k,u_{k+1}]$ on $I_k$.  Both are non-degenerate, since $u_k$ and
$u_{k\pm1}$ are of opposite types, and their quotient is $\gamma_k$.  The
coefficient $\alpha_k$ is then determined by the identity.  The initial values follow
as in Proposition~\ref{prop:reduction}, using the representations of $u_0$
and $u_1$ at $0$ and $-h$.  Inserting $u_{k-1}=\gamma_ku_{k+1}+\alpha_ku_k$
into $\varphi=\mathcal P_{k-1}u_{k-1}+\mathcal Q_{k-1}u_k$ gives the recursion.  Expanding the
recursion gives one term for each word: a letter $R$ on $I_k$ comes from
a factor $\gamma_k$ if $k\ge1$ and from $\mathcal Q_0$ if $k=0$; a final letter $D$ on $I_m$
comes from $u_m(1)$; and a pair $DD$ comes from a factor $\alpha_k$.
The moduli follow from those of $\gamma_k$, $\mathcal Q_0$ and $u_m(1)$, which are given by the WKB
forms, and from $u_{m+1}(1)=(\sin p(1))^{-1/2}(1+O(h))$.
\end{proof}

Thus the whole problem reduces to the normalized switching coefficients
$\tilde\alpha_k$.  Under \textup{(H1)}, we have $m=1$ and $\varsigma_0=1$.  In
cases (a$_+$) and (a$_-$), comparing with \eqref{eq:reduction} and
Proposition~\ref{prop:switch} gives
\[
 \tilde\alpha_1=\ii\,e^{-\frac2h\mathcal U_+}(1+O(h))
 -\ii\,e^{-\frac2h\mathcal U_-}(1+O(h)),\qquad
 \mathcal U_\pm=-\ii\Bigl(J_\pm-\int_0^{x_c}p\Bigr),
\]
and $\Re\mathcal U_\pm=U_\pm$ by Lemma~\ref{lem:actions}; these are the
functions \eqref{eq:calU} with $m=1$.  In case (b$_\pm$) the corresponding term
is replaced by an exponentially smaller remainder.

\begin{lemma}\label{lem:phases}
Let $H$ be a history with word $w$ and ends $t_k$, and put
\[
 \Phi_H=\ii\,\Theta_w-2\sum_{k:\,w_{k-1}=w_k=D}\mathcal U_{k,t_k},
\]
with $\mathcal U_{k,t}$ as in \eqref{eq:calU}.  On every open subset of
$\Sweep^\circ\setminus \mathcal E_0$ on which $m$ and the admissible ends are fixed,
$\Phi_H$ is holomorphic and $\Re\Phi_H=\Pot{\M}-\log\alpha-2c(H)$.  In
particular, $\Pot{\M}-2c(H)$ is harmonic there.
\end{lemma}

\begin{proof}
The integrand $p(z,E)$ and the ends $x_{k,t}$, which solve
$V(x_{k,t})=E\mp2\alpha$, depend holomorphically on $E$.  The
crossing points $x_k$ depend only on $\sigma=-\Im E$, and they enter only as
endpoints of integrals.  At a pair $DR$ or $RD$ the solution follows one root
through $x_k$, so the integrand of $\Theta_w$ does not change there and $x_k$ is
not an endpoint.  At a pair $DD$ the type changes from $\varsigma_{k-1}$ to
$\varsigma_k=-\varsigma_{k-1}$, so $\partial_{x_k}(\ii\Theta_w)
=2\ii\varsigma_{k-1}p(x_k)$, while $\partial_{x_k}(-2\mathcal U_{k,t})
=-2\ii s_kp(x_k)$.  These cancel because $s_k=\varsigma_{k-1}$.  Hence $\Phi_H$
is holomorphic.  The formula for its real part follows from Theorem~\ref{thm:recursion} and
$\Re\mathcal U_{k,t}=U_{k,t}$.
\end{proof}

\subsection{Existence under a switching hypothesis}
We say that $E\in\Sweep^\circ\setminus \mathcal E_0$ satisfies \textup{(Sw)} if, for some choice of the local WKB solutions in
Theorem~\ref{thm:recursion}, there is $\eta>0$ such that, for every $k$,
\[
 \tilde\alpha_k=\sum_{t\in T_k}\kappa_{k,t}\,e^{-\frac2h\mathcal U_{k,t}}
 \bigl(1+o(1)\bigr)+r_k ,
\]
where:
\begin{enumerate}[(i)]
\item $T_k$ is a set of admissible ends at $x_k$, $\kappa_{k,t}\ne0$, and
$\mathcal U_{k,t}$ is given by \eqref{eq:calU};
\item $|r_k|\le e^{-\frac2h(C(E)+\eta)}
+\delta_{k1}e^{-\frac2h(A_0+\eta)}+\delta_{km}e^{-\frac2h(A_m+\eta)}$;
\item every admissible end $t\notin T_k$ satisfies $U_{k,t}>C(E)$, or
$k=1$ and $U_{1,t}\ge A_0$, or $k=m$ and $U_{m,t}\ge A_m$.
\end{enumerate}
Condition (iii) guarantees that the switches omitted from $T_k$ are never
needed in a history of minimal cost: if a history uses such a switch at
$x_1$, replacing $w_0=D$ by $R$ does not increase its cost, and similarly at
$x_m$.

Let $Z$ be the set of $E\in\Sweep^\circ\setminus \mathcal E_0$ at which the set
of admissible ends is not locally constant.  It is relatively closed in
$\Sweep^\circ\setminus \mathcal E_0$, and it has zero area.  To see the latter, note that the arcs
$\Gamma_{k,t}$ depend continuously on $E$ and cross $\R$ transversally
at $x_k$, and that transversal crossings persist.  Hence admissibility can change only
where an end $x_{k,t}$ lies on $\R$, that is, $E\mp2\alpha\in
V(\Omega\cap\R)$, or where $\Gamma_{k,t}$ is tangent to $\R$ at some
$z\in\Omega\cap\R$.  In the latter case $f_x(z)=0$ and $f(z)=\sigma$, so $\sigma$ is a
critical value of $f$ on $\Omega\cap\R$.  Both sets have zero area. The first is contained in a countable union of analytic arcs. As for the second, on
each real interval in $\Omega$ the real-analytic function $f$ has isolated
critical points unless it is constant, and a constant interval contributes
only one critical value; thus the second set is contained in a countable
union of horizontal lines.

\begin{theorem}\label{thm:conditional}
Assume \textup{(H2$'$)}, that \textup{(Sw)} holds for almost every
$E\in\Sweep$, and the non-resonance condition that, on every open chart with fixed admissible ends, the costs of distinct histories do not agree identically on any open subset.  Then $\mu_N$ converges weakly to a probability
measure $\mu$ with $\Pot{\mu}=F$ on $\C\setminus(\mathcal E_0\cup Z\cup\partial\Sweep)$,
hence almost everywhere.  Its support is contained in
$\mathcal E_0\cup Z\cup\partial\Sweep$ together with a countable family of
real-analytic arcs in $\Sweep^\circ\setminus(\mathcal E_0\cup Z)$, locally finite there,
on which two histories have equal minimal cost.
\end{theorem}

\begin{proof}
We insert (Sw) into the sum of Theorem~\ref{thm:recursion} and expand.  A term that
contains only main parts of the $\tilde\alpha_k$ corresponds to a history all
of whose switches belong to the sets $T_k$; it equals a non-zero constant
times $e^{\frac1h\Phi_H}(1+o(1))$, with $\Phi_H$ as in Lemma~\ref{lem:phases},
so its modulus is a non-zero constant times
$\exp(\frac1h(\Pot{\M}-\log\alpha-2c(H)))(1+o(1))$.  By (iii), the minimal
cost of such histories is $C(E)$.

Now consider the terms that contain a remainder.  We split every factor
$\tilde\alpha_k$ of a word into its main part and the parts of its remainder,
so that each term carries, at each switch, either a main part or one part of
a remainder.  If a term contains a part of some $r_k$ bounded by
$e^{-\frac2h(C+\eta)}$, then, since all costs are non-negative, its modulus
is $O(\exp(\frac1h(\Pot{\M}-\log\alpha-2C-2\eta)))$.  Otherwise its only
remainder parts are boundary parts: a part of $r_1$ bounded by
$e^{-\frac2h(A_0+\eta)}$, which occurs only with $w_0=w_1=D$, and the
analogous part of $r_m$.  Replacing $w_0$ (and $w_m$) by $R$ and removing
these factors leaves a word all of whose switches carry main parts.  This word
defines a history $H'$ whose ends lie in the sets $T_k$, and the term is
$O(e^{-2\eta/h})$ times $\exp(\frac1h(\Pot{\M}-\log\alpha-2c(H')))$.  Since
$c(H')\ge C$, every remainder term is
$O(e^{\frac1h(\Pot{\M}-\log\alpha-2C-2\eta)})$.

Hence $p_N$ is a finite exponential sum with non-vanishing coefficients, up
to an exponentially smaller error.  On each chart with fixed admissible ends, the non-resonance assumption makes the
real parts of distinct phases distinct off a set of zero area, so
Proposition~\ref{prop:expsum} applies, with $v=\Pot{\M}-2C$ almost everywhere
on $\Sweep$.  Outside $\Sweep$ we use Proposition~\ref{prop:exterior} as
before.  We apply Corollary~\ref{cor:charts} with the closed null set
consisting of $\mathcal E_0\cup Z\cup\partial\Sweep$ and the horizontal lines through
$V(0)$ and $V(1)$.  On a component $W$ inside $\Sweep$ we take as phases
the functions $\log\alpha+\Phi_H$ for all histories $H$ with admissible ends; they are
locally holomorphic by Lemma~\ref{lem:phases}.  On a component that is not
simply connected they may be multivalued, but their real parts are
single-valued, and this suffices for Corollary~\ref{cor:charts}, whose proof
is local.  The real parts
$\Pot{\M}-2c(H)$ differ pairwise by the non-resonance assumption, and their
maximum is $\Pot{\M}-2C=v$ almost everywhere.  This gives $\Pot{\mu}=F$ off
$\mathcal E_0\cup Z\cup\partial\Sweep$, and the description of the support; outside
$\Sweep$, $\Pot{\mu}=\Pot{\M}$ is harmonic.
\end{proof}

\begin{proposition}\label{prop:Swmonotone}
Under \textup{(H1)} and \textup{(H2)}, \textup{(Sw)} holds at every
$E\in\Sweep^\circ\setminus(\mathcal E\cup \mathcal E_0)$, with $\kappa_{1,\pm}=\pm\ii$ and
$\mathcal U_{1,\pm}=\mathcal U_\pm$.  Hence Theorem~\ref{thm:conditional}
recovers the formula of the Main Theorem; the non-resonance condition is
verified in the proof of the Main Theorem, and Lemma~\ref{lem:finite} gives
the finiteness of the set of arcs.
\end{proposition}

\begin{proof}
This is Proposition~\ref{prop:switch}, rewritten with the normalization of
Theorem~\ref{thm:recursion}.  In case (a$_\pm$) the end belongs to $T_1$.  In
case (b$_\pm$) with exit through $\Re z=0$, the proof of
Proposition~\ref{prop:switch} gives $|r_1|\le Ch\,e^{-\frac2h(A_L+c)}$ and
$U_\pm\ge A_L+c$; exit through $\Re z=1$ gives the same with $A_R$.  The
factor $Ch$ is absorbed by taking $\eta<c$ in \textup{(Sw)}.
\end{proof}

\subsection{Stokes regular points}
We now prove (Sw) in the case of several crossing points, under a geometric
condition on the real-momentum arc.  Assume that $V$ is holomorphic and univalent on a
convex domain $\Omega$, and let $E\in\Sweep^\circ\setminus \mathcal E_0$ with
$S\subset V(\Omega)$.  Then $\Gamma=\beta(S)$ is an arc from $x_+=\beta(w_+)$
to $x_-=\beta(w_-)$ through $x_1,\dots,x_m$, and the momentum is real exactly
on $\Gamma$.

\emph{End values and orientation.}  For a crossing point $x_k$ and
$t\in\{+,-\}$, the end value $p_{k,t}\in\pi\Z$ is the value at $w_t$ of the
continuation of $p$ along $S$ from $V(x_k)$; we have $p_{k,+}\in2\pi\Z$ and
$p_{k,-}\in\pi+2\pi\Z$.
The arc $\Gamma_{k,t}$ leaves $x_k$ into the half-plane
$\{\varepsilon_{k,t}\Im z>0\}$, and $\varepsilon_{k,-}=-\varepsilon_{k,+}$.  For
$\varepsilon=\pm1$, let $\pi_k^\varepsilon=p_{k,t}$, where $t$ is determined by
$\varepsilon_{k,t}=\varepsilon$.  The computation in the proof of
Lemma~\ref{lem:costpositive} gives $\varepsilon_{k,+}=\eta_k\varsigma_{k-1}$,
where $\eta_k=\mathrm{sgn}\sin p(x_k)$ and $\varsigma_{k-1}=s_k$.  Hence
\begin{equation}\label{eq:orient}
 \varsigma_{k-1}\bigl(\pi_k^--\pi_k^+\bigr)=\pi .
\end{equation}
If $0<p(x_k)<\pi$, then $\pi_k^\varepsilon=0$ for $\varepsilon=\varsigma_{k-1}$
and $\pi_k^\varepsilon=\pi$ for $\varepsilon=-\varsigma_{k-1}$.  In the monotone
case $\varsigma_0=1$, and the end values are $0$ above the real axis and $\pi$
below it.

\emph{Heights.}  The heights defined below refer to curves
$C_j=\{c_j(y)+\ii y:-Y^-\le y\le Y^+\}$, $0\le j\le m+1$, where $Y^\pm>0$;
we write $Y^\varepsilon$ for $Y^+$ or $Y^-$ according to $\varepsilon=\pm1$.
For $1\le k\le m$ and $\varepsilon=\pm1$, let $\ell_k^\varepsilon$ be the
horizontal segment at height $\varepsilon Y^\varepsilon$ between $C_{k-1}$ and
$C_{k+1}$, and let $z_k^\varepsilon$ be its endpoint towards which
$|e^{\frac\ii h\varsigma_{k-1}\int p}|$ grows.  Put
\begin{equation}\label{eq:heightaction}
 H_k^\varepsilon=\varsigma_{k-1}\,\Im\int_{x_k}^{z_k^\varepsilon}
 \bigl(p-\pi_k^\varepsilon\bigr)\,dz ,
\end{equation}
the integral being taken along $C_k$ up to height $\varepsilon Y^\varepsilon$
and then along $\ell_k^\varepsilon$, with $p$ continued from $x_k$.  The quantities
$H_k^\varepsilon$, $\varepsilon_{k,t}$ and the costs are unchanged if $p$ is
replaced by $-p$ (which reverses all $\varsigma_k$) or by $p+2\pi l$.

\begin{definition}\label{def:stokesregular}
The point $E$ is \emph{Stokes regular} if there are $Y^\pm>0$ and
$C^1$ functions $c_j$, $0\le j\le m+1$, on a neighbourhood of $[-Y^-,Y^+]$,
with $C_j\subset\Omega$, $c_0(0)=0$, $c_j(0)=x_j$ $(1\le j\le m)$,
$c_{m+1}(0)=1$ and $c_0<c_1<\dots<c_{m+1}$, such that:
\begin{enumerate}
\item[(S1)] every $C_j$ avoids $x_\pm$ and is canonical with a positive
margin for some branch $\pm p+2\pi l$ of the momentum, continued along $C_j$
from the real axis;
\item[(S2)] for $1\le j\le m$, $C_j$ meets $\Gamma$ only at $x_j$, and does so
transversally;
\item[(S3)] $\Gamma\subset\{-Y^-<\Im z<Y^+\}$, and $\Gamma$ lies strictly
between $C_0$ and $C_{m+1}$, with one exception: a subarc of $\Gamma$ that
joins $x_1$ (resp.\ $x_m$) to an end of $\Gamma$ without passing through
another crossing point may meet $C_0$ (resp.\ $C_{m+1}$), but only at a
single point and transversally;
\item[(S4)] $H_k^\varepsilon>C(E)$ whenever the part of $\Gamma$ that
leaves $x_k$ into $\{\varepsilon\Im z>0\}$ reaches another crossing point
before it reaches an end of $\Gamma$.  Moreover, $U_{1,t}\ge A_0$ (resp.\
$U_{m,t}\ge A_m$) if the switch through $x_{1,t}$ (resp.\ $x_{m,t}$) is
admissible and the subarc from $x_1$ (resp.\ $x_m$) to $x_{1,t}$ (resp.\
$x_{m,t}$) meets $C_0$ (resp.\ $C_{m+1}$).
\end{enumerate}
\end{definition}

We call the region $B$ bounded by $C_0$, $C_{m+1}$ and the lines
$\Im z=\pm Y^\pm$ the \emph{box}; it lies in $\Omega$ by convexity.  For $m=1$
we have $x_1=x_m$.  In this case Lemma~\ref{lem:cells} below, which uses
\textup{(S2)} for $C_1$, places the two end pieces at $x_1$ in $R_0$ and
$R_1$; the one in $R_0$ may exit through $C_0$, and the one in $R_1$ through
$C_2$.  Apart from this, the proof below does not use $C_1$, and Stokes
regularity is then a strengthened form of the situation of
Section~\ref{sec:proof}, with $C_0$ and $C_2$ in place of the vertical
segments through $0$ and $1$.  The extra condition $U\ge A_0$ in (S4) is
automatic in the monotone case: the proof of
Proposition~\ref{prop:switch}(b$_\pm$) gives $U_\pm\ge A_L+c$ for an exit
through $\Re z=0$ and $U_\pm\ge A_R+c$ for an exit through $\Re z=1$.  In
general we do not know whether it follows from the other conditions, because
$\Im z$ need not be monotone along an exiting end piece; it does hold in the
example of Theorem~\ref{thm:circle}.

\begin{lemma}\label{lem:cells}
Let $E$ be Stokes regular, and let $R_j$ be the region between $C_j$ and
$C_{j+1}$, $0\le j\le m$.  Then $\Gamma\setminus\{x_1,\dots,x_m\}$ consists of
the following pieces:
\begin{itemize}
\item an \emph{end piece} in $R_0$ from $x_1$ to one end of $\Gamma$;
\item for $1\le j\le m-1$, a \emph{bridge} in $R_j$ from $x_j$ to $x_{j+1}$,
which meets the real axis only at its endpoints;
\item an end piece in $R_m$ from $x_m$ to the other end of $\Gamma$.
\end{itemize}
An end piece may leave its cell only through $C_0$ or $C_{m+1}$, as allowed in
(S3).  Consecutive pieces leave $x_j$ into opposite half-planes, and the
part of each piece inside its cell lies in a single half-plane.  In
particular, the crossing points are met along $\Gamma$ in their natural
order, and only $x_1$ and $x_m$ can have admissible switches: $x_1$ through
the end of the piece in $R_0$, and $x_m$ through the end of the piece in
$R_m$.
\end{lemma}

\begin{proof}
Since $C_0$ and $C_{m+1}$ pass through $0$ and $1$, the real points of $\Gamma$
between them are $x_1,\dots,x_m$.  Removing these points splits $\Gamma$ into
$m+1$ open arcs.  By (S2) each arc lies in one cell, apart from the exits
allowed in (S3).  An arc joining $x_i$ to $x_j$ lies in a cell whose boundary
contains $C_i$ and $C_j$, so $|i-j|=1$.  By transversality the two arcs at
$x_j$ lie in $R_{j-1}$ and in $R_j$.  The only crossing curve on the boundary
of $R_0$ is $C_1$, so the arc in $R_0$ at $x_1$ runs to an end of $\Gamma$;
likewise for $R_m$.  If $m\ge2$, an arc in $R_j$ with $1\le j\le m-1$ cannot run
to an end of $\Gamma$, since both ends are already taken, so it joins $x_j$ to
$x_{j+1}$.  A bridge meets $\mathbb R$ only at its endpoints, so it lies in one
half-plane; inside the cells, an end piece meets $\mathbb R$ only at its
crossing point.  Since $\Gamma$ crosses $\mathbb R$ at each $x_j$, consecutive
pieces leave $x_j$ into opposite half-planes.  The statement on admissible
switches follows from Definition~\ref{def:admissible}.
\end{proof}

\begin{theorem}\label{thm:stokesregular}
Let $E$ be Stokes regular.  Then \textup{(Sw)} holds at $E$, with
\begin{equation}\label{eq:kappaU}
 \kappa_{k,t}=\ii\,\varepsilon_{k,t},\qquad
 \mathcal U_{k,t}=-\ii\varsigma_{k-1}\Bigl(\int_{x_k}^{x_{k,t}}p\,dz
 -p_{k,t}\,x_{k,t}\Bigr),
\end{equation}
the integral being taken along $\Gamma_{k,t}$; this agrees with
\eqref{eq:calU}, since $s_k=\varsigma_{k-1}$.  For $m\ge2$, $T_1$ consists of
the end of the piece in $R_0$ unless that piece meets $C_0$, and similarly for
$T_m$; moreover, $T_k=\emptyset$ for $1<k<m$.  For $m=1$, $T_1$ consists of the ends
whose pieces do not meet $C_0\cup C_2$.
\end{theorem}

The term $p_{k,t}x_{k,t}$ in \eqref{eq:kappaU} matters.  Since
$\Im(p_{k,t}x_{k,t})=\Im\int_{x_k}^{x_{k,t}}p_{k,t}\,dz$, it leaves the real
part $\Re\mathcal U_{k,t}=U_{k,t}$ unchanged.  But replacing it by $p_{k,t}(x_{k,t}-x_k)$, that
is, writing $\int(p-p_{k,t})\,dz$, multiplies $e^{-\frac2h\mathcal U_{k,t}}$
by $e^{2\ii\varsigma_{k-1}p_{k,t}x_k/h}$, and for generic $x_k$ this factor
differs from $1$ when $p_{k,t}\ne0$; see the numerical test in Section~\ref{sec:numgeneral}.
Under \textup{(H1)}, \eqref{eq:kappaU} reduces to the values of
$\mathcal U_\pm$ given after Theorem~\ref{thm:recursion}.

\begin{proof}
\emph{Solutions.}  Lemma~\ref{lem:segment} and its proof apply to a canonical
curve with a positive margin in place of a vertical segment: the lens is
formed by the curves $c_j(y)+sw_1\sin\frac{\pi(y+Y^-)}{Y^++Y^-}+\ii y$.  Let $u_j$ be
the WKB solution of type $\varsigma_j$ normalized at $c_j(0)$ along $C_j$,
$0\le j\le m+1$.  If $C_j$ is canonical for the branch $q=\eta p+2\pi l$
($\eta=\pm1$), we apply Theorem~\ref{thm:WKB} with $q$ and multiply the
solution by a constant times $e^{-2\pi\ii l\eta\varsigma_jz/h}$.  This turns its
WKB form into the one written with $p$, and the factor is constant on $h\Z$.
These solutions may be used in Theorem~\ref{thm:recursion}
because their representations hold near $x_j$; as the $c_j$ are defined on
a neighbourhood of $[-Y^-,Y^+]$, the representations hold up to the heights
$\pm Y^\pm$.  For each turning point inside $B$ we choose a cut joining it to
the side of $B$ that lies in the same half-plane as the corresponding end
piece (the top side for an end piece in $\{\Im z>0\}$, the bottom side
otherwise).  The cut is chosen inside the half-cell that contains that end
piece, and disjoint from the other pieces of $\Gamma$ and from the $C_j$;
this is possible because an end piece does not separate its half-cell.  Off
the cuts, $p$ is single-valued on $B$, and $\Im p$ vanishes only on $\Gamma$.

\emph{Fourier modes.}  Put
\[
 G_k(z)=W[u_{k+1},u_{k-1}](z)\,e^{-\frac{2\ii}h\varsigma_{k-1}\pi_k^+z}.
\]
The exponential factor is $h$-periodic because
$\varsigma_{k-1}\pi_k^+\in\pi\Z$, and it equals $1$ at every lattice point;
consequently $G_k=W[u_{k+1},u_{k-1}]$ on $h\Z$.  The two solutions
$u_{k-1}$ and $u_{k+1}$ have the same type $\varsigma_{k-1}$.  If $f=\hat fe^{\ii cz/h}$ and $g=\hat ge^{\ii cz/h}$,
then $W[f,g](z)=e^{\ii c}e^{2\ii cz/h}W[\hat f,\hat g](z)$.  We apply this with
$c=\varsigma_{k-1}\pi_k^+$ and $\hat u_{k\pm1}=e^{-\ii cz/h}u_{k\pm1}$, whose
WKB forms are $(\sin p)^{-1/2}e^{\frac\ii h\varsigma_{k-1}\int(p-\pi_k^+)}$ up to
constant factors.  Hence, on the top line, $G_k$ is $\pm1$ times the Wronskian
of these forms.  On the bottom line, the forms written with $p-\pi_k^-$
give, by \eqref{eq:orient},
$G_k=\pm e^{\frac{2\ii}h\varsigma_{k-1}(\pi_k^--\pi_k^+)z}(\cdots)
=\pm e^{2\pi\ii z/h}(\cdots)$.  Thus Lemma~\ref{lem:lattice} applies to $G_k$
exactly as it does to $\W$ in the proof of Proposition~\ref{prop:switch},
with items (i) and (ii) on the top line and items (iii) and (iv) on the
bottom line.  In the monotone case $\pi_k^+=0$ and $G_k=\W$.

\emph{Reduction to one orientation.}  Replacing $p$ by $-p$ or by $p+2\pi l$
only relabels the solutions.  Under these relabellings the right-hand sides of
\eqref{eq:kappaU} and the normalized coefficients $\tilde\alpha_k$ change by
the same factors, as one checks from the definitions; for example, $p\mapsto
p+2\pi l$ multiplies both $\tilde\alpha_k$ and $e^{-\frac2h\mathcal U_{k,t}}$ by
$e^{-4\pi\ii l\varsigma_{k-1}x_k/h}$.  The substitution $V\mapsto-V$,
$E\mapsto-E$ maps $p$ to $\pi-p$ and solutions $\psi$ to
$e^{\ii\pi z/h}\psi$.  It reverses every type and exchanges $w_+$ and $w_-$,
hence the labels $\pm$ of the ends; each arc $\Gamma_{k,t}$ keeps its
half-plane, and $\tilde\alpha_k$ and
$e^{-\frac2h\mathcal U_{k,t}}$ change by the same factor
$e^{2\pi\ii\varsigma_{k-1}x_k/h}$.  It therefore suffices to compute the main
terms in the case $0<p(x_k)<\pi$ and $\varsigma_{k-1}=1$, which is the
situation of Proposition~\ref{prop:switch}.

\emph{End crossings.}  Let $m\ge2$ and consider $x_1$; the case of $x_m$ is
symmetric.  Let $P_0$ be the end piece at $x_1$, with end $x_t=x_{1,t}$.
Suppose first that $P_0$ does not meet $C_0$.  By Lemma~\ref{lem:cells} it
lies in one half-plane, say $\{\Im z>0\}$, and the bridge in $R_1$ lies in
$\{\Im z<0\}$.  Let $\ell$ be the horizontal segment at a height $Y'<Y^+$ above
$\Gamma$ running from $C_0$ to $C_2$.  It meets no piece of $\Gamma$, so $\Im p\ne0$
for the momentum continued along it, and the sign of $\Im p$ is constant on
$\ell$.  Along $C_0$, which does not meet $\Gamma$, this sign is the one on
$I_0$.  Hence $u_0$ grows to the right along $\ell$, and
Theorem~\ref{thm:horizontal} continues it to $C_2$.  The closed curve formed
by $C_0$, $\ell$, $C_2$ and $[0,x_2]$ encloses $x_t$ and no other turning
point.  We now argue exactly as in case (a$_\pm$) of the proof of
Proposition~\ref{prop:switch}.  Near $\ell\cap C_2$ the continued
momentum is $2p_t-p$, with $p_t=p_{1,t}=\pi_1^+$, and the amplitude factor
$(\sin p)^{1/2}$ is multiplied by $\pm\ii$.  Lemma~\ref{lem:wronskian} gives $G_1$ on
the top line, and Lemma~\ref{lem:lattice}(i) gives the corresponding modes.
(If $P_0$ lies below the real axis, we use the bottom line and (iii).)  In the
normalization of Theorem~\ref{thm:recursion} the contribution to
$\tilde\alpha_1$ is $\kappa_{1,t}e^{-\frac2h\mathcal U_{1,t}}(1+O(h))$, with
$\kappa_{1,t}$ and $\mathcal U_{1,t}$ as in \eqref{eq:kappaU}.  In the
reference orientation this is the computation of
Proposition~\ref{prop:switch}, and comparison with the formulas after
Theorem~\ref{thm:recursion} gives $\kappa=\ii$ for the end above the real axis
($p_t=0$) and $\kappa=-\ii$ for the end below it ($p_t=\pi$).  The general case
follows from the reduction above.  Finally, since $\varsigma_0=s_1$ (both equal
$-\mathrm{sgn}\Im p$ on $I_0$), $\Re\mathcal U_{1,t}=U_{1,t}$.

\emph{Exits.}  Suppose now that $P_0$ meets $C_0$, at the single point $\zeta_*$
of (S3), and let $\{\varepsilon\Im z>0\}$ be the half-plane into which $P_0$
leaves $x_1$.  Then the end $x_t$ of $P_0$ lies beyond $C_0$, and the closed
curve formed by $C_0$, $\ell_1^\varepsilon$, $C_2$ and $[0,x_2]$ encloses no
turning point.  The half-plane $\{\varepsilon\Im z>0\}$ contains no piece of
$\Gamma$ in $R_1$ (Lemma~\ref{lem:cells}).  Hence on $\ell_1^\varepsilon$ the
sign of $\Im p$ is the one on $I_1$, and $u_2$ grows to the left there.  The
argument given below for sides without turning points therefore applies
verbatim and gives a contribution $O(he^{-\frac2hH_1^\varepsilon})$ to
$\tilde\alpha_1$, with $z_1^\varepsilon$ on $C_0$.  We deform the path in
\eqref{eq:heightaction} to $[x_1,0]$ followed by $C_0$, which is possible
because the region in between contains no turning point.  This gives
\[
 H_1^\varepsilon=A_0+c,\qquad
 c=\varsigma_0\Im\int_0^{c_0(\varepsilon Y^\varepsilon)+\ii\varepsilon
 Y^\varepsilon}\bigl(p-\pi_1^\varepsilon\bigr)\,dz\quad\text{along }C_0 .
\]
By the branch changes described above, we may replace $p$ by the branch
$\pm p+2\pi l$ for which $C_0$ is canonical.  The momentum at $\zeta_*$ is the
same whether it is continued along $P_0$ or along $C_0$; it is real, and it
lies in $(0,\pi)$ because $C_0$ is canonical for $p$.  Along $P_0$ the momentum moves
monotonically from $p(x_1)$ to $p_t=\pi_1^\varepsilon$ without passing a
multiple of $\pi$, so $0<p(x_1)<\pi$ and $p_t\in\{0,\pi\}$.  From
$\varepsilon_{1,+}=\varsigma_0$ and the canonical inequality along $C_0$ one
gets $c\ge\nu_0Y^\varepsilon>0$, where $\nu_0$ is the margin of $C_0$.  Thus
$r_1$ satisfies (ii) of \textup{(Sw)}.  If the switch is admissible, then
$U_{1,t}\ge A_0$ by (S4), so (iii) holds as well; if it is not admissible, it
does not occur in any history.

\emph{Sides without turning points.}  Let $1\le k\le m$, and let the part of
$\Gamma$ that leaves $x_k$ into $\{\varepsilon\Im z>0\}$ be a bridge.  By
Lemma~\ref{lem:cells}, the parts of $R_{k-1}$ and $R_k$ in this half-plane
contain at most bridges and no turning point, so the segment
$\ell_k^\varepsilon$ meets neither $\Gamma$ nor a cut.  On this segment,
$u_{k-1}$ and $u_{k+1}$ have the same type and the same branch.  Whichever of
them grows along $\ell_k^\varepsilon$ is continued by
Theorem~\ref{thm:horizontal} to the endpoint $z_k^\varepsilon$, where the
other is known.
Lemma~\ref{lem:wronskian} then gives
$|W[u_{k+1},u_{k-1}]|\le Ch|u_{k+1}u_{k-1}|$ near $z_k^\varepsilon$.  By
Lemma~\ref{lem:lattice}(ii) or (iv) applied to $G_k$, the normalizations of
Theorem~\ref{thm:recursion}, and Cauchy's theorem in the cells, the
corresponding modes contribute $O(he^{-\frac2hH_k^\varepsilon})$ to
$\tilde\alpha_k$, which is
$O(e^{-\frac2h(C+\eta)})$ by (S4).  For $1<k<m$ both half-planes are of this
kind.

\emph{The case $m=1$.}  The arguments for end crossings and exits apply to
both end pieces at $x_1$, with $C_2=C_{m+1}$; this is
Proposition~\ref{prop:switch} with $C_0,C_2$ in place of the vertical
segments.
\end{proof}

\begin{corollary}\label{cor:stokesregular}
Assume that \textup{(H2$'$)} holds, that almost every $E\in\Sweep$ is Stokes
regular (for some convex domain on which $V$ is univalent; the domain may
depend on $E$), and that the non-resonance condition of
Theorem~\ref{thm:conditional} is satisfied.
Then $\mu_N\to\mu$ with $\Pot{\mu}=F$ almost everywhere, that is,
$\Pot{\mu}=\Pot{\M}-2C$ almost everywhere on $\Sweep$ and $\Pot{\mu}=\Pot{\M}$
outside $\Sweep$.
\end{corollary}

\subsection{Vertical canonical curves and the circular arc}\label{sec:vertical}
We keep the assumptions of the previous subsection and put
\[
 L_\sigma=\{z\in\Omega:f(z)=\sigma\},\qquad
 \mathcal K=\{z\in L_\sigma:|r(z)-\tau|\ge2\alpha\}.
\]
Thus $L_\sigma=\Gamma\cup \mathcal K$ and $\Gamma\cap \mathcal K=\{x_+,x_-\}$.  Since
\[
 \mathfrak c(z)=\frac{E-V(z)}{2\alpha}=\frac{\tau-r(z)}{2\alpha}
 -\ii\,\frac{\sigma-f(z)}{2\alpha},
\]
$\mathfrak c$ takes values in $(-\infty,-1]\cup[1,\infty)$ exactly on $\mathcal K$.  Hence the
principal branch $p_0=\arccos\mathfrak c$ is holomorphic on $\Omega\setminus \mathcal K$, with
$0<\Re p_0<\pi$ there.  As in Lemma~\ref{lem:geometry}(iii), the sign of
$\Im p_0$ is that of $\sigma-f$.

\begin{lemma}\label{lem:vertical}
Every vertical segment in $\Omega\setminus \mathcal K$ is canonical for $p_0$ with a
positive margin, and Lemma~\ref{lem:segment} applies to it.
\end{lemma}

\begin{proof}
For a vertical segment the canonical condition reads
$\eta\le\Re p_0\le\pi-\eta$.  On the compact segment $\Re p_0$ is continuous
with values in $(0,\pi)$, so the condition holds with a positive margin.
Moreover, $\sin p_0$ vanishes only at the points of $\mathcal K$ where $\mathfrak c=\pm1$;
this does not use univalence.
\end{proof}

\begin{proposition}\label{prop:vertical}
Let $Y^\pm>0$ with $[0,1]+\ii[-Y^-,Y^+]\subset\Omega$, and put
$V_x=\{x+\ii y:-Y^-\le y\le Y^+\}$.  Suppose:
\begin{enumerate}
\item[(V1)] $V_{x_j}\cap L_\sigma=\{x_j\}$ and $\Re V'(x_j)\ne0$, for
$1\le j\le m$;
\item[(V2)] $V_0\cap \mathcal K=V_1\cap \mathcal K=\emptyset$;
\item[(V3)] $\Gamma$ satisfies (S3) with $C_0=V_0$ and $C_{m+1}=V_1$.
\end{enumerate}
Then (S1)--(S3) hold with $C_0=V_0$, $C_j=V_{x_j}$ and $C_{m+1}=V_1$.  Hence
$E$ is Stokes regular if, in addition, (S4) holds for these curves.
\end{proposition}

\begin{proof}
By (V1) and (V2) the segments avoid $\mathcal K$, so (S1) follows from
Lemma~\ref{lem:vertical}; along each segment $p_0$ is one of the branches
$\pm p+2\pi l$.  The tangent direction of $\Gamma$ at $x_j$ is $1/V'(x_j)$,
which is vertical exactly when $V'(x_j)$ is purely imaginary, so (V1) gives
(S2).  Condition (S3) is just (V3).
Since $\mathcal K$ and $L_\sigma$ are closed, (V1) and (V2) also hold for slightly longer
segments, as required in Definition~\ref{def:stokesregular}.
\end{proof}

\begin{corollary}\label{cor:Imb}
Suppose that $\Re V'\ne0$ on the convex domain $\Omega$.  Then:
\begin{itemize}
\item $V$ is univalent on $\Omega$, and each vertical line meets $L_\sigma$
at most once;
\item $\Gamma$ is the graph of a function over an interval of the real axis,
and the crossing points are met along $\Gamma$ in their natural order;
\item conditions (V1) and (V3) hold for every $Y^\pm$ such that $\Gamma$ lies
in $\{-Y^-<\Im z<Y^+\}$.
\end{itemize}
Hence $E$ is Stokes regular as soon as there exist such $Y^\pm$ for which
$[0,1]+\ii[-Y^-,Y^+]\subset\Omega$, $V_0\cap \mathcal K=V_1\cap \mathcal K=\emptyset$ and
(S4) holds.
\end{corollary}

\begin{proof}
The function $\Re V'$ has a constant sign on $\Omega$, so the segment-average
argument used after (H2), applied now to the real part of $V'$, proves
univalence.  Along a vertical segment
$\partial_y f=-\Re V'\ne0$, so $f$ is strictly monotone there.  The
tangent $1/V'$ of $L_\sigma$ has real part $\Re V'/|V'|^2\ne0$, so $L_\sigma$,
and with it $\Gamma$, is a graph over the real axis.  Such a graph meets each
of $V_0$ and $V_1$ at most once, and transversally.  The part of $\Gamma$ over
$\{x<0\}$ (resp.\ $\{x>1\}$) lies on the subarc from $x_1$ (resp.\ $x_m$) to
an end, which gives (V3).
\end{proof}

The condition $\Re V'\ne0$ makes $r$ strictly
monotone, in a strong sense, while $f'$ may vanish, so that the curve
$V([0,1])$ may have many horizontal tangents.  The remaining conditions to be
checked are \textup{(V2)}, the size of $\Omega$ and the heights in (S4).  The
circular arc of Theorem~\ref{thm:circle} below is a union of two pieces of
this type, and for it all these conditions can be verified.

\begin{theorem}\label{thm:circle}
Let $V(x)=-\ii e^{3\pi\ii x/2}$ and $\alpha>0$, and let $\mathcal E_c$ consist of $\mathcal E_0$ together with the
points with $\sigma>0$ and $\tau=2\alpha$.  Then every
$E\in\Sweep^\circ\setminus \mathcal E_c$ is Stokes regular, with vertical curves
$C_j$.  (The points with $\sigma=0$ lie in $\mathcal E_0$, since $f(1)=0$.)
\end{theorem}

\begin{proof}
\emph{Domains.} Write $B=\ii V=e^{3\pi\ii z/2}$. For
$-1/3<a<0$, $V$ is univalent on the strip
$\Omega_a=\{a<\Re z<a+4/3\}$. In the rotated energy coordinate
$w=\ii E$, the omitted ray of $B(\Omega_a)$ has angle $3\pi a/2$,
and the segment $\ii S$ lies on the line $\Re w=\sigma$. If $\sigma<0$, the ray
misses this line; if $\sigma>0$, the height of their intersection,
$\sigma\tan(3\pi a/2)$, tends to $-\infty$ as $a\downarrow-1/3$.
Thus $S\subset V(\Omega_a)$ for some $a$, and every position rectangle over
$[0,1]$ that we need is contained in $\Omega_a$.

\emph{The level set.}  Since $V'=\frac{3\pi\ii}2V$, we have
$\Re V'=\frac{3\pi}2f$ and $f(x+\ii y)=e^{-3\pi y/2}\cos\frac{3\pi
x}2$.  Hence $L_\sigma$ is the graph of
\[
 \gamma(x)=\frac2{3\pi}\log\frac{\cos(3\pi x/2)}\sigma
\]
over the set where $\cos(3\pi x/2)/\sigma>0$.  Each component of $L_\sigma$ is
a cap over an interval between consecutive zeros of $\cos(3\pi x/2)$,
intersected with $\Omega_a$.  On a full cap $\gamma$ attains its maximum
$y_\sigma=\frac2{3\pi}\log\frac1{|\sigma|}$ at the midpoint, is strictly
monotone on either side, and tends to $-\infty$ at the ends.  Along $L_\sigma$ we have $r=\sigma\tan(3\pi x/2)$, which is
strictly monotone in $x$ on each cap.  Hence $\Gamma$ is the graph of $\gamma$
over an interval $[x^L,x^R]$ inside one cap, and the crossing points are the
zeros of $\gamma$ in $[x^L,x^R]\cap(0,1)$.  Each vertical line meets $L_\sigma$
at most once, and $\Re V'(x_j)=\frac{3\pi}2\sigma\ne0$; hence (V1) holds.  By
Proposition~\ref{prop:vertical} it remains to check (V2), (V3) and (S4).

\emph{Case $\sigma<0$.}  The cap lies over $(\frac13,1)$, with its top at
$x=\frac23$.  The zeros of $\gamma$ are $x_a=\frac2{3\pi}\arccos\sigma$ and
$x_b=\frac43-x_a$.  The segments $V_0$ and $V_1$ do not meet $L_\sigma$, because
$\cos0/\sigma<0$ and $\cos(3\pi/2)=0$.  Taking $Y^\pm$ larger than the heights
of $\Gamma$, we obtain (V1)--(V3).  If $m=1$, condition (S4) is void.

Let $m=2$.  Then $\Gamma$ contains the graph over $[x_a,x_b]$, which is a
bridge in the upper half-plane, and the end pieces lie below the real axis.
Hence (S4) concerns the upper sides at $x_1=x_a$ and $x_2=x_b$.  The only
solutions of $f=\sigma$ in $[0,1]$ are $x_a$ and $x_b$, which are crossing
points, so $p=p_0$ on $[0,1]$, $\varsigma_0=1$ and $\varsigma_1=-1$.  From
$x_a$ the bridge runs to the right, where $r$ decreases; hence it heads
towards $w_+$, and $\pi_1^+=0$.  Similarly, $\pi_2^+=\pi$.  Let $Y^+>y_\sigma$; then the
line $\Im z=Y^+$ does not meet $L_\sigma$.  As $y\to\infty$,
$V(x+\ii y)\to0$ uniformly, so
\[
 p_0(x+\ii y)=p_\infty+O(e^{-3\pi y/2}),\qquad
 p_\infty=\arccos\frac{E}{2\alpha},
\]
and $0<\Re p_\infty<\pi$ because $E/(2\alpha)=(\tau-\ii\sigma)/(2\alpha)$
is not real.  The horizontal parts of the paths in \eqref{eq:heightaction}
have length less than $1$, and the integrands are bounded on them.  Hence
\[
 H_1^+\ge\Re p_\infty\,Y^+-M,\qquad H_2^+\ge(\pi-\Re p_\infty)\,Y^+-M,
\]
with $M$ independent of $Y^+$, and (S4) holds for large $Y^+$.

\emph{Case $\sigma>0$.}  The cap lies over $(-\frac13,\frac13)$ (intersected
with $\Omega_a$), with its top at $x=0$, and $f>\sigma$ exactly below it,
that is, for $y<\gamma(x)$.  The cap meets $(0,1)$ only at $x_a\in(0,\frac13)$,
so $m=1$, and (S4) concerns only exits.  The segment $V_1$ does not meet
$L_\sigma$, and $V_0$ meets it only at the top $\ii y_\sigma$, where $r=0$.
Since $r$ increases along the cap, the end $x^L$ of $\Gamma$ is $x_+$, and the
top belongs to $\Gamma$ exactly when $\tau<2\alpha$; note that
$\tau>\sqrt{1-\sigma^2}-2\alpha>-2\alpha$.

If $\tau>2\alpha$, then $x^L>0$.  We take $Y^+\in(\max\{\gamma(x^L),0\},y_\sigma)$, an interval that is
non-empty since $\sigma<1$, and $Y^-$ large; then (V1)--(V3) hold.

If $\tau<2\alpha$, we take $Y^\pm$ large.  Then (V2) holds, and the end piece
$P_0$, the graph over $[x^L,x_a]$, meets $V_0$ once, transversally, at the
top.  It remains to check that $U_{1,+}\ge A_0$ when the switch along $P_0$ is
admissible, that is, when $\gamma>0$ on $[x^L,x_a)$.  Here $p=p_0$ on
$[0,x_a]$, $\varsigma_0=1$ and $p_{1,+}=0$.  Hence
\[
 U_{1,+}-A_0=\Im\int_0^{x_+}p_0\,dz,
\]
where the integral is taken along $[0,x_a]$ followed by $P_0$.  Now replace this path by the
vertical segment from $0$ to $\ii\gamma(x^L)$ followed by the horizontal
segment to $x_+$.  The two paths enclose a region below the cap, and the new
path also lies below the cap, except for its endpoint $x_+$.  On this region
$f>\sigma$, so $p_0$ is holomorphic with $\Im p_0<0$, and therefore
\[
 U_{1,+}-A_0=\int_0^{\gamma(x^L)}\Re p_0(\ii y)\,dy
 -\int_{x^L}^0\Im p_0\bigl(x+\ii\gamma(x^L)\bigr)\,dx>0 .\qedhere
\]
\end{proof}

\begin{corollary}\label{cor:circle}
Let $V(x)=-\ii e^{3\pi\ii x/2}$ and $\alpha>0$, and assume that no two histories
have costs that coincide on an open set.  Then $\mu_N\to\mu$, where
$\Pot{\mu}=\Pot{\M}-2C$ almost everywhere on $\Sweep$ and $\Pot{\mu}=\Pot{\M}$
outside $\Sweep$.
\end{corollary}

\begin{proof}
The potential $V$ is entire with $V'\ne0$, and for $\sigma\ne0$ the segment $S$ avoids $0$, so
the inverse $\frac2{3\pi\ii}\log(\ii w)$ continues along $S$ and (H2$'$) holds
with $\Omega=\C$.  By Theorems~\ref{thm:circle} and~\ref{thm:stokesregular},
\textup{(Sw)} holds on $\Sweep^\circ\setminus \mathcal E_c$, and
Theorem~\ref{thm:conditional} applies.
\end{proof}

For other circular arcs $V(x)=B+Ae^{\ii\theta x}$,
$0<\theta<2\pi$, the level curves still admit a logarithmic-cosine
parametrization after translating and rescaling their arguments.  However,
$V_0$ or $V_1$ may then cross $\mathcal K$ below or above an end piece of
$\Gamma$, and tilted curves $C_0$, $C_{m+1}$ are needed there.  We do not
pursue this here.

\begin{conjecture}\label{conj:general}
For every real-analytic $V\colon[0,1]\to\C$ satisfying \textup{(H2$'$)},
\textup{(Sw)} holds for almost every $E\in\Sweep$, with
$\kappa_{k,t}$ and $\mathcal U_{k,t}$ as in \eqref{eq:kappaU}.  Consequently,
if no two histories have costs that coincide on an open set, then
$\mu_N\to\mu$, where
\[
 \Pot{\mu}=\Pot{\M}-2C\quad\text{a.e. on }\Sweep,\qquad
 \Pot{\mu}=\Pot{\M}\quad\text{outside }\Sweep.
\]
We further conjecture that the support is a locally finite union of analytic
arcs, including near the exceptional energies.
\end{conjecture}

Since $C\ge0$, the almost-everywhere formula gives
$\Pot{\mu}\le\Pot{\M}$ everywhere whenever the conditional limit exists; to
see this, apply the disk-mean inequality to $\Pot{\mu}$ and use the continuity
of $\Pot{\M}$ (Lemma~\ref{lem:PotM}). Exterior equality is already part of
Theorem~\ref{thm:conditional}, via Proposition~\ref{prop:exterior}. Neither
fact by itself shows that the whole support is a finite union of arcs; this
assertion is proved in the monotone case.

\emph{Remarks on the hypotheses.}  Stokes regularity is not automatic for
univalent $V$, and two obstructions are apparent.  First, by
Lemma~\ref{lem:cells}, Stokes regularity forces $\Gamma$ to visit the crossing
points in their natural order, with bridges alternating between the
half-planes.  This ordering is the local picture near a fold of $V([0,1])$,
and it holds whenever $\Re V'\ne0$ (Corollary~\ref{cor:Imb}).  For curves that
wind, it may fail when $m\ge3$; the switching coefficients at interior
crossing points then cannot be computed cell by cell.  Second, part of
$\mathcal K$ may pass over $0$ or $1$, below or above an end piece of
$\Gamma$; then the vertical segments are not canonical, and tilted curves are
needed.  Moreover, Theorem~\ref{thm:stokesregular} does not cover domains on
which the potential fails to be univalent.  For the circular arc, the
formula of Conjecture~\ref{conj:general} is proved under non-resonance
(Corollary~\ref{cor:circle}).  For $\Re V'\ne0$ the conjecture reduces to
\textup{(V2)}, the height conditions (S4), the size of $\Omega$ and
non-resonance; in general it remains open.

\begin{addendum}[Partial results on Conjecture~\ref{conj:general}]
Assume that \textup{(H2$'$)} holds, that \textup{(Sw)} holds for almost every
$E\in\Sweep$ with $\kappa_{k,t}$ and $\mathcal U_{k,t}$ as in
\eqref{eq:kappaU}, and that no two histories have costs that coincide on
an open set. Then $\mu_N\to\mu$ with $\Pot{\mu}=F$ almost everywhere
(Theorem~\ref{thm:conditional}), and the following hold.
\begin{enumerate}[(1)]
\item $\Pot{\mu}\ge\log\alpha$ everywhere on $\C$. In particular, $\mu$ has
no atoms, not even at the exceptional energies.
\item Every $E\in\Sweep^\circ\setminus\mathcal E_0$ has a neighbourhood in
which $\supp\mu$ is contained in a finite union of compact real-analytic
arcs. This includes the energies in $Z$, where the set of admissible ends
changes, so $Z$ can be omitted from the description of the support in
Theorem~\ref{thm:conditional}.
\item Let $a\in\C$ be arbitrary, and let $\Phi_1,\dots,\Phi_n$ be
holomorphic near $0$ and such that for almost every $s$ near $0$ we have
$F(a+s^2)=\Re\Phi_i(s)$ for some $i$. Then for some $\delta>0$ the set
$\supp\mu\cap\{|E-a|<\delta\}$ is contained in the image under
$s\mapsto a+s^2$ of finitely many compact real-analytic arcs. The same
holds, with the arcs themselves in place of their images, if
$\Phi_1,\dots,\Phi_n$ are holomorphic near $a$ and $F(E)=\Re\Phi_i(E)$ for
some $i$, for almost every $E$ near $a$.
\end{enumerate}
Independently of these assumptions, $\mathcal E_0\cup\partial\Sweep$ lies
on the union of finitely many horizontal lines and the two curves
$V([0,1])\pm2\alpha$. By (3), the last assertion of the conjecture at an
exceptional energy reduces to continuing the history phases there. Such a
continuation is possible at a simple end of a layer boundary, as follows. Let $E_0\in\Sweep^\circ$ be such that
$V(x)\in\{E_0-2\alpha,E_0+2\alpha\}$ for exactly one $x=x_0\in[0,1]$, with
$x_0\in(0,1)$. Assume that $f'\ne0$ at every solution of $f=\sigma$ in
$[0,1]$, that $f(0)\ne\sigma\ne f(1)$, that \textup{(H2$'$)} holds also at
$E_0$, and that the local inverse of $V$ at $V(x_0)$ continues
holomorphically along the frozen segment of $E_0$. Then every history phase
at a regular energy near $E_0$, on either side of the layer boundary, is
one of finitely many functions holomorphic on a fixed neighbourhood of
$E_0$. Here the action $\int_0^1p\,dx$ is taken along a polygon that passes the
turning point on the side away from it; the switch through the nearby
turning point uses the cubic germ of the action there, and the opposite
switch at the new crossing uses the continued inverse. No square-root
coordinate is needed. Hence, near such $E_0$, $\supp\mu$ is contained in a
finite union of compact real-analytic arcs. The same conclusion holds near
every $E_0\in\Sweep^\circ$ with no turning point on $[0,1]$ for which $f'\ne0$ at every
solution of $f=\sigma$ in $[0,1]$ and $f(0)\ne\sigma$, provided that \textup{(H2$'$)}
holds at $E_0$ and that, when $x=1$ is a crossing point of $E_0$, the inverse of
$V$ continues along both segments from $V(1)$. This case includes the height
$f(1)=\sigma$, where a crossing point enters or leaves through $x=1$; all
switching primitives are then anchored at $x=0$. The remaining exceptional
energies (critical heights of $f$, the height $f(0)$, non-simple contacts,
and $\partial\Sweep$), as well as \textup{(Sw)} itself, remain open.

Stokes regularity is not implied by \textup{(H2$'$)}. For
$V(x)=e^{4\pi\ii x}$ and $\alpha=1$, condition \textup{(H2$'$)} holds with
$\Omega=\C$ at every energy of an open subset of $\Sweep^\circ$ of positive
area, but no energy is Stokes regular, since $V(0)=V(1)$. Hence
Theorem~\ref{thm:stokesregular} cannot give \textup{(Sw)} under
\textup{(H2$'$)} alone.
\end{addendum}

\subsection{Numerical evidence for several crossings}\label{sec:numgeneral}
For the circular arc we take $\alpha=0.4$. Where $f$ has two preimages inside
the same horizontal frozen segment, the real-momentum arc has an upper bridge
between them and two lower end pieces. Figure~\ref{fig:circle} compares
the computed spectrum with the prediction from the actions. A more sensitive
test uses the (complex) Dirichlet determinant itself. At
$E=-0.4+0.9\ii$ and at $0.2-0.5\ii$ there is a single crossing, with
opposite orientations in the two cases, and at $0.05+0.7\ii$ there are two crossings.  The
second energy is an exit case, in which the end piece leaves through $V_0$;
it therefore also tests the exit part of Theorem~\ref{thm:stokesregular}.
Substituting the switching coefficients given by
\eqref{eq:kappaU} into Theorem~\ref{thm:recursion} gives a relative error
below $5\times10^{-4}$ in the determinant at $N=1600$ in these examples.
Replacing $p_{k,t}x_{k,t}$ by $p_{k,t}(x_{k,t}-x_k)$, which changes only the
phase, alters the result already at leading order. The comparison thus tests
the phase as well as the logarithmic potential, but it does not establish
non-resonance on an open set.

The spiral $V(x)=-\ii\exp((0.5+1.9\pi\ii)x)$ with $\alpha=0.6$
illustrates a further issue: up to three frozen segments can overlap.
For instance, at $E=0.199-0.836\ii$, numerical continuation of a
real-momentum arc returns to the real position axis at $x<0$. Assigning
to that arc the cost of an admissible switch would give a negative value,
but Definition~\ref{def:admissible} excludes it. This example lies outside the
scope of the circular-arc corollary.

\begin{figure}[tbp]
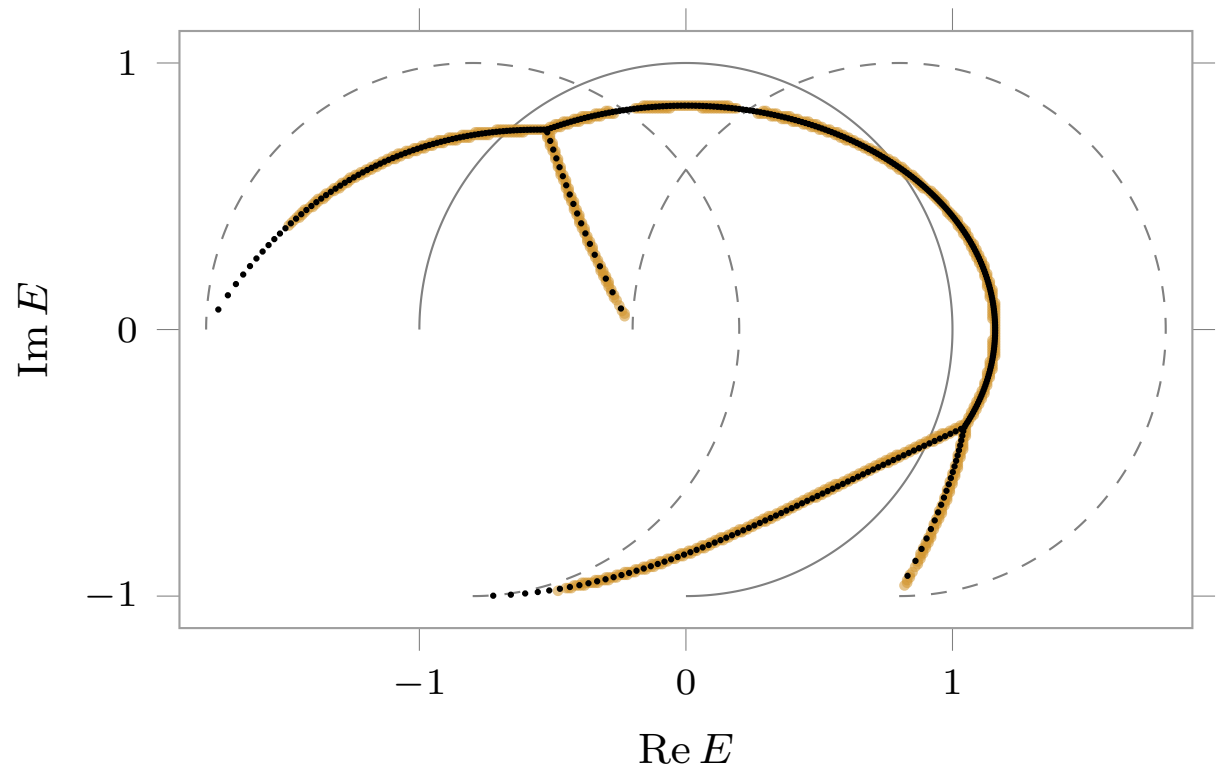

\centering

\caption{Circular potential $V(x)=-\ii e^{3\pi\ii x/2}$, $\alpha=0.4$, $N=400$.
Black dots are eigenvalues. Gold points approximate the set where the
proposed potential $F$ is not harmonic, detected by its discrete Laplacian on a
$301\times301$ grid. Grid points within $0.03$ of the dashed curves
$V(x)\pm2\alpha$ were excluded to avoid boundary artefacts. The solid
grey curve is $V([0,1])$.}\label{fig:circle}
\end{figure}

\section{Open questions}\label{sec:outlook}
For a general analytic potential, the unresolved step is to compute the
switching coefficients when the real-momentum arcs cannot be separated
by canonical curves. Even for the circular arc, it remains to exclude
persistent ties between the costs of distinct histories. The conditional
result localizes the support away from the exceptional energies, but for
general potentials no global finiteness assertion follows from this local
argument.

The estimates here are pointwise in energy. At an endpoint of a frozen
segment, a turning point can reach the real position interval; at a corner,
it can reach a Dirichlet boundary. Uniform asymptotics across such energies
would require a separate turning-point analysis; see \cite{FK19,Fed21,Fed22,KO26}
for related settings.

For a variable adjacent coefficient $c(x)$, the frozen equation becomes
$2c(x)\cos p=E-V(x)$, and the amplitude and connection calculations would
have to be redone. A related question is whether junctions between
analytic pieces permit switches of only algebraic cost, and whether this
changes the limit for successive piecewise-linear approximations of an
analytic potential. The Charlier case \cite{BlStC} offers a concrete
comparison between saddle-point phases and transfer-matrix actions when
the product of the adjacent coefficients varies.

\section*{Acknowledgments}
The author thanks A.~Kuijlaars, M.~Duits, the late S.~Naboko,
F.~\v{S}tampach and, especially, D.~Giandinoto for discussions that
helped to shape this paper over almost a decade.

\section*{Statement on the use of AI}
OpenAI's Codex and Anthropic's Claude were used in the substantial
revision and strengthening of an earlier manuscript, in checking
references, in implementing the numerical calculations, in formalizing
proofs in the Lean proof assistant, and in language editing and
typesetting checks. Responsibility for the mathematical content of the
paper rests with the author.

\end{document}